\documentclass[12pt,a4paper,final]{article}
\usepackage[margin=24.9mm]{geometry}

\usepackage[utf8]{inputenc}
\usepackage{amsmath,amsfonts,amssymb,bef_alex}

\usepackage{todonotes}

\numberwithin{equation}{section}

\usepackage[color,notref,notcite]{showkeys}

\def\bm{\mathbf}

\numberwithin{figure}{section}

\newcommand{\STEP}[1]{\vspace{0.2em}\noindent \emph{#1}}
\def\ti{{\times}}
\newcommand{\Div}{\mathop{\mathrm{div}}}
\newcommand{\weak}{\rightharpoonup}
\newcommand{\bbin}{\text{\footnotesize\raisebox{-0.1em}{$\;\overset{\,\mafo{bdd}\;}{\in}\,$}}}

\newcommand{\st}{\mathrm{st}}
\newcommand{\mfP}{\mathfrak P}
\newcommand{\INIv}{v^0}
\newcommand{\INIw}{w^0}

\begin{document}
\title{On two coupled degenerate parabolic equations\\ motivated by 
 thermodynamics\thanks{The research was partially supported by Deutsche
   Forschungsgemeinschaft (DFG) via the Collaborative Research Center SFB\,910
   ``Control of self-organizing nonlinear systems'' (project no.\ 163436311),
   subproject A5 ``Pattern formation in coupled parabolic systems''.}}
\author{Alexander Mielke\thanks{Weierstra\ss-Institut f\"ur Angewandte Analysis
    und Stochastik, Mohrenstra\ss{}e 39, D--10117 Berlin,
    \texttt{alexander.mielke@wias-berlin.de} 
and
Institut f\"ur Mathematik, Humboldt-Universit\"at zu Berlin,
Rudower Chaussee 25, D--12489 Berlin (Adlershof), Germany}}
\date{13. Mai 2022}
\maketitle
\vspace*{-2em} 
\begin{abstract} 
  We discuss a system of two coupled parabolic equations that have degenerate
  diffusion constants depending on the energy-like variable. The dissipation of
  the velocity-like variable is fed as a source term into the energy equation
  leading to conservation of the total energy. The motivation of studying this
  system comes from Prandtl's and Kolmogorov's one and two-equation models for
  turbulence, where the energy-like variable is the mean turbulent kinetic
  energy.

  Because of the degeneracies there are solutions with time-dependent support
  like in the porous medium equation, which is contained in our system as a
  special case.  The motion of the free boundary may be driven by either
  self-diffusion of the energy-like variable or by dissipation of the
  velocity-like variable. The cross-over of these two phenomena is exemplified
  for the associated planar traveling fronts. We provide existence of suitably
  defined weak and very weak solutions. After providing a thermodynamically
  motivated gradient structure we also establish convergence into steady state
  for bounded domains and provide a conjecture on the asymptotically
  self-similar behavior of the solutions in $\R^d$ for large times.
\end{abstract}

% Keywords: free boundary, growing support, gradient system, energy
% conservation, momentum conservation, porous medium equation,
% energy-dissipation estimates, entropy estimates 

% MSC: 35K65 % Degenerate parabolic equations
%      35K40 % Second-order parabolic systems
%      80M30 % Variational methods applied to problems in thermodynamics and heat transfer
%      49S05 % Variational principles of physics 

% WIAS 
% Analysis of Partial Differential Equations and Evolutionary Equations
% Variational methods
% Systems of partial differential equations: modeling, numerical analysis and simulation

\vspace{-2em}

\section{Introduction}
\label{se:SimpleModel}

On a smooth domain $\Omega \subset \R^d $ we consider the degenerate parabolic system 
\begin{subequations}
  \label{eq:SM01}
\begin{align}
  \label{eq:SM01.a}
& \dot v = \DIV\big( \eta(w) \nabla v\big), &\text{for }&(t,x)\in {]0,\infty[} \ti \Omega,
\\
  \label{eq:SM01.b}
&\dot w = \DIV \big( \kappa(w) \nabla
  w\big) + \eta(w) |\nabla v|^2 &\text{for }&(t,x)\in {]0,\infty[} \ti \Omega,
\\[0.3em]
  \label{eq:SM01.c}
& 0=\eta(w) \nabla u\cdot \rmn , \quad 0= \kappa(w) \nabla w \cdot \rmn&
\text{for } &(t,x)\in {]0,\infty[} \ti \pl\Omega ,
\end{align}
\end{subequations}
where $v(t,x)\in \R$ can be considered as a shear velocity and $w(t,x)\geq 0$
is an internal energy. Here the functions $w\mapsto \eta(w)$ and $w\mapsto
\kappa(w)$ describe the viscosity law for $v$ and the energy-transport coefficient
for $w$. Throughout this work we will mainly restrict to the choice
\begin{equation}
  \label{eq:Specif}
  \eta(w)=\eta_0 w^\alpha  \quad \text{and} \quad \kappa(w)= \kappa_0 w^\beta,
\end{equation}
where $\alpha,\beta,\eta_0,\kappa_0>0$ are given parameters.

The main feature of the model is that the shearing dissipation
$\eta(w)|\nabla v|^2$ is feeding into the energy equation such that in addition
to the total momentum $\calV(v,w)$ also the total
energy $ \calE(v,w)$ are conserved along solutions:
\[
\calV(v,w):=\int_\Omega v(x) \dd x \quad \text{and} \quad 
\calE(v,w):=\int_\Omega \big(\frac12v^2 + w\big) \dd x 
\]

One difficulty of the coupled system is that the viscosity coefficient
$\eta(w)$ and the energy-transport coefficient $\kappa(w)$ can be unbounded, but
this problem will play a minor role in our work. The main emphasis is on the 
degeneracies arising from the fact that 
$\eta(0)=\kappa(0)=0$ and that the solutions of our interest have a nontrivial
support. Thus, we are deriving a theory for solutions that have
$(v(t,x),w(t,x)=(0,0))$ in regions of the $Q_T=[0,T]\ti \Omega$ of full
measure. In particular, we are interested in the free boundary arising at the
boundary of the time-dependent support of $w(t,\cdot)$. There is already some
existence theory for related models motivated by turbulence in fluids, see 
\cite{GLLMT03TSUE,LedLew07RANS,DruNau09EWSS} for stationary models and  
\cite{Naum13DPPT, MieNau15GTEW, BulMal19LDAK, MieNau18?EGTW} for time-dependent
models. However, there it is either assumed that $\INIw(x)\geq \underline w>0$ for
all $x\in \Omega$ or that $\INIw(x)=0$ for all $x \in \pl\Omega$ and
$\INIw(x) >0$ a.e.\ in $\Omega$ such that $\int_\Omega \log \INIw(x) \dd x >
- \infty$, see Remark \ref{re:BuliMalek}. 

We provide a preliminary existence theory for our coupled system in Section
\ref{se:ExistTheory} which allows for solution with nontrivial support, $\Omega
\setminus\mafo{sppt}(w(t))$ has nonempty interior. However, we emphasize that this is
just for completeness and we rather focus on the growth behavior of the
support, i.e.\ the moving free boundary. Moreover, we emphasize that the
present paper is not a typical paper in applied analysis, but rather a paper in
modeling. Many of the 
statements in this paper are in fact conjectures, and only a few results
are formulated rigorously as propositions or theorems. Nevertheless, we believe
that the degenerate coupled system is relevant in applications and open up new
avenues for developing the tools in applied analysis, in particular in the
field of free-boundary problems. The system specific enough to analyze it in
more detail, it is close enough to the porous medium
equation (PME) to lend some of the tools
from there, but it displays a richer structure of nontrivial effects stemming
from the coupling between the two scalar equations.

To start with we remark that \eqref{eq:SM01} contains the PME, when restricting
to the case $v\equiv 0$: 
\begin{equation}
  \label{eq:I.PME}
  \dot w = \DIV \big( \kappa(w) \nabla   w\big) \quad \text{in }{]0,\infty[} \ti
  \Omega,
\qquad \kappa(w) \nabla w \cdot \rmn=0 \quad \text{on } {]0,\infty[} \ti \pl\Omega,
\end{equation}
which is known for its solutions with time-dependent support. For
$\Omega=\R^d$ and $\kappa(w)=(\beta{+}1)w^\beta$ we have the celebrated 
self-similar Barenblatt solutions (see \cite[Eqn.\,(1.8)]{Vazq07PMEM}) 
\begin{equation}
  \label{eq:PME.SimiSol}
  w(t,x)= \frac1{(t{+}t_*)^{d\delta}} \Big(\!\max\Big\{ C- \frac{ k|x|^2}{(t{+}t_*)^{2\delta}}
\,,\,0\, \Big\}\Big)^{1/\beta}  \text{ with } \delta=\frac1{2{+}d\beta},\ 
 k=\frac{\delta \beta}{2(\beta{+}1)}.
\end{equation}
Here $C>0$ determines the conserved total mass $\int_{\R^d}
w(t,x)\dd x = E_0$, see Section \ref{su:SimilPME}.

For $v\neq 0$ there is a true coupling between the two scalar equations, and
its structure is discussed in Section \ref{se:ThermoDyn}. In addition to the
conservation laws for momentum and energy and symmetries, we show that all
entropies of the form $\calS(v,w)=\int_\Omega \sigma(w(x)) \dd x$ with
nondecreasing and concave $\sigma$ are growing along solutions $(v,w)$ of
\eqref{eq:SM01}. Moreover, we establish a gradient
structure the the coupled system: For given $\sigma$ with $\sigma'(w) > 0$ and
$\sigma''(w)< 0$ there exists a state-dependent Onsager operator
$\bbK=\bbK^*\geq 0$ describing the 
dissipation mechanisms:
\begin{align}
&  \label{eq:I.GradFlow}
  \binom{\dot v}{\dot w} = \bbK(v,w)\rmD\calS(v,w) =
  \bbK(v,w)\binom{0}{\sigma'(w)}  \quad \text{with }
\\[0.3em] \nonumber
&  \bbK(v,w)\binom{\zeta}{\xi} = 
  \binom{-\DIV\big(p_1(w)(\nabla \zeta-\xi\nabla v)\big)} 
        {-p_1(w)\nabla v\cdot \nabla \zeta +p_1(w)|\nabla v|^2 \xi -
          \DIV\big(p_2(w)\nabla \xi\big) }  
\end{align}
for suitably chosen functions $w\mapsto p_j(w)$, see Section
\ref{su:Thermodyn}.

Sections \ref{su:Scalings} and \ref{su:ExaSimiSol} 
are devoted to scaling invariances and self-similar
solutions of the coupled system \eqref{eq:SM01}. We argue that even for
general $\eta$ and $\kappa$ there are solutions of the form 
$(v(t,x),w(t,x))=\big( V(x/(t{+}1)^{1/2}), W((x/(t{+}1)^{1/2})\big)$ where
$(V,W)$ may attain nontrivial limits for $y\to \pm \infty$. For
$\eta(w)=\kappa(w)=w$ an explicit family of solutions with nontrivial support
of $W$ is provided in Example \ref{ex:SimiExplicit}.

In Section \ref{su:Thermodyn} we consider the case of bounded $\Omega$ and
exploit the gradient structure for the case $\alpha=\beta=1/2$ to show that
most solutions converge exponentially to constants states
$(v(t),w(t)) \to (\wh v_{V_0,E_0} \bm1_\Omega, \wh w_{V_0,E_0}\bm1_\Omega)$,
where $\wh v_{V_0,E_0}$ and $\wh w_{V_0,E_0}$ are given explicitly in terms of
$V_0=\calV(\INIv)$ and $E_0=\calE(\INIv,\INIw)$. The exponential decay rate is
quite explicit. However, we also show that the decay does not hold for all
solutions: for instance, because of non-uniqueness we may have
$v(t,x)=\INIv(x)$ while $w\equiv 0$, which is certainly not decaying to the
thermodynamic equilibrium.

We also compare our model to the plasma model discussed \cite{RosHym85ANS,
  RosHym86PDAM, HymRos86ANPE} for the mass density $\rho\geq 0$ and the
temperature $\theta\geq 0$:  
\begin{equation}
  \label{eq:Rosenau}
  \rho_t = \DIV\big( \rho^\gamma \phi_1(\rho,\theta)\nabla \rho \big) \quad \text{and}
\quad  \big(\rho\theta\big)_t= \DIV\big(\rho^\delta
\phi_2(\rho,\theta)\nabla \theta + \theta \rho^\gamma
\phi_1(\rho,\theta)\nabla\rho \big), 
\end{equation}
see Section \ref{su:Rosenau} for more details. \medskip
 
% Section 3:         se:SM.TravFronts 
%3.1 Steady states  su:SteadyStates
%3.2 Traveling fronts su:GenPlanarW
% 3.3 Conjectured behavior near boundary of growing supports su:Conj.GrowSppt
  
Section \ref{se:SM.TravFronts} is devoted to steady states and traveling
fronts. Because of the degeneracy it is obvious that all function of the form
$(v,w)=(\INIv,0)$ are steady states, which we call trivial steady
states. Nontrivial steady states are necessarily spatially constant, i.e.\
$(v,w)=(v_*,w_*)=$ const., which provides, for bounded domains, a unique steady
state $(\wh v_{V_0,E_0},\wh w_{V_0,E_0})$ as introduced above. 

In Section \ref{su:GenPlanarW} we study planar traveling fronts of the form
\[
\big(v(t,x),w(t,x)\big) = \big(V(z), W(z)\big) \quad
\text{with } z= x_1{+}c_\rmF t, 
\]
where $ c_\rmF \in \R$
is the front speed. It is well-known that the planar fronts play an important
role in the theory of the PME (cf.\ \cite[Sec.\,4.3]{Vazq07PMEM}), and we
expect a similar role for our coupled system \eqref{eq:SM01}, in particular,
for the understanding of the propagation of the boundary of the support. 
Inserting this ansatz into \eqref{eq:SM01} and assuming $V(z)=W(z)=0$ for
$z\leq 0$, which simulates a 
support propagating with front speed $c_\rmF$, we obtain after integrating each
equations once (see Section \ref{se:SM.TravFronts} for details) the two ODEs
\[
c_\rmF\,V=\eta(W) V', \quad c_\rmF\, \big(W-\frac12 V^2) = \kappa(W) W'.
\]  
We analyze all solutions of this system, for
the different cases occurring for the choices in \eqref{eq:Specif}. To highlight
one of the results, we consider the case $\eta(w)=w$ and $\kappa(w)=\kappa_0
w$, i.e.\ $\alpha=\beta=1$. For 
$\kappa_0 \geq 1/2$ all traveling fronts have the form
  $(V(z),W(z))=(0,\frac1{\kappa_0}\,c_\rmF z) $ for $z\geq 0$, which
  corresponds to the case of the pure PME with $v\equiv0$. These solutions
  still exists for $\kappa_0 \in {[0,1/2[}$, but now additional, truly coupled
  solutions exists: 
\[
 (V(z),W(z))= \Big( 2\sqrt{(1{-}2\kappa_0)\,c_\rmF z\,}\:  ,\:
 2c_\rmF z \Big) \quad  \text{for } z\geq 0. 
\]
For these solutions, the propagation of the support of $w$ is not only driven
by self-diffusion as  for the PME, but is is driven also by the
generation of $w$ via the source term $\eta(w)|\nabla u|^2$. This is best seen
in the limit $\kappa_0\to 0$, where self-diffusion disappears but
propagation is still possible. In particular we obtain $c_\rmF= \max\{
\kappa_0, 1/2 \} W'(0)$, which again shows that for $\kappa_0<1/2 $ the
propagation speed is no longer dominated by self-diffusion alone.  

In Section \ref{su:Conj.GrowSppt} we conjecture that the typical behavior of
$w$ near the boundary of its support is given by $w(t,x)= w_0 (z_+)^\gamma$
with $\gamma = \max\{ 1/\alpha, 1/\beta\}$, which clearly shows that the front
is driven by the $v$-diffusion in case of $\beta >\alpha$. In the critical case
$\alpha=\beta$ the switch between the two regimes occurs for
$\eta_0=2\kappa_0$.\medskip

Our definitions of weak and very weak solutions are given in Section
\ref{se:WeakVeryWeakSol} and are based on a reformulation of the coupled system
\eqref{eq:SM01} in terms of \eqref{eq:SM01.a} and a conservation law
\eqref{eq:ConsLaw.e} for the energy density $e=\frac12+w$, thus following the
ideas in \cite{FeiMal06NSET, BuFeMa09NSFS}. This allows us to avoid defect
measures. The notion of very weak solutions is based on the 
\emph{weak weighted gradient} $G_{w^\alpha} v$, where $w^\alpha \nabla v$ is
defined in terms of the distributional form of $ \nabla (w^\alpha v ) - v
\nabla (w^\alpha) $, thus 
avoiding any derivatives of $v$ but using $\nabla w$ instead, see Definition
\ref{de:WeightedGrad}, where $G_a v \in \rmL^1(\Omega)$ is defined for $a \in
\rmW^{1,q}(\Omega)$ and $v \in \rmL^{q^*}(\Omega)$.  Section
\ref{su:NonuniqueVWS} provides an explicit example for nonuniqueness of very
weak solutions. 

Before showing existence of solution, we provide a series of natural a
priori bounds for strong solutions satisfying $w>0$ in
$\ol Q_T=[0,T]\ti \ol\Omega$. Section \ref{suu:IntegralEstim} establishes
$\rmL^p$ bounds for $v$ and, in the case $\alpha=\beta$, also
for $w$. Section \ref{su:ComparisonEst} provides comparison. Here, we also show
that the case $\eta=\kappa$ is very special, as for this case an estimate
$|\INIv(x)| \leq M_*\INIw(x)$ for all $x\in \Omega$ propagates for positive time,
i.e.\ we have $|v(t,x)| \leq M_x w(t,x)$ in all of $Q_T$. The crucial
dissipation estimates are discussed in Section \ref{su:DissipEstim}. In
particular, for $\alpha\in {]0,1[}$ and bounded $\Omega$ we obtain
$\int_{Q_T} |\nabla v|^2 \dd x \dd t \leq
C(\Omega,\alpha,\INIv,\INIw)$.

In Section \ref{se:ExistTheory} we develop our (rather preliminary) existence
theory for bounded $\Omega$. For this we approximate the initial data be smooth
functions $(\INIv_\eps,\INIw_\eps)$ which additionally satisfy
$\INIw_\eps(x)\geq \eps$. Using the comparison principles derived earlier we
find classical solutions $(v_\eps,w_\eps):Q_T\to \R^2$ still satisfying
$w_\eps(t,x)\geq \eps$. For passing to the limit $\eps\to 0^+$ we use the
appropriate a priori bounds proving spatial and temporal compactness such that
a suitable version of the Aubin-Lions-Simon lemma provides strong convergence. 
So far, we are only able to treat the case $\alpha \in {]0,1[}$, where one can
exploit the $\rmL^2$ a priori bound for $\nabla v_\eps$, which was also used in
\cite{Naum13DPPT}. This approach allows us to construct weak solutions.  
For $\alpha=1$ we only able to handle the case $\eta\equiv \kappa$ and we are
only able to establish very weak solutions, where the weak weighted gradient
$G_{w^\alpha} v$ is well-defined but $\nabla v$ may be not. From the
expected behavior of the solutions near the boundary of the support, it is
clear that gradients may have a blow-up, so that the difference between weak
and very weak solutions may be essential. 

Hence, our existence theory is different from the one developed in
\cite{BerKam90SDPE,DalGia99WSSC} for the plasma model \eqref{eq:Rosenau},
because we enforce global Sobolev regularity, while the latter as for local
regularity of $\theta$ on the support of $\rho$ only.\medskip

Section \ref{se:Scaling} provides a few conjectures concerning
the longtime behavior in the case $\Omega=\R^d$ and $\eta(w)=\eta_0 w^\alpha$
and $\kappa(w)=\kappa_0 w^\beta$.  The whole system does not have any self-similar
solution, however we expect that in many cases $v$ and $w$ behave self-similar
in the limit $t\to \infty$. In these cases we expect that $v(t)$ converges to
$0$ in $\rmL^2(\Omega)$ while the momentum is conserved
$\calV(v(t))=\int_\Omega v(t,x)\dd x =\calV(\INIv)$. Then
$\int_\Omega w(t,x)\to \calE(\INIv,\INIw)$ and we expect that $w$ behaves like
the solution of the PME obtained from \eqref{eq:SM01} for $v\equiv 0$, but now
the total energy is fixed to $\calE(\INIv,\INIw)$.

Finally, in Section \ref{su:DerivFluids} we show how our coupled model
\eqref{eq:SM01} is motivated by models from turbulence modeling, where $w$
plays the role of the mean turbulent kinetic energy, such that $e=\frac12v^2+w$
denotes the total kinetic energy. Our model is obtained when the solutions
$\bfu$ of the Navier-Stokes equation are assumed to be parallel flows, namely
$\bfu(t,x)= (0,...,0,v(t,x_1,...,x_d)^\top \in \R^{d+1}$ with $d\in \{1,2\}$.
Prandtl's model for turbulence (cf.\ \cite{Pran46UNFA, Naum13DPPT}) is
discussed in Section \ref{suu:Kolmo1eqn} relating to our case
$\alpha=\beta=1/2$, while Kolmogorov's two-equation model (\cite{Kolm42ETMI,
  Spal91KTEM}) is discussed in Section \ref{suu:Model} relating to our case
$\alpha=\beta=1$. For the rich theory of these models we refer to
\cite{Lewa97MACT, Naum13DPPT, ChaLew14MNFT, MieNau15GTEW, BulMal19LDAK,
  MieNau18?EGTW} and the references therein.\vspace{-0.5em}

\section{The model and its thermodynamical formulation}
\label{se:ThermoDyn}\vspace{-0.5em}

Here we discuss the basic properties of system \eqref{eq:SM01},
namely  the conservation laws for total linear momentum and the total
energy, the symmetries and scalings, as well as exact similarity solutions. 
Section \eqref{su:Thermodyn} provides gradient structures which
allows us to show convergence into steady state for the case $\alpha=\beta=12$
and $\Omega$ bounded, see Theorem \ref{th:Cvg.SteadySt}.\vspace{-0.5em}

\subsection{Conservation laws}
\label{su:ConservLaws} 

We first observe that the divergence structure of the equation for $v$ and the
no-flux boundary condition provide
the conservation of the integral over $v$, namely 
\[
\calV(v,w):= \int_\Omega v(x) \dd x. 
\]
We call this conserved quantity the a momentum because in the thermodynamical
interpretation below $v$ should be considered as a velocity, and it should not
be mistaken for a concentration of a diffusing species. 

In fact, $w$ should be considered as an internal energy such that the
energy density 
\[
e=  \frac12 v^2 + w
\]
plays an important role. It satisfies a conservation law without source
term, namely 
\begin{equation}
  \label{eq:ConsLaw.e}
  \dot e = \DIV\!\big( \kappa(w) \nabla w + \eta(w) v \nabla v\big) = 
\DIV\!\big( \eta(w) \nabla e + (\kappa(w){-}\eta(w)) \nabla w\big).
\end{equation}
Integration over $\Omega$ and exploiting the no-flux boundary conditions 
gives conservation of the total energy 
\[
\calE(v,w):= \int_\Omega \big(\,\frac12 v(x)^2 + w(x)\, \big) \dd x = \text{const.}
\]

\subsection{Symmetries and scaling properties}
\label{su:Scalings}

The full set symmetries of system \eqref{eq:SM01} are given for
$\Omega=\R^d$. For subsets $\Omega \neq \R^d$ only those symmetries survive
that are valid for $\Omega$. These symmetries hold for general functions $\eta$
and $\kappa$. 
\\[0.3em]
\emph{Euclidean symmetry:} For all $x_*\in \R^d$ and $Q\in \rmO(d):=\bigset{R\in
  \R^{d\ti d}}{ R^\top R=I}$ and solutions $(v,w)$ of \eqref{eq:SM01},
the rigidly moved pair
\[
\big(\wh v^{Q,x_*},\wh w^{Q,x_*}\big)(t,x):= \big(v(t, Qx{+}x_*),  w(t, Qx{+}x_*)\big) 
\]
is a solution again.
\\[0.3em]
\emph{Time and Galilean invariance:} For all $t_*\geq 0$ and $V_*\in \R$ the
time and velocity shifted pair
\[
\big(\wt v^{t_*,V_*},\wt w^{t_*,V_*}\big)(t,x):= \big(v(t{+}t_*, x)+V_*\,, \;
w(t{+}t_*,x)\big)  
\] 
is a solution again.\medskip

Now we discuss scaling properties. Observing that all terms involve either one
time derivative or two spatial derivatives we have the following invariance of
\eqref{eq:SM01}:
\\[0.3em]
\emph{Scaling S1 (parabolic scaling):} Assume $\Omega=\R^d$.  If the pair
$(v,w)$ is a solution of \eqref{eq:SM01} and $\lambda>0$, then the pair
$(v_\lambda,w_\lambda)$ is a solution as well, where
\[
(v_\lambda,w_\lambda)(t,x)= (v,w)(\lambda^2 t,\lambda x).
\]

A more complex symmetry occurs if the viscosity $\eta$ and the diffusion
constant $\kappa$ are of the same power-law type.
\\[0.3em]
\emph{Scaling S2 (nonlinear scaling):} Assume that $\eta(w)=\eta_0 w^\alpha$
and $\kappa(w)=\kappa_0 w^\alpha$ for some $\alpha>0$ and
$\eta_0,\kappa_0\geq 0$.  If the pair $(v,w)$ is a solution of
\eqref{eq:SM01} and $\mu,\lambda>0$ and $\lambda=1$ in the case
$\Omega \neq \R^d$, then the pair
$\big( v^{\mu,\lambda},w^{\mu,\lambda} \big) $ is a solution as well, where
\[
\big( v^{\mu,\lambda},w^{\mu,\lambda} \big)(t,x)=  \big( \mu\,
v(\mu^{2\alpha}\lambda^2 t, \lambda x), \mu^2 w(\mu^{2\alpha}\lambda^2 t,
\lambda x) \big).
\]
Note that the energy density $e=\frac12 v^2 + w$ scales similarly to $w$ and
the total conserved quantities satisfy 
\[
\calV(v^{\mu,\lambda},w^{\mu,\lambda})= \frac{\mu}{\lambda^d}\calV(v,w) \quad
\text{and} \quad 
\calE(v^{\mu,\lambda},w^{\mu,\lambda})= \frac{\mu^2}{\lambda^d}\calE(v,w).
\]
The main observation is that the two conserved functionals scale
differently. Hence, it is not possible to have exact similarity solutions with both,
$\calV(v,w)$ and $\calE(v,w)$ being finite and different from $0$. As
nontrivial solutions satisfy $\calE(v,w)>0$ the only choice for similarity
solutions is $\calV(v,w)=0$, see also Section \ref{se:Scaling}.  

\subsection{Exact similarity solutions}
\label{su:ExaSimiSol}

In general, the scaling symmetries can be used to transform into so-called
scaling variables via
\begin{equation}
  \label{eq:ScalTrafo}
\begin{aligned}
& \tau=\log(t{+}t_*),  \ \  y=(t{+}t_*)^{-\delta}x \quad 
 \text{and} \quad \\
&\big(v(t,x),w(t,x)\big) = \big( (t{+}t_*)^{-\gamma/2} \,\wt
 v(\tau,y)\, , \, (t{+}t_*)^{-\gamma}\, \wt w (\tau , y) \big).
\end{aligned}
\end{equation}

\paragraph*{Parabolic Scaling S1.}
According to scaling S1 have to choose $\gamma=0$ and
$\delta=1/2$ and arrive at the transformed parabolic equation 
\begin{equation}
  \label{eq:TraS1ParaSyst}
 \begin{aligned}
& \pl_\tau \wt v - \frac12\, y{\cdot} \nabla \wt v= \DIV\!\big( \eta(\wt w) \nabla
\wt v\big), \quad 
\pl_\tau \wt w- \frac12\,y{\cdot} \nabla \wt w= \DIV\!\big( \kappa(\wt w)\nabla
\wt w\big) + \eta(\wt w) \big|\nabla \wt v\big|^2.
\end{aligned} 
\end{equation}
Exact similarity solutions are steady states of this coupled system; however 
the existence of nontrivial steady states (i.e.\ with $\wt v\not\equiv 0$ and
$\wt w\not\equiv 0$) is largely open.  

Note that \eqref{eq:TraS1ParaSyst} cannot have steady states with compact support
(or more generally finite energy), because $\wt e=\wt w + \frac12\wt v^2$
satisfies $\pl_\tau \wt e   - \frac12\, y{\cdot} \nabla \wt e= \DIV( \cdots)$
and integrating over $\R^d$ and setting $\wt\calE(\tau):= \int_{\R^d}  \wt
e(\tau,y)\dd y$ gives $\frac{\rmd}{\rmd \tau} \wt\calE(\tau)=-\frac12\int_{\R^d} y{\cdot} \nabla \wt e\dd x=-\frac d2 \wt\calE
(\tau)$.

For the case $d=1$ one can obtain some results for steady states $(\wt
v(\tau,y),\wt w(\tau,y))=(V(y),W(y)$ for $y \in\R$, because the problem reduces to
the ODE system
\begin{equation}
  \label{eq:ParScaI.1D}
  -\frac y2 V' = \big( \eta(W) V'\big)', \ \  
  -\frac y2  W' = \big( \kappa(W) W'\big)' + \eta(W)
  (V')^2 \  \text{ for } y\in \R.
\end{equation}

We first introduce $c(y):=\eta(W(y))V'(y)$ such that the first equation reduces
to $c'(y)= -\big(y/\eta(W(y))\big) c(y)$. Thus, $c$ cannot change sign, which
implies that $V'$ doesn't change sign because $\eta(W)\geq 0$. Hence, for any
solution $(V,W)$ the function $y \mapsto V(y)$ must be monotone.  This implies
that for any nontrivial solution $V$ cannot lie in $\rmL^p(\R)$ for any
$p<\infty$. Nevertheless, the function $W$ may be integrable (or even have
compact support); then we can integrate the second equation in
\eqref{eq:ParScaI.1D} over $\R$ and find
\[
\frac12\int_\R W(y) \dd y = \int_\R \eta(W(y))  \big(V'(y)\big)^2 \dd y.
\]
Typical solutions of \eqref{eq:ParScaI.1D} will be such that $W$ and $V'$ are
even functions. By Galilean invariance we may assume $V(0)=0$ such that $V$ is
odd and $E=W+\frac12 V^2$ is even again. A special case is given by 
$\eta \equiv \kappa$, because by \eqref{eq:ConsLaw.e} we then have $-\frac y2
\,E' =\big( \eta(W)E'\big)'$ which has the unique even solution solution
$E\equiv E(0)$.   

\begin{theorem}[Similarity solutions]
\label{th:SimiSol} Consider the case $\eta\equiv \kappa$. Then, for each pair
$(V_-,V_+) \in \R^2$ and each $E_0> \frac12 \max\{V_-^2,V_+^2\}$ there exists a
unique solution $(V,W)$ of \eqref{eq:ParScaI.1D} satisfying 
\[
V(y) \to V_\pm \text{ for }x \to \pm\infty \quad \text{and} \quad W(y)+ \frac12
V(y)^2 =E_0 .
\]
In particular, if $V_-\neq V_+$, then $y\mapsto V(y)$ is strictly monotone and
$W(y)> E_0 - \frac12 \max\{V_-^2,V_+^2\}>0$.
\end{theorem}
\begin{proof} It suffices to solve the scalar ODE for $V$ namely 
\[
 -\frac y2 V' = \big( \eta\big( E_0{-}\frac12 V^2\big) V'\big)' \text{ for
 }y\in \R, \qquad  V(y) \to V_\pm \text{ for }y\to \pm\infty.
\]
Then, setting $W(y)=E_0-\frac12 V^2$, the pair $(V,W)$ solves \eqref{eq:ParScaI.1D}.

The existence and uniqueness of $V$ follows by applying
  \cite[Thm.\,3.1]{GalMie98DMSS} or \cite{MieSch21?ESPS}. For this may assume
  $V_- \leq V_+$ and define $A \in \rmC^{1,\mafo{Lip}}_\mafo{loc}(\R)$ via  
\begin{align*}
&A(V)= \int_{V_-}^V \eta(E_0{-}\frac12 v^2) \dd v \text{ for } V\in [V_-, V_+ ] 
 \text{ and }
\\
& A'(V)=\eta\big(E_0{-}\frac12
 V_\pm)^2\big)+ (V{-}V_\pm)  
    \text{ for } \pm V \geq \pm V_\pm .
\end{align*}
Thus, $A$ is uniformly convex with upper and lower
quadratic bound, and the above-mentioned results are applicable. Since the $V$ is
monotone its range stays inside $[V_-,V_+]$,  hence the choice of $A$ outside
of $  [V_- , V_+ ]$ is irrelevant.
\end{proof} 

\begin{figure}
\centering 
\includegraphics[height=6em]{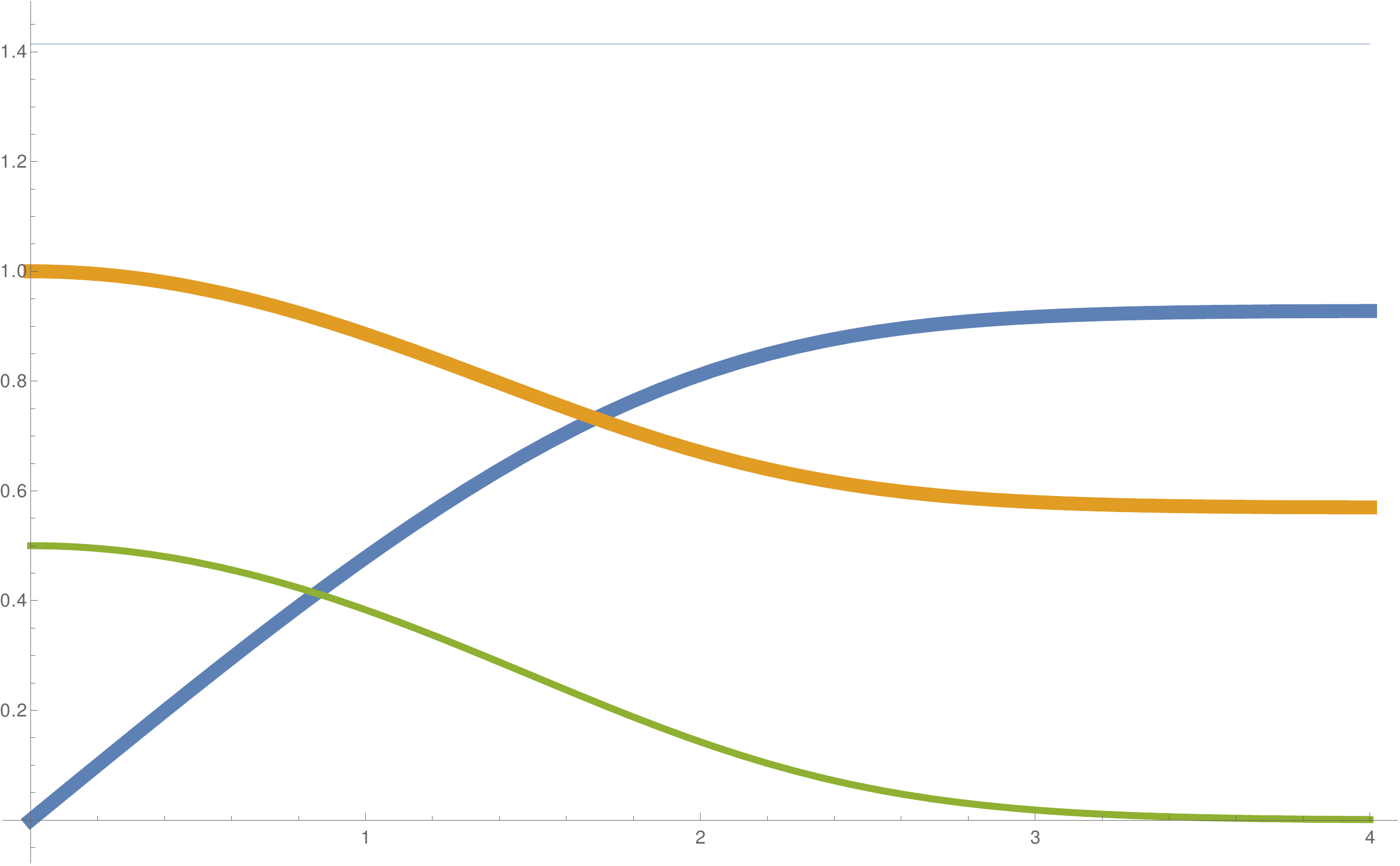} \quad
\includegraphics[height=6em]{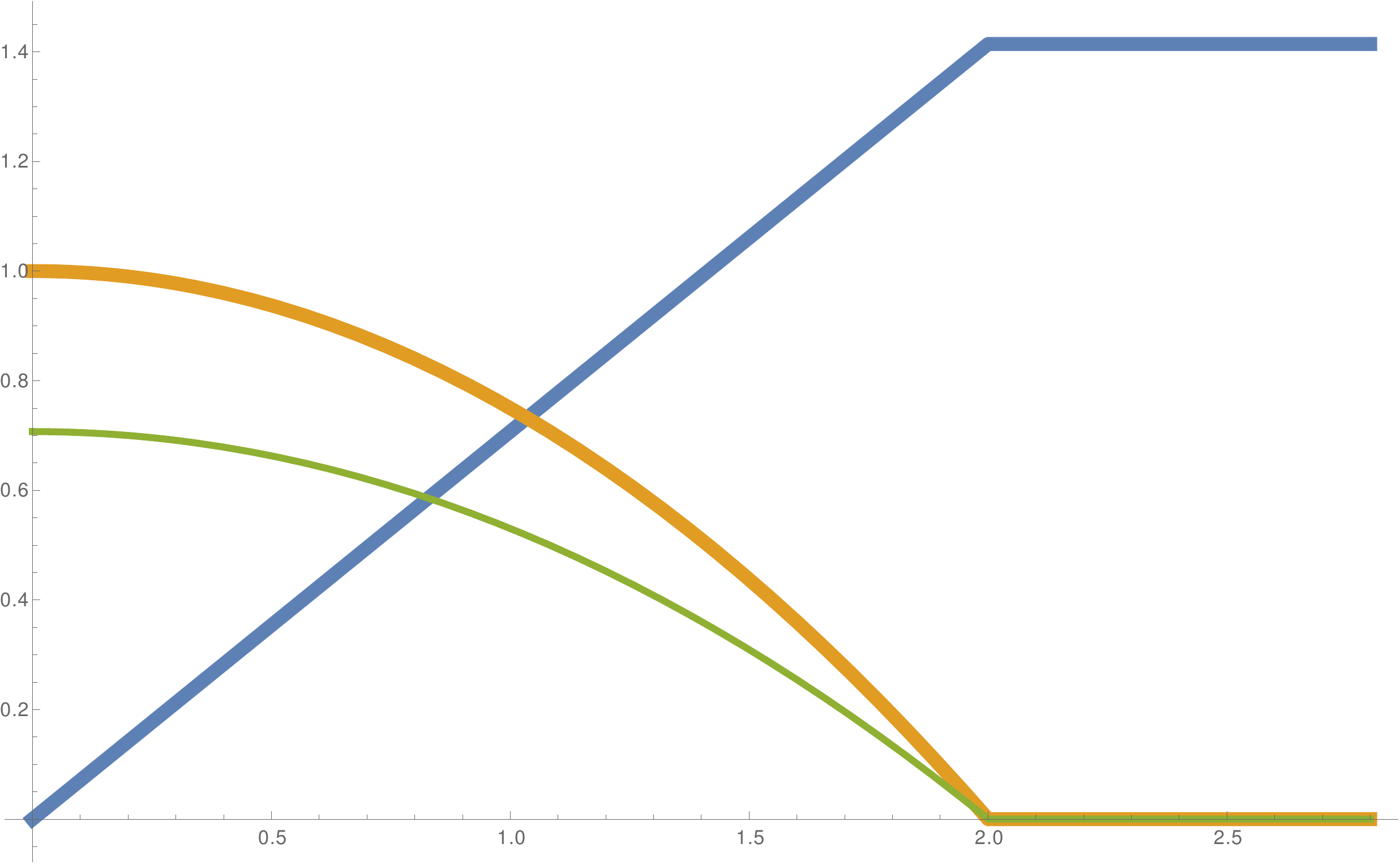} \quad
\includegraphics[height=6em]{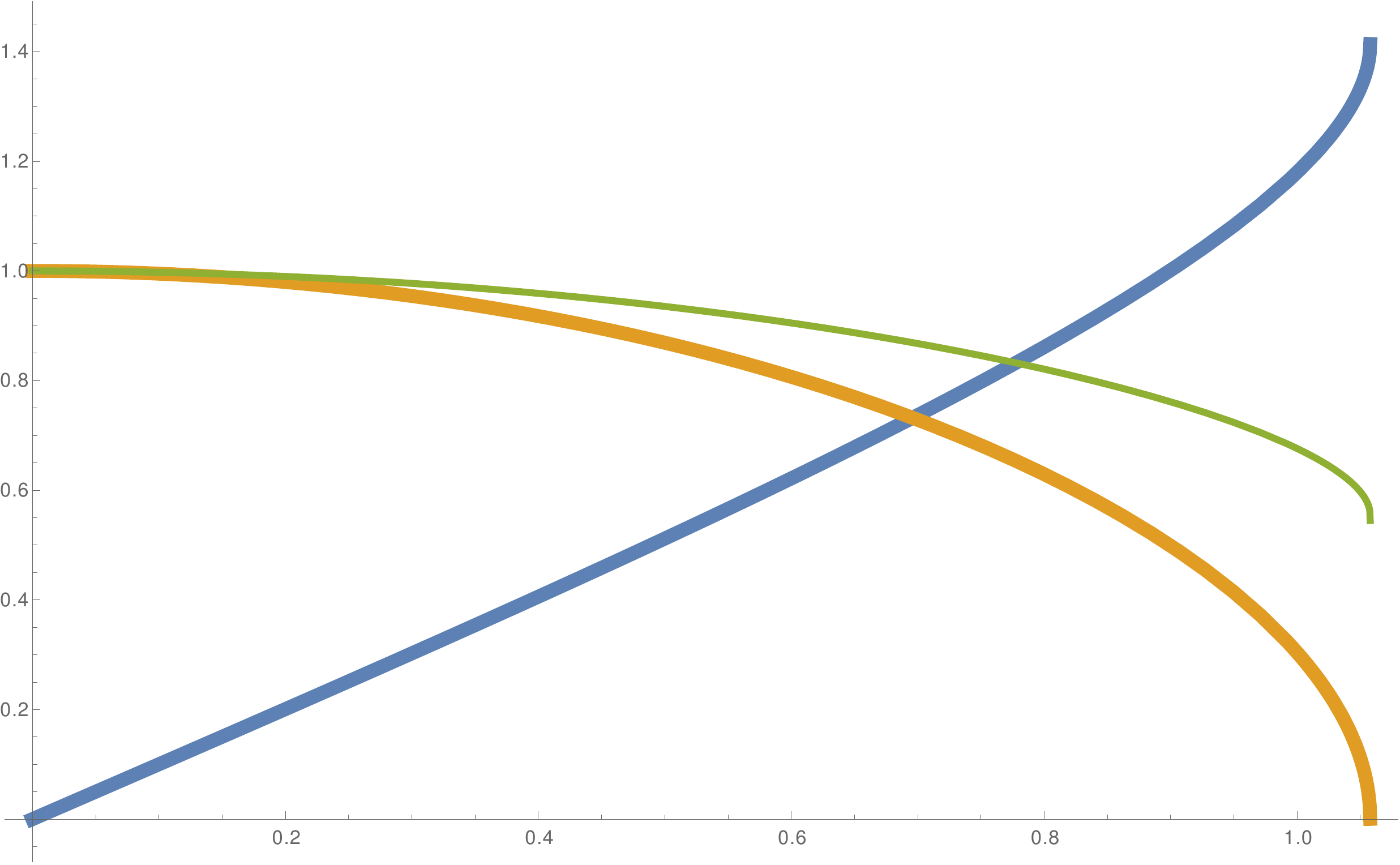} 

\caption{Solutions $y\mapsto (V(y),W(y))$ (blue,orange) of \eqref{eq:ParScaI.1D} for
  $V(0)=W'(0)=0$, $W(0)=1$, and   $V'(0)\in \{0.5,1/\sqrt2, 1.0\}$ (left,
  middle, right). The green curve displays $W(y)V'(y)$, which for $V'(0)=1.0$
  has a positive limit at $y_*$ where $W(y_*)=0$.}
\label{fig:Example.w}
\end{figure}

For an illustration of these solutions we refer to the left plot in Figure
\ref{fig:Example.w}.

Unfortunately the previous result does not apply to the degenerate case. For
this we would need $E_0 =\frac12 \max\{V_-^2,V_+^2\}$. which implies $W(y)\to
0$ for $y \to \infty$ or $y \to -\infty$. We expect that in the case
$V_-=-V_+$ and $E_0 =\frac12 V_-^2=\frac12V_+^2$ we still have a solution
$(V,W)$ and that $W$ has compact support. For the special case
$\eta(w)=\kappa(w)=w$ this can be confirmed by an explicit solution.

\begin{example}[The case $\eta(w)=\kappa(w)=w$]
\label{ex:SimiExplicit}
In this case we exploit that $E(y)=B^2$ for all $y$, where $B$ is
arbitrary.  Indeed, \eqref{eq:ParScaI.1D} has a one-parameter family of
explicit solutions:
\begin{equation}
  \label{eq:ExplicSol01}
  \big(V(y),W(y) \big)= 
\left\{ \ba{cl} \big( y/\sqrt{2}\;,\; B^2 {-}y^2/4\big ) &\text{for }|y|\leq 2B, \\
\big(\pm \sqrt2\,B\;,\; 0\big)& \text{for } \pm y\geq 2B, \ea\right.
\end{equation}
However, this solution is untypical even for our special case
$\eta(w)=\kappa(w)=w$.  To see this, we solve \eqref{eq:ParScaI.1D} as an
initial-value problem for $y \in {[0,\infty[}$ with
$(V(0),W(0),W'(0)))=(0,B^2,0)$ and $V'(0)>0$. As $V'(y)>0 $ and
$B^2=W(y)+\frac12 V(y)^2$ the solutions stay bounded with $W(y) \in [0,B^2]$
and $V(y)\in [0,\sqrt2 B]$ as long as they exist.

Starting with $V'(0)\in {]0,1/\sqrt 2[}$ we find smooth solutions with
$W(\xi)\geq W(+\infty)>0$. These are the solutions given by Theorem
\ref{th:SimiSol}.  When starting with $V'(0)>1/\sqrt 2$ the solution $(V,W)$
reaches the point $(\sqrt2 \,B,0)$ at a point $y_*\in {]0,2B[}$ with a
square-root type behavior, see Figure \ref{fig:Example.w} for some plots. In
particular, $W(y)V'(y)$ remains bounded from below by a positive constant,
which means that the solution cannot be extended by $(V(y),W(y))=(\sqrt2\,B,0)$
for $y\geq y_*$.
\end{example} 

Returning to the case of general $\eta$ and $\kappa$, it remains an open
question to discuss whether for all pairs 
$(V(\infty),W(\infty)) \in {]0,\infty[}$ there exists a unique solution $(V,W)$
($V$ odd and $W$ even) of
\eqref{eq:ParScaI.1D} that attain these limits for $\xi \to \infty$. Moreover,
one may prescribe the limits $(V(\infty),0)$ and the integral $\int_0^\infty
W(y)\dd y \in {]0,\infty[}$. Clearly none of these solutions will have finite
energy $\calE(V,W)=\int_\R(W{+}\frac12V^2)\dd y$. 
\bigskip

\paragraph*{Parabolic Scaling S2.} 
In the case $\eta(w)=\eta_0w^\alpha$ and
$\kappa(w)=\kappa_0 w^\alpha$, it is natural to search for similarity solutions
induced by the \emph{nonlinear scaling S2}. We again can use the transformation
\eqref{eq:ScalTrafo}, where we are no longer forced to use $\gamma=0$ because
we can exploit the scaling properties of $\eta$ and $\kappa$. It suffices to
chose $2\delta + \alpha \gamma=1$ to obtain an equation that is autonomous with
respect to $\tau$:
\begin{equation}
  \label{eq:TraS2ParaSyst}
 \begin{aligned}
 \pl_\tau \wt v\: -\frac\gamma2 \wt v\:  - \delta\, y\!\;{\cdot} \nabla
 \wt v\: &= \DIV\!\big( \eta_0 \wt w^\alpha\, \nabla
\wt v\big), \quad \text{ where } 2\delta+ \alpha \gamma=1, 
\\
\pl_\tau \wt w- \gamma \wt w - \delta\,y\!\;{\cdot} \nabla \wt w&= 
 \DIV\!\big( \kappa_0\wt w^\alpha\!\nabla \wt w\big) + \eta_0\wt w^\alpha
 \big|\nabla \wt v\big|^2. 
\end{aligned} 
\end{equation}
Again, the existence for nontrivial steady state solutions is totally open. 

However, we can say something for finite-energy solutions. If we look for
solutions respecting energy conservation, i.e.\ $\calE(\wt v(\tau),\wt w(\tau)) = 
\calE(v(t),w(t))=\calE(v(0),w(0)) \in {]0,\infty[}$, then we additionally have
to impose $\gamma=d\delta $. Together with $2\delta+\alpha\gamma =1$ we obtain
\[
 \delta= \frac1{2{+}d\alpha} \quad \text{and} \quad \gamma= \frac d{2{+}d\alpha} . 
\]
 
Again we can show that there are no nontrivial steady states $(\wt v, \wt w)$
with compact support. To see this, we test the steady-state equation for
$\wt v$ in \eqref{eq:TraS2ParaSyst} by $|\wt v|^{-\theta} \wt v$ for
$\theta\in {]0,1[}$. Integrating by parts the convective part on the left-hand
side and the divergence term on the right-hand side leads to the relation
\[
\big({-}\frac\gamma2 + \frac{d\delta}{2{-}\theta}  \big) \int_{\R^d} |\wt
v|^{2-\theta} \dd y  = - \int_{\R^d}  \eta_0 \wt w^\alpha (1{-}\theta)|\wt
v|^{-\theta} |\nabla \wt v\big|^2 \dd y.  
\] 
The prefactor on the left-hand side equals $\theta
d/(2(2{-}\theta)(2{+}d\alpha))>0$ whereas the right-hand side is
non-positive. Thus, we conclude $\wt v=0$ for all compactly supported steady
states. For $\wt v\equiv 0$ the system reduces to the scaled PME and it is
well-known that all similarity solutions are given by
\eqref{eq:PME.SimiSol}. 

In Section \ref{se:Scaling} we provide some evidence for our conjecture that
all finite energy solutions $(\wt v,\wt w)$ of \eqref{eq:TraS2ParaSyst} with
$\gamma =d\delta$ convergence to the corresponding steady state
$(0,W_{E_0})$.

\subsection{Gradient structure and convergence into steady state}
\label{su:Thermodyn}

We now show that the coupled system can be generated by a gradient system
$(\bfQ,\calS,\bbK)$, where $\bfQ=\rmL^2(\Omega)\ti \rmL^1(\Omega)$ is the state
space, $\calS$ is an entropy functional, and $\bbK$ is the Onsager operator
satisfying $\bbK=\bbK^*\geq 0$. The latter 
defines the dual entropy-production potential $\calP^*(v,w;\zeta_v,\zeta_w):=
\frac12 \big\langle \bfzeta, \bbK\bfzeta\big\rangle$. 
The aim is to show that the coupled system \eqref{eq:SM01} can be written
in the form 
\[
\binom{\dot v}{\dot w} = \bbK(v,w) \rmD\calS(v,w)
\]
for a suitable choice of $\calS$ and $\bbK$, see \cite{Pele14VMEG,Miel16EGCG}
for the general theory on gradient systems.
For this, we consider entropies in the form
$\calS(w)=\int_\Omega \sigma(w)\dd x $ and obtain, along solutions,
\begin{equation}
  \label{eq:EEP}
  \frac\rmd{\rmd t} \calS(w(t)) = \int_\Omega \Big\{{-}\sigma''(w)\kappa(w)|\nabla
w|^2 + \sigma'(w) \eta(w) |\nabla v|^2 \Big\} \dd x=:\mfP(v,w) .
\end{equation}
Thus, we have entropy production whenever $\sigma$ is nondecreasing and
concave. 

For finding suitable Onsager operators $\bbK$ we consider dual
entropy-production potentials in the form  
\[
\calP^*(v,w;\zeta_v,\zeta_w) := \frac12\int_\Omega \Big\{ 
p_1(w) \big|\nabla\zeta_v{-}\zeta_w \nabla v \big|^2 
+  p_2(w)|\nabla \zeta_w|^2 \Big\} \dd x 
\] 
with suitable mobilities $p_1$ and $p_2$. Here $\zeta_v$ and $\zeta_w$ are the
variables dual to $v$ and $w$, respectively. The conservation laws for $\calV$
and $\calE$ are reflected in the properties
\[
\calP^*(v,w;\rmD\calV(v,w) \equiv 0 \quad \text{and} \quad
\calP^*(v,w;\rmD\calE(v,w)) \equiv 0, 
\]
where we use $\rmD\calV(v,w)=(1,0)^\top$ and $\rmD\calE(v,w)=(v,1)^\top$. 

From the formula of $\calP^*$ we calculate $\bbK$ via $\bbK(u,k)=
\rmD^2_{\bfzeta}\calP^*(v,w;\bfzeta)$, which results in   
\begin{align*}
\bbK(v,w)\binom{\zeta_v}{\zeta_w} &= \binom{-\DIV\!\big(p_1(w)(\nabla \zeta_v{-}\zeta_w
  \nabla v) \big)}{ -p_1(w) \nabla v\cdot \nabla\zeta_v + p_1(w)|\nabla
  v|^2\varkappa -\DIV(p_2(w)\nabla \zeta_w\big)} \\
&= \bma{cc}-\DIV\!\big(p_1(w)\nabla \Box\big) & \DIV\!\big( \Box 
p_1(w)\nabla v\big)\qquad\qquad \\ - p_1(w)\nabla v\cdot \nabla\Box & \Box
p_1(w) |\nabla  v|^2  - \DIV\!\big(p_2(w)\nabla \Box\big)  \ema \binom{\zeta_v}{\zeta_w},
\end{align*}
where $\Box$ indicates the position into which the corresponding component of
$\bfzeta=(\zeta_v,\zeta_w)^\top$ has
to be inserted. 

Now, calculating $\bbK(v,w)\rmD\calS(v,w)$ with
$\rmD\calS(v,w)=(0,\sigma'(w))^\top$  we see that we obtain our coupled problem
\eqref{eq:SM01} if the relations 
\[
\eta(w)=\sigma'(w)p_1(w) \quad \text{and} \quad \kappa(w)= -\sigma''(w)p_2(w)
\quad \text{for all }w>0
\]
hold. Moreover, we see that the above entropy entropy-production relation
\eqref{eq:EEP} takes the general form
\begin{align*}
 \frac\rmd{\rmd t} \calS(v(t),w(t)) &=\big\langle \rmD\calS(v,w), \bbK(v,w)
 \rmD\calS(v,w) \big\rangle = 2 \calP^* \big( v,w;(0,\sigma'(w))^\top \big)
\\
&=\int_\Omega \Big(p_1(w)\big|0-\sigma'(w)\nabla v\big|^2  
   + p_2(w)\big|\nabla \sigma'(w)\big|^2 \Big)  \dd x =\mfP(v,w).
\end{align*}

The gradient structure for general choices of $\sigma$ can be used to obtain a
priori estimates, see Section \ref{su:DissipEstim}. Moreover, it can be used
to prove convergence into steady state on bounded domains $\Omega \subset
\R^d$. For this we observe that taking an increasing and strictly convex
$\sigma$ given the momentum $V_0:=\calV(v(0),w(0))$ and the 
initial energy $E_0:=\calE(u(0),w(0))$ there is a unique maximizer of the
entropy $\calS(v,w)$ on all states in $\rmL^2(\Omega)\ti \rmL^1_{\geq
  0}(\Omega)$ satisfying the constraints $\calV(v,w)=V_0$ and
$\calE(v,w)=E_0$, namely the spatially solutions 
\[
\wh v_{V_0,E_0}= \frac1{|\Omega|} V_0 \quad \text{and} \quad 
\wh w_{V_0,E_0}= \frac1{|\Omega|} W_0 \text{ with }W_0:=\big( E_0 -
\frac1{2|\Omega|} V_0^2\big) \geq 0.  
\]
By integrating the equation for $w$ over $\Omega$ and using $\eta(w)|\nabla
v|^2\geq 0$ and the no-flux boundary conditions, we easily obtain
\begin{equation}
  \label{eq:IntEnerMonotone}
  0\leq \int_\Omega w(0,x)\dd x \leq \int_\Omega w(t_1,x)\dd x \leq 
\int_\Omega w(t_2,x)\dd x \leq W_0 \quad \text{for }0<  t_1< t_2<\infty,
\end{equation}
where we used energy conservation for the last estimate.

Moreover, under natural conditions all solutions of $\mfP(v,w)=0$ are given by
the constant function pair $(v,w)\equiv (\INIv,\INIw)$ or by pairs of the form
$(\wt v(\cdot),0)$ for arbitrary $\wt v \in \rmL^2(\Omega)$. The latter
solutions will be excluded by assuming $\int_\Omega w(0,x)\dd x >0$. 
Hence, in good situations one
can hope for convergence of all solutions into the unique thermal equilibrium
state $(\wh v_{V_0,E_0},\wh w_{V_0,E_0})$ depending on the constraint given by
the initial conditions. 

The following results provides a first results and explains the main idea in the
simplest case, but we expect that this method generalizes to more general
situations, see \cite{MieMit18CEER,HHMM18DEER} for related convergence results
for systems of diffusion equations. It is surprising that the exponential decay
holds globally for our degenerate parabolic system and that the decay rate
$\Lambda$ is even close to the lower bound for the rate dictated by the
linearization at the steady state. 

\begin{theorem}[Convergence to steady state]\label{th:Cvg.SteadySt}
Consider a bounded Lipschitz domain $\Omega \subset \R^d$ and assume that 
\begin{equation}
  \label{eq:eta.kap.low}
  \exists\, c_\eta,c_\kappa>0\ \forall\, w\geq 0:\quad \eta(w)\geq
  c_\eta w^{1/2} \text{ and } \kappa(w) \geq c_\kappa w^{1/2}. 
\end{equation}
Let $\lambda_N:=\lambda_\rmN(\Omega)>0$ be the first nontrivial eigenvalue of
the Neumann Laplacian on $\Omega$, then all solutions converge exponentially to
a spatially constant equilibrium. More precisely, solutions with
$\calV(v(0))=V_0$, $\calE(v(0),w(0))=E_0$, and $w(0,x)\geq c>0$ a.e.\ in
$\Omega$ satisfy
\begin{equation}
  \label{eq:ExpDecay}
\begin{aligned}
& \int_\Omega \Big(\frac12\big|v(t,\cdot){-}\wh v_{V_0,E_0}\big|^2 +
\big|\sqrt{w(t,\cdot)} {-} \sqrt{\wh w_{V_0,E_0}}  \big|^2 \Big)\dd x \\
&\leq 
\ee^{-\Lambda t} \int_\Omega \Big(\frac12\big|v(0,\cdot){-}\wh v_{V_0,E_0}\big|^2 +
\big| \sqrt{w(0,\cdot)} {-} \sqrt{\wh w_{V_0,E_0}}  \big|^2 \Big)\dd x \text{
  for all }t\geq 0,
\end{aligned}
\end{equation}
where $\Lambda=\sqrt{\wh w_{V_0,E_0}}\, \min\{ c_\eta,c_\kappa\}\,  \lambda_\rmN $. 
\end{theorem}
\begin{proof} We use the entropy density function $\sigma(w)=\sqrt{w}$ and
  assume $V_0=0$ by Galileian invariance.  We introduce a relative entropy as
  an the auxiliary functional, which depends on the initial values
  $(v(0),w(0))$ via $(\wh v_{0,E_0},\wh w_{0,E_0})$ where
  $E_0=\calE(v(0),w(0))$: 
\[
\calH(v,w)= \calE(v,w) - 2\sqrt{\wh w_{V_0,E_0}} \calS(v,w) +\wh
w_{V_0,E_0}|\Omega|= \int_\Omega \!
\Big(\frac12 v^2 + \big(\sqrt{w} {-}\sqrt{\wh w_{V_0,E_0}}  \big)^2  \Big) \dd x ,
\] 
where we heavily used $\sigma(w)=w^{1/2}$ such that the linear term $w$ inside
$\calE$ equals $\sigma(w)^2$. 
Clearly, $\calH$ is a Liapunov function, and along solutions \eqref{eq:EEP}
provides the relation 
\begin{equation}
  \label{eq:dt.calH}
  \frac\rmd{\rmd t} \calH(v(t),w(t))=-2\sqrt{\wh w_{V_0,E_0}} \:\mfP(v(t),w(t)).
\end{equation}

We now show that $\mfP$ can be estimated from below by $\calH$ as follows. To
simplify notation we introduce $z=\sqrt{w}$ and $z_0=\sqrt{\wh w_{V_0,E_0}}$. 
First, we can use the explicit form of $\mfP$, $\sigma(w)=\sqrt w$,
$\nabla w=2z\nabla z$, 
and the lower estimates for $\eta$ and $\kappa$ and arrive at 
\[
\mfP(v,w)\geq \int_\Omega\!\big(\,\frac{c_\eta}2 \big| \nabla v\big|^2 +c_\kappa
\big|\nabla z\big|^2 \, \big) \dd x .
\] 
Here, we use the assumption $w(0,\cdot)\geq c>0$ with implies by Section
\ref{su:ComparisonEst} (C2) that $w(t,x)\geq c>0$ for all $(t,x)\in
{[0,\infty[} \ti \ol \Omega$.  

Secondly, setting $\ol z:=\frac1{|\Omega|}\int_\Omega z\dd x >0$  
we exploit energy conservation $\calE(v,z^2)=E_0$ along solutions and find 
\[
E_0=|\Omega|\,z_0^2 =|\Omega|\big(\ol z^2 {+} A^2\big) \text{ with }
|\Omega| A^2:=\int_\Omega \! \big( \frac12 v^2+z^2-\ol
z^2\big) \dd x =\int_\Omega \! \big( \frac12 v^2+(z {-} 
\ol z)^2 \big) \dd x .
\] 
Clearly, we have $\ol z = \big(z_0^2{-}A^2\big)^{1/2} \in [0,z_0]$, and
$z_0 \leq \big(z_0^2{-}A^2\big)^{1/2} +A $ yields $0\leq z_0-\ol z
\leq A$. 

Thus, for all $(v,z)$ satisfying $\calE(v,z^2)=E_0$ we have the following chain
of estimates
\begin{align*}
\calH(v,z^2)&=\int_\Omega \big(\frac12v^2 + (z{-}z_0)^2 \big) \dd x 
 =\int_\Omega \big( \frac12 v^2 +z^2- \ol z^2 +(z_0{-}\ol z)^2 \big) \dd x 
\\
&= |\Omega| A^2 + |\Omega|(z_0{-}\ol z)^2 \leq 2|\Omega| A^2
\\
& =\int_\Omega \big( v^2 + 2(z{-}\ol z)^2 \big) \dd x 
\leq \frac1{\lambda_N} \int_\Omega \big(|\nabla v|^2 + 2|\nabla z|^2 \big) \dd x,
\end{align*}
where we used the Poincar\'e--Wirtinger inequality for the Neumann
Laplacian on $\Omega$ and our assumption $\int_\Omega v\dd x =V_0=0$. 

Combining all the estimates we obtain the lower bound
\[
\mfP(v,w)\geq \frac{\lambda_\rmN}2 \min\{c_\eta, c_\kappa\} \,\calH(v,w) .
\] 
Together with \eqref{eq:dt.calH} and the definition of $\Lambda$ we find
$\frac\rmd{\rmd t} \calH(v(t),w(t)) \leq -\Lambda\calH(v(t),w(t))$, which
implies the desired exponential decay estimate \eqref{eq:ExpDecay}. 
\end{proof}

We emphasize that the global decay rate $\Lambda$ obtained in Theorem
\ref{th:Cvg.SteadySt} is optimal up to a possible factor of $2$, because the
linearization of the 
coupled system at the steady state $(\wh v_{V_0,E_0}, \wh w_{V_0,E_0})$ reads 
\[
\dot{\wt c} =\eta(\wh w_{V_0,E_0}) \Delta_\rmN \wt c, \quad \dot{\wt w} =\kappa
(\wh w_{V_0,E_0}) \Delta_\rmN \wt w, \quad \int_\Omega \wt c\dd x = 0 =
\int_\Omega \wt w \dd x. 
\] 
Thus, the decay rate $\Lambda$ cannot be larger than $2\min\{\eta(\wh
w_{V_0,E_0}), \kappa (\wh w_{V_0,E_0})\} \lambda_\rmN(\Omega)$, because
\eqref{eq:ExpDecay} measures quadratic distances from equilibrium.

\begin{example}[Worse decay if supports are disjoint]\slshape
\label{ex:DisjointSppt}
We want to emphasize that some assumption on the positivity of $w(t,x)$ is
necessary in Theorem \ref{th:Cvg.SteadySt}, since otherwise the stated decay
may not hold. As an example consider the case $\Omega ={]{-}\ell, \ell[}$ for
fixed $\ell \gg 1$, $\eta(w)=\kappa(w) =\frac32w^{1/2}$ such that
$c_\kappa=c_\eta=3/2$, and initial data $(\INIv,\INIw) \in 
\rmH^1(\Omega)^2$:
\[
  \INIv(x)= \mafo{sign}(x) \big(\max\{0,2|x|{-}\ell\} \big) \quad \text{and}
  \quad \INIw(x)= \Big(\frac1{15} \max\big\{ 1{-} x^2, 0 \big\} \Big) ^2 \text{ for
  } x \in \Omega.
\]
We set $T_\ell = (\ell/2)^{5/2}-1$ and observe that for $t\in [0,T_\ell]$ we
have the explicit solution 
\[
  v(t,x)= \INIv(x) \ \text{ and } \ w(t,x)= \frac1{(t{+}1)^{2/5}}
  \Big(\frac1{15} \max\big\{ 1- \frac{x^2}{(t{+}1)^{4/5}}, 0 \big\} \Big) ^2 .
\]  
Indeed, we have $\mafo{sppt}(v(t,\cdot))= \ol\Omega\setminus{]{-}\ell/2,\ell/2[}$
and $\mafo{sppt}(w(t,\cdot))= [{-}(t{+}1)^{2/5}, (t{+}1)^{2/5}]$, such that the
supports are disjoint for $t< T_\ell$. Hence, it is easy to see that both
equations are satisfied because of $\eta(w)=0$ on the support of $v(t,\cdot)$:
in the equation $v$ the dynamics is trivial with $\dot v=0$, and in the
equation for $w$ we simply have the similarity solution for the PME with
$\beta=1/2$, see \eqref{eq:PME.SimiSol}. 

To see that this solution contradicts the expoential decay estimate, it is
sufficient to calculate the involved terms only in there main term in $\ell$,
which we indicate by $\sim \ell^\alpha$:
\[
\lambda_\rmN =\frac{\pi^2}{4\ell^2} \sim \ell^{-2}, \quad V_0 = 0 ,\quad  
E_0 \sim \ell^{2}, \quad \wh v_{V_0,E_0} = 0, \quad 
\wh w_{V_0,E_0} \sim \ell^{2}.
\]
With this we find the exponential decay rate $\Lambda \sim \ell^{-1}$. Because
in \eqref{eq:ExpDecay} the integrals over $\Omega$ on both sides are of order
$\ell^3$ for all $t\in [0,T_\ell]$ and $T_\ell \sim \ell^{5/2}$, we
easily obtain a contradiction because $\Lambda T_\ell \sim \ell^{3/2}$ such
that $\ee^{-\Lambda T_\ell} \ell^3$ is smaller than $1$, whereas it
is still of order $\ell^3$ for the given solution. 

Nevertheless, we conjecture that the above decay estimate can be extended to
solutions with $\mafo{sppt}(\INIw) \subsetneqq \ol\Omega$. The point is that
solutions are not unique because of the nonlinearity $w^{1/2}|\nabla
v|^2$. Constructing solutions by approximating the initial conditions from
above as sketched in Section \ref{se:ExistTheory} one may obtain solutions
satisfying $\mafo{sppt}(w(t,\cdot)) \supset\mafo{sppt}(\nabla v(t,\cdot))$ for
all $t>0$, which then satisfy the exponential decay estimate
\eqref{eq:ExpDecay}.
\end{example}

\subsection{A related plasma model}
\label{su:Rosenau}

In a series of papers starting with \cite{RosHym85ANS, RosHym86PDAM, HymRos86ANPE}
a model of the diffusion of the mass density $\rho\geq 0$ and the heat
transport for the temperature in a plasma is developed: 
\begin{equation}
  \label{eq:Rosenau2}
  \rho_t = \DIV\big( \rho^\gamma \phi_1(\rho,\theta)\nabla \rho \big) \quad \text{and}
\quad  \big(\rho\theta\big)_t= \DIV\big(\rho^\delta
\phi_2(\rho,\theta)\nabla \theta + \theta \rho^\gamma 
\phi_1(\rho,\theta)\nabla\rho \big), 
\end{equation}
that conserves mass and energy, namely $\int_\Omega \rho(t,x)\dd
x=M_0$  and $\int_\Omega \rho(t,x)\theta(t,x)\dd x =E_0$. 

Note that the model is such that it allow one to consider a constant
temperature $\theta(t,x)=\theta_*$, and it is then sufficient to study the 
remaining PME for $\rho$:
\begin{equation}
  \label{eq:Ros.dens}
  \dot\rho = \DIV\big( \rho^\gamma \phi_1(\rho,\theta_*)\nabla \rho \big).
\end{equation}
Thus, asymptotic self-similar behavior follows if $\phi_0(,\theta_*)\neq
0$. Moreover, if $\phi(\rho,\theta_*)$ is constant, we even have Barenblatt
solutions as in \eqref{eq:PME.SimiSol}, see \cite[Eqn.\,(5)-(7)]{RosHym86PDAM}.  

To have better comparison with our model, we can introduce the energy density
$e = \rho\theta$ and obtain the system,
\begin{equation}
  \label{eq:RosSyt.r.e}
\begin{aligned}
& \dot\rho = \DIV\big( \rho^\gamma \varphi_1(\rho,e)\nabla \rho \big)  
   \quad \text{and} \quad  
\\
&  \dot e = \DIV\Big(\rho^{\delta-1} \varphi_2(\rho,e)\nabla e  + 
   e \big(\rho^{\gamma-1}\varphi_1(\rho,e) -
   \rho^{\delta-2}\varphi_2(\rho,e)\big)  
 \nabla\rho \Big),  
\end{aligned}
\end{equation}
where $\varphi_j(\rho,e)=\phi_j(\rho,e/\rho)$. Thus, our case
$\eta(w)=\kappa(w)$ corresponds here to the case $\rho^\gamma
\varphi_1=\rho^{\delta-1}\varphi_2$ leading to the special system
\[
\dot\rho = \DIV\big( \rho^\gamma \varphi_1(\rho,e)\nabla \rho \big)  
   \quad \text{and} \quad \dot e = \DIV\big( \rho^\gamma
   \varphi_1(\rho,e)\nabla e \big)  
\] 
which obviously has solutions with $e(t,x)=\rho(t,x)\theta_*$ if $\rho$
satisfies \eqref{eq:Ros.dens}. 

For an existence theory for the model \eqref{eq:Rosenau2} we refer to
\cite{BerKam90SDPE,DalGia99WSSC}. However, to the best of the knowledge of the
author a systematic study of the motion of the moving front does not exist.

\section{Steady states and traveling fronts}
\label{se:SM.TravFronts}

Here we provide a few special solutions that will highlight the coupling
between the two degenerate equations. To obtain a first feeling about the
nontrivial interaction between the two equations we study some simple explicit
solutions, namely steady states and traveling fronts.

\subsection{Steady states} 
\label{su:SteadyStates}

Steady states are all the solutions of the coupled degenerate elliptic system
\begin{equation}
  \label{eq:EllipticSyst}
\begin{aligned}   0 &= \DIV\!\big( \eta(w) \nabla v\big), &\quad 
  0 &= \DIV\!\big( \kappa(w) \nabla w\big) + \eta(w) |\nabla v|^2 &&\quad
  \text{in }  \Omega,\\
 0&=\eta(w) \nabla v\cdot \rmn &0&= \kappa(w) \nabla w\cdot \rmn
  &&\quad \text{on }\pl \Omega. 
\end{aligned}
\end{equation}
We cause of $\eta(0)=0$ it is easy to construct steady states in the form 
\begin{equation}
  \label{eq:TrivialStSt}
  (v,w) = (V_\st , 0) \quad \text{where $ V_\st$ is arbitrary.} 
\end{equation}
Thus, there is an infinite-dimensional family of \emph{trivial} steady states,
but this family is exceptional. If we assume $\INIw \not\equiv 0$. i.e.\
$\int_\Omega \INIw \dd x>0$, we find   $\calE(v(t),w(t))\geq
\calE(\INIv,\INIw)\geq 
\int_\Omega \INIw\dd x>0$ and conclude that these solutions are not relevant
any more. 

For bounded domains we have the following uniqueness result for nontrivial
steady states. 

\begin{proposition}[Steady states for $\Omega$ bounded] 
\label{pr:StStOmegaBDD} Assume that $\Omega \subset \R^d$ is a bounded
Lipschitz domain. Then, all steady states are given either by the trivial ones
in \eqref{eq:TrivialStSt} or by the spatially constant ones, namely
\begin{equation}
  \label{eq:ConstStSt}
(v,w)\ \equiv \  \big(\frac1{|\Omega|} V_0\,, \frac1{|\Omega|} W_0\big) ,  
\end{equation}
where $V_0$ and $W_0$ are uniquely given by the conserved quantities
$\calV(v,w)=V_0$ and $\calE(v,w)=E_0 = W_0+\frac1{2|\Omega|}V_0^2$, see Section
\ref{su:Thermodyn}.
\end{proposition}
\begin{proof} We use the pressure function $\wh\pi(x)=\Pi(w(x))$ with
  $\Pi(w)=\int_0^w \eta(u)  \dd u$ and observe that for a steady state $(v,w)$
  the function $\wh\pi$ satisfies the linear Neumann problem
\[
-\Delta \wh\pi = \eta(w)|\nabla v|^2 \ \text{ in }\Omega, \qquad \nabla \wh\pi \cdot
\rmn =0 \ \text{ on } \pl\Omega.
\] 
The classical solvability condition for the Neumann problem requires
$\int_\Omega \eta(w)|\nabla v|^2 \dd x =0$. Since the integrand is nonnegative
we conclude $\eta(w)|\nabla v|^2=0$ a.e.\ in $\Omega$ and $\wh\pi = \pi_*=$const. 
Because of the strict monotonicity of $\Pi$ we obtain $w=w_*=$const. and deduce
either $w_*=0$ given the trivial steady states  \eqref{eq:TrivialStSt} or
$w_*>0$ and $v=v_*=$ const. 
\end{proof}

In the case of unbounded $\Omega$ we have to allow for solutions with infinite
energy and the result is less complete. 

\begin{proposition}[Steady states for $\Omega=\R^d$] 
\label{pr:StStOmegaRd} Assume $\Omega =\R^d$ with $d \in \{1,2\}$. 
Then, all steady states (i.e.\ solutions of \eqref{eq:EllipticSyst}) are given
either by the trivial ones in \eqref{eq:TrivialStSt} or by the spatially
constant ones, namely $(v,w)\ \equiv \  (v_*, w_*) = $ const.
\end{proposition}
\begin{proof} The pressure $\wh\pi$ introduced in the previous proof still
  satisfies $-\Delta \wh\pi = \eta(w)|\nabla v|^2\geq 0$. Moreover, we know
  $\wh\pi\geq 0$. Hence, $-\wh\pi$ is a subharmonic function that is bounded from
  above.  For $d=1$ the function $-\wh\pi:\R\to {]{-}\infty,0]}$ is convex and
  hence can only be bounded if it is constant.  For $d=2$ we invoke
  \cite[Cor.\,2.3.4]{Rans95PTCP} which shows that bounded subharmonic functions
  on $\R^2$ are constant.  In both cases we conclude $\eta(w)|\nabla v|^2=0$
  a.e.\ in $\R^d$ and the result concerning $v$ follows.
\end{proof}

It is unclear whether the last result is still true in $\R^d$ with $d\geq 3$.

\subsection{Traveling fronts}
\label{su:GenPlanarW}

The importance of traveling fronts in the PME arises from the fact that they
can be used as comparison functions and that they serve as models for the local
behavior near the boundary of the support of $w$.
By isotropy of our
system \eqref{eq:SM01} it is sufficient to study the one-dimensional case $x\in
\R^1$. We start from the traveling-wave ansatz
\[
v(t,x)= V(x{+}c_\rmF t) \quad \text{and} \quad w(t,x)= W(x{+} c_\rmF t),
\]
where $c_\rmF$ is the front speed. We obtain a coupled  system of ODEs for 
$z=x+c_\rmF t$  with the unknowns
\begin{equation}
  \label{eq:TravFront}
   c_\rmF V'=\big( \eta(W)  \,V'\big)' \quad \text{and} \quad 
   c_\rmF W' = \big(\kappa(W) \,W'\big)' + \eta(W) (V')^2. 
\end{equation}
Clearly, the speed $c_\rmF$ needs to be determined together with
the nontrivial solution $(V,W)$ of \eqref{eq:TravFront}. However, we observe
that both right-hand sides in \eqref{eq:TravFront} contain two derivatives
while both left-hand sides contain only one derivative and one factor
$c_\rmF$. Hence, we can rescale solutions in such a way that for a solution 
$(c_\rmF,V,W)$ also $(\lambda c_\rmF, V(\lambda\,\cdot),W(\lambda\,\cdot))$ is
a solution for all $\lambda \in \R$. Since the case $c_\rmF=0$ leads to steady
states that were investigated already in Section \ref{su:SteadyStates}, 
it suffices to consider the case $c_\rmF=1$, only. 

To analyze the solution set of \eqref{eq:TravFront}, we integrate the first
equation obtaining the integration constant
$v_*$ and substitute the result for $\eta(W)V'$ in the second equation. Then,
the second equation can also be integrated with an integration constant $w_*$:
\begin{equation}
  \label{eq:TravWaves}
  V-v_*=  \eta(W)V', \qquad W-w_*= \kappa(W)W' +\frac12(V{-}v_*)^2. 
\end{equation}
The integration constants were chosen such that $(V,W)\equiv (v_*,w_*)$ is the
only constant solution. The system can be analyzed in the $(V,W)$ phase
plane. To see whether \eqref{eq:TravWaves} has other solutions that are defined
for all $z\in \R$, we treat for the three cases $\alpha=\beta$,
$\alpha < \beta$, and $\alpha>\beta$ separately.\medskip

\underline{\em Case $\alpha=\beta$:} For $\eta(w)=
w^\alpha$ and $ \kappa(w)= \kappa_0 w^\alpha $, there is the following
helpful observation: the parabola  
\[
\bbP_{\kappa_0}: V\mapsto \wh w(V):= w_* + \frac{1}{2(1{-}2\kappa_0)}\,(V{-}v_*)^2
\]
is invariant by the flow of \eqref{eq:TravWaves}. In the case $\kappa_0 < 1/2$
all points starting above 
$\bbP_{\kappa_0}$ stay above, which means they cannot reach $W=0$ in finite time
and hence exist for all $z\in \R$. 

\begin{figure}
\centering 
\includegraphics[width=12em]{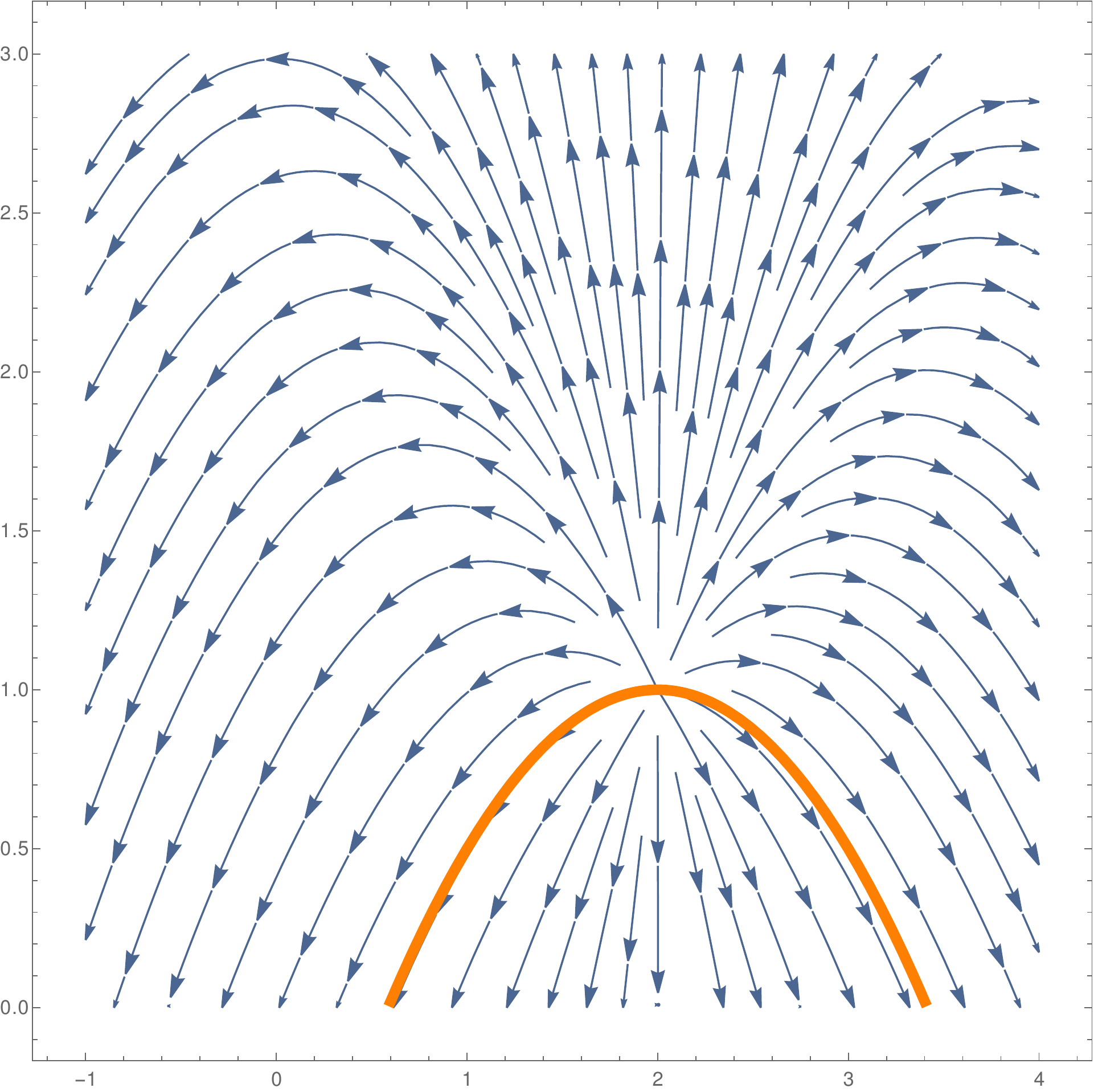}\quad
\includegraphics[width=12em]{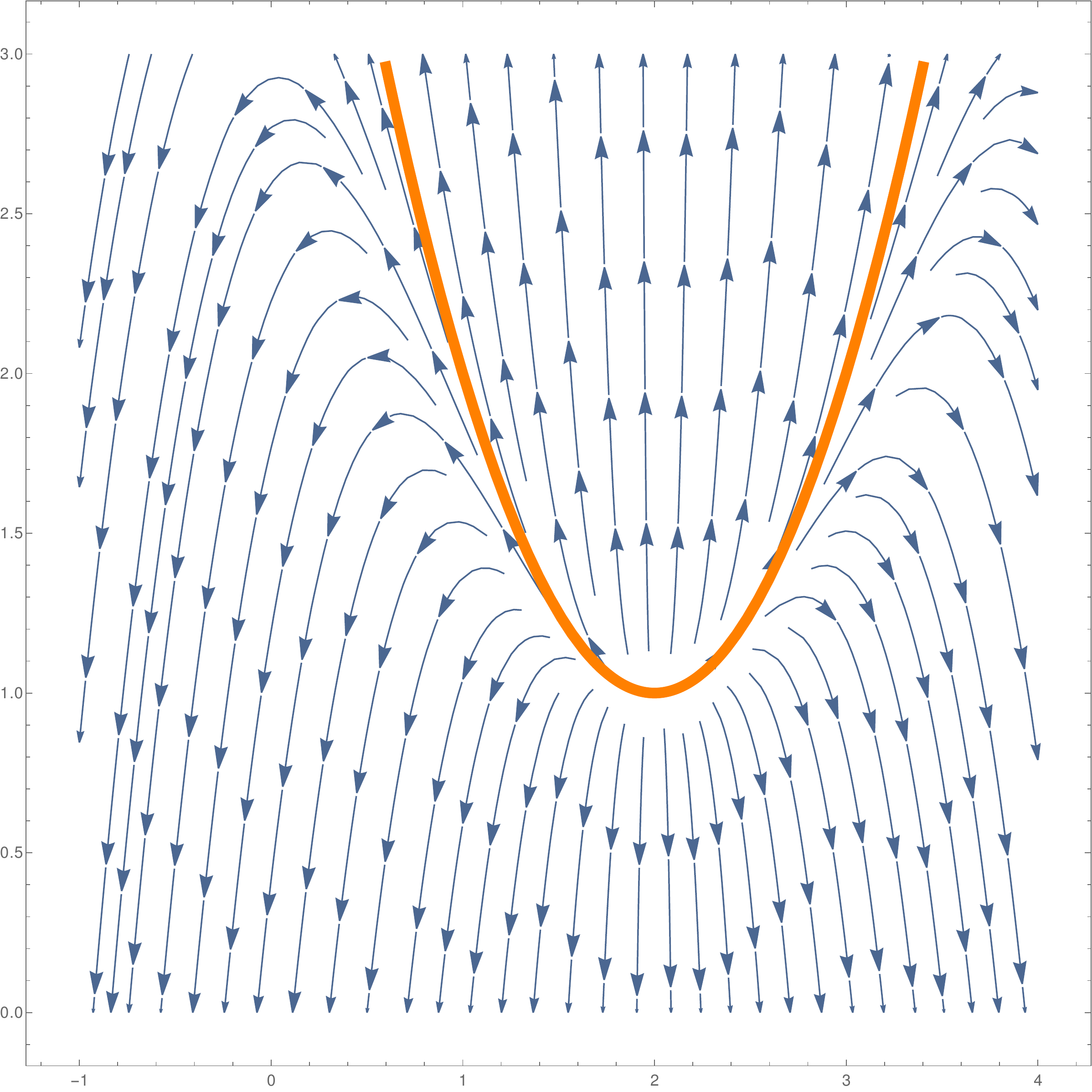}\quad
\begin{minipage}[b]{0.22\textwidth}
\caption{\small\raggedright Phase portraits for \eqref{eq:TravWaves} with
  $\kappa(w)=\kappa_0\eta(w)$ including the parabola $\bbP_{\kappa_0}$.\newline 
 Left: $\kappa_0=1$.\newline Right: $\kappa_0=1/4$.}
\label{fig:VW.parab} \end{minipage}
\end{figure}

For $w_*>0$, the solutions lying above $\bbP_{\kappa_0} $ behave as follows
(where $\rho_*=w_*^{-\alpha}$)  
\begin{align*}
&(V(z),W(z))= (v_* ,w_*) + \big(c_1\ee^{\rho_* z}, 
\frac{c_1^2\, \ee^{ 2\rho_* z} }{2(1{-}2\kappa_0)} + c_2 \ee^{\rho_* z/\kappa_0} \big) +
\text{h.o.t.\ for }z \to -\infty , 
\\
&(V(z),W(z)) = \big( c_3
z^{\kappa_0/\alpha}, (\tfrac\alpha{\kappa_0} \;\!z )^{1/\alpha} \big) +
\text{l.o.t.\ for }z\to \infty,  
\end{align*}
where $c_1,c_3\in \R$ and $c_2> 0$. The two solutions lying exactly on 
$\bbP_{\kappa_0}: W=\wh w(V) $ have a slightly different asymptotics. 
Still in the case $\kappa_0<1/2 $ one can show that all solutions starting below of
$\bbP_{\eta_0}$  reach $W=0$ in finite time and cannot be extended for all $z
\in \R$, see Figure \ref{fig:VW.parab} (right picture). 

For the case $\kappa_0\geq 1/2$, it can be shown that only one solution satisfies
$w(z)\geq 0$ for all $z$, namely the one with $V \equiv v_*$ and $W(z)>w_*$, see
Figure \ref{fig:VW.parab} (left picture). \medskip 

The situation $w_*=0$ is special, as now we can construct solutions with
$(v(z),w(z))=(v_*,0)$ for $z\leq 0$ (by Galileian invariance we can set $v_*=0$
subsequently).  These solutions are in particular interesting, because they
provide solutions with time-dependent support.  Of course, for all $\kappa_0$
we have the pure PME traveling wave
$(V(z),W(z))=\big(0, (\frac\alpha{\kappa_0} z_+ )^{1/\alpha}\big)$ where $z_+ =
\max\{0,z\}$.

For $\kappa_0 <1/2$ the two solutions lying on $\bbP_{\kappa_0}$ have the
explicit form
\[
(V(z),W(z)) = \left\{\ba{cl} \big(\pm \sqrt{2{-}4\kappa_0}\,(2\alpha
  z)^{1/(2\alpha)} , (2\alpha z)^{1/\alpha} \big)&\text{for } z\geq 0, \\ 
 \big(\ 0\ ,\ 0 \ \big)& \text{for } z\leq 0.  \ea\right. 
\]
These solutions will serve as the prototype of solutions with time-dependent
support. 

As above there are more traveling waves from the solutions lying above the
parabola $\bbP_{\kappa_0}$, which is now touching the axis $w=0$ in the origin
$(V,W)=(0,0)$. All these other solutions have the asymptotics 
\[
(V(z),W(z))= \big(c (z_+)^{\kappa_0/\alpha}, (\tfrac{\alpha}{\kappa_0}
z_+)^{1/\alpha} + \tfrac{c^2}2 (z_+)^{2\kappa_0/\alpha}\big) +
\text{h.o.t.\quad for }  z\to 0^+, 
\] 
where $c\in \R$ is a parameter for choosing the individual solutions above
$\bbP_{\kappa_0}$.\medskip

\underline{\em Case $\alpha < \beta$:} In the case $\eta(w)=w^\alpha$ and
$\kappa(w)=w^\beta$ the ODE reads
\begin{equation}
  \label{eq:VW.ODE.SpecI}
  V' = V \,W^{-\alpha}, \quad W'=\big(W{-}\frac12 V^2\big) W^{-\beta}, \qquad
  (V(0)),W(0))=(0,0).  
\end{equation}
In the cases with $\alpha\neq \beta$,  we may assume $\eta_0=\kappa_0$ without
loss of generality. 
 
\begin{figure}\centering \includegraphics[height=10em]{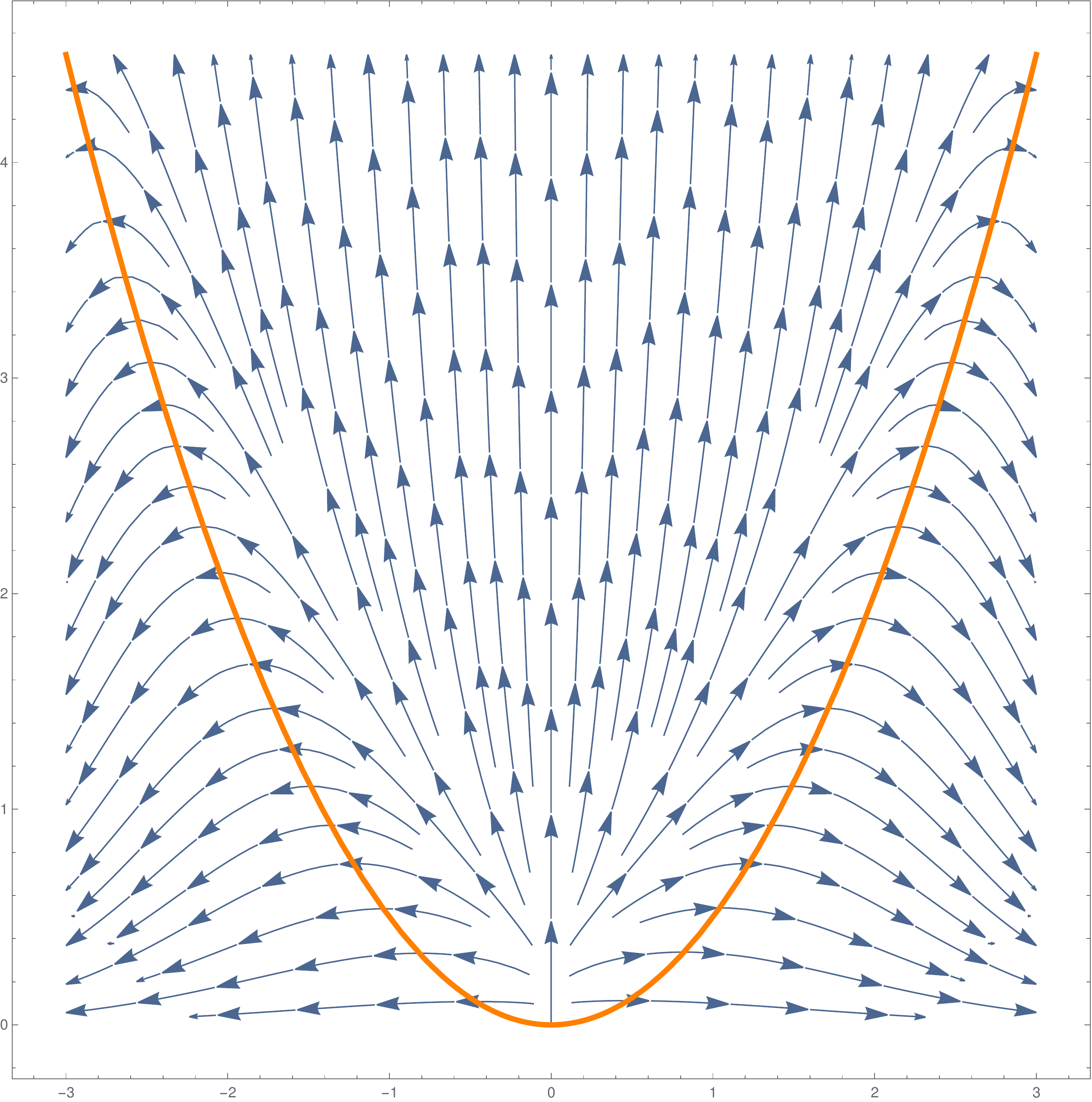} \quad 
\includegraphics[height=10em]{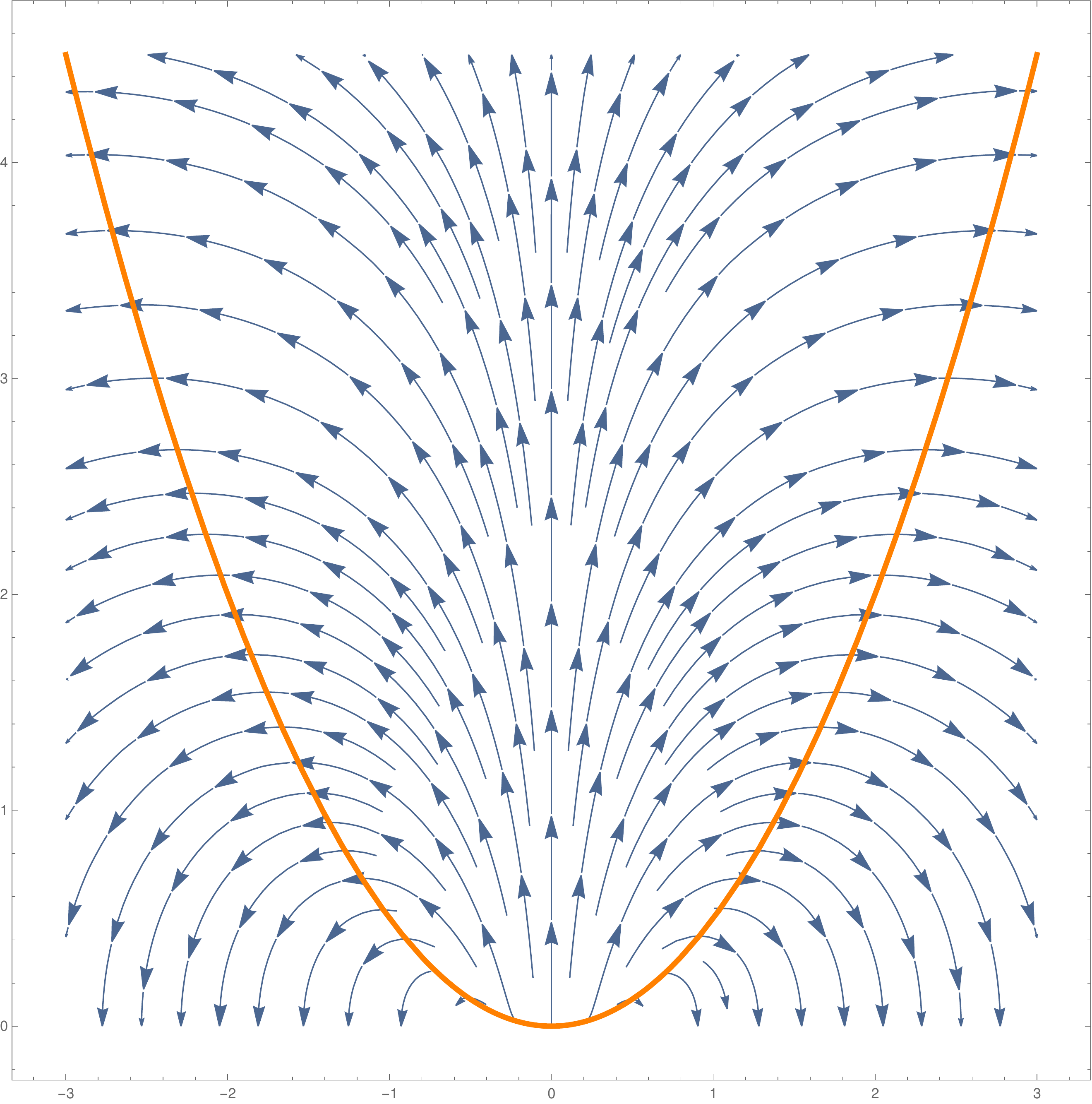}\quad
\begin{minipage}[b]{0.34\textwidth}\raggedright
\caption{Phase planes for \eqref{eq:VW.ODE.SpecI} with $(\alpha,\beta)=(2,1)$
  (left) and  $(\alpha,\beta)=(1,2)$ (right).}
\end{minipage}
\label{fig:PTF.SpecI}
\end{figure}

As usual, there is the trivial solution
with $V\equiv 0$, but all nontrivial solutions of
\eqref{eq:VW.ODE.SpecI} only exist on a subinterval of $\R$, see the right
phase plane in Figure \ref{fig:PTF.SpecI}. Nontrivial solutions first have to
lie above the parabola $W=\frac12 V^2$ to allow for $W'>0$, but after finite
time they return to the parabola and then $W'$ remains negative until
$W(z_*)=0$ is reached in finite time, where the solution ceases to exist.

The solutions can be constructed in the form $W(z)=\wt w(V(z))$ where $\wt w$
satisfies the ODE $\wt w{}'(V)= \wt w(V)-\frac12 V^2/(V \wt
w(V)^{\beta-\alpha})$. From this one sees that all solutions with $\wt w(0)=0$
satisfy the expansion 
\[
\wt w(V)= \frac12V^2 + 2^{\alpha-\beta}V^{2(1+\beta-\alpha)} + \text{h.o.t.}
\]
Inserting this into $V' =V\wt w(V)^{-\alpha}$ we obtain the expansion 
\begin{equation}
  \label{eq:VW.a.leq.b}
  \big(V(z),W(z)\big) = \Big((\alpha 2^{1+\alpha})^{1/(2\alpha)}\,z_+^{1/(2\alpha)}\, ,
  \:(2\alpha)^{1/\alpha}\, z_+^{1/\alpha} \Big)  + \text{h.o.t.}
\end{equation}
We observe that the local front behavior only depends on the smaller of the two
values $\alpha$ and $\beta$. 
\medskip
 
\underline{\em Case $\alpha > \beta$:} In this case we again obtain a one-parameter
family of traveling fronts for \eqref{eq:VW.ODE.SpecI} with half-line support.

\begin{proposition}[Fronts for $\alpha >\beta$]
\label{pr:Front.a.ge.b} 
Assume $\eta(w)=w^\alpha$ and $\kappa(w)=w^\beta$ with $\alpha >\beta$. Then
for all $v_\infty$ there exists a unique traveling front $(V,W)$ solving
\eqref{eq:VW.ODE.SpecI},
$(V(z),W(z))=(0,0)$ for $z\leq 0$, and $V(z)\to v_\infty$ for $z\to
\infty$. Moreover, for $v_\infty\neq 0$, the functions
$\mafo{sign}(v_\infty) V(\cdot)$ and $W(\cdot)$ are strictly increasing on
${[0,\infty[}$ and we have the expansion
\[
W(z)\approx (\beta z)^{1/\beta} \text{ for }z\ll 1 \text{ and for }
z\gg 1, \quad V(z) \approx c  \exp\big(-\frac\beta{\alpha{-}\beta} \!\;
z^{-(\alpha-\beta)/\beta}  \big)  \text{ for } z\ll 1. 
\]  
Moreover, there are two traveling fronts $(\pm V_\rmB,W_\rmB)$ with the expansion 
\[
(V_\rmB(z),W_\rmB(z)) \approx \big((2\alpha\, z)^{1/(2\alpha)} ,
\frac12(2\alpha\,z)^{1/\alpha } \big) \text{ for } z\gg 1. 
\]
All other solutions are not defined for all $z\in \R$, see Figure
\ref{fig:PTF.SpecI} (left).
\end{proposition}  

Again we see that the smaller of the two exponents $\alpha$ and $\beta$
dominates the local behavior of the traveling fronts. 

\subsection{Conjectured behavior near boundary of growing supports}
\label{su:Conj.GrowSppt}

Here we collect conjectured consequences of the above established behavior of
 traveling fronts. Throughout we assume that we are considering sufficiently
smooth solutions $(v,w)$ that have the property that
$S(t):=\mafo{sppt}(v(t))=\mafo{sppt}(w(t)) \Subset \Omega$.

We conjecture that $t\mapsto S(t)$ has similarly good properties as the support
of solutions of the PME, see e.g. \cite{Vazq07PMEM}. In particular, $t\mapsto
S(t)$ is non-decreasing and the boundary $\pl S(t)$ becomes smooth after a
suitable waiting time. However, in our coupled system the growth of the support
can be steered by different mechanisms depending on the relative size of
$\eta(w)$ and $\kappa(w)$ for $w\ll 1$. 

To explain the conjectured behavior in more detail, we consider a point $x_*\in
\pl S(t_*)$ and assume that $\pl S(t_*)$ is smooth. Without loss of generality
we may assume $x_*=0$ and that the the outer normal vector to $S(t_*)$ at
$x_*$ is given by $\rmn_*=-e_1$. 

In the PME \eqref{eq:I.PME} with $\kappa(w)=\kappa_0 w^\beta$ the typical
behavior (after waiting time) is that $w(t_*,x)=w_0 (x_1)_+^{1/\beta} +
\text{h.o.t.}$ The support $S(t)$ is then growing with propagation speed
$c_\rmF=\kappa_0w_0^\beta/\beta$. 

The behavior for the coupled system depends strongly on the exponents $\alpha$
and $\beta$ in $\eta(w)=\eta_0 w^\alpha$ and $\kappa(w)=\kappa_0 w^\beta$. In
all cases we will address the question of local integrability of
$\nabla v^\gamma$ and $\nabla w^\gamma$ near the boundary of $S(t)$. These
integrability properties will nicely fit together with the a priori estimates
to be derived below, see \eqref{eq:DissEstAll} in Proposition
\ref{pr:DissipEstim}. \medskip

\STEP{\underline{$\alpha>\beta$: support is driven by $w$ as in PME.}} In the
case $\alpha>\beta$ and $w\ll 1$ the energy-transport coefficient $\kappa(w)$ is
much bigger than the viscosity coefficient $\eta(w)$. Hence, $w$ will diffuse
fast and $v$ will try to keep up by following the growing support. As in the
PME the conjectured behavior (after waiting times) is
\[
(v(t_*,x),w(t,x)) = \big( 0, w_0 \:(x_1)_+^{1/\beta}\big) + \text{h.o.t.},
\]
where $v$ vanishes faster than $O(|x|^m)$ for any $m\in \N$, see Proposition
\ref{pr:Front.a.ge.b}. Again the front speed is solely controlled by $w$ alone,
namely $c_*=\kappa_0w_0^\beta/\beta$. 

The expansion does not give any information about the integrability of $\nabla
v^\gamma$, however we see that $\nabla w^\gamma \in
\rmL^p(B_r((t_*,0)))$ for $\gamma>\beta\,(1{-}1/p)$. \medskip 
        
\STEP{\underline{$\alpha<\beta$: support is driven by $v$.}} Now $\eta(w)\gg
\kappa(w)$ for $w\ll 1$, hence hence, $v$ can easily diffuse to the boundary of
the support and pile up there. We conjecture that the typical behavior (after
waiting times) is given by \eqref{eq:VW.a.leq.b}, namely 
\[
(v(t_*,x),w(t,x)) = \big( v_0 \:(x_1)_+^{1/(2\alpha)} , \frac12 v_0^2\:
  (x_1)_+^{1/\alpha}\big) + \text{h.o.t.}.
\]
The corresponding propagation speed is then given by $c_*=\eta_0
v_0^{2\alpha-2}/(\alpha 2^{\alpha+1})$. 

The reason of this behavior is that in the equation for $w$ the energy transport
via $\kappa$ can be neglected and the growth of the support is controlled by
the source term $\eta(w)|\nabla v|^2$. This leads to a equipartition of
energy near the boundary of the support giving $w\approx \frac12v^2$, or using
$e=\frac12v^2+w$ we have $\frac12 v^2 \approx \frac12 e \approx w$. 

In this case, we see that $\nabla v^\gamma \in \rmL^p(B_r((t_*,0)))$ for
$\gamma>2\alpha\,(1{-}1/p) $ and $\nabla w^\gamma \in
\rmL^p(B_r((t_*,0)))$ for $\gamma>\alpha\,(1{-}1/p)$.
\medskip 
        
\STEP{\underline{Critical case $\alpha=\beta$.}} We now consider
$\eta(w)=w^\alpha$ and $\kappa(w)=\kappa_0w^\alpha$ and will see that the both
previous cases appear, because $\kappa(w)/\eta(w)= \kappa_0$ may be large or
small depending on $\kappa_0$. 

For $\kappa_0>1/2$ the support is driven by $w$ as in the PME, however, now $v$
can follow fast enough. The conjectured behavior is 
\[
(v(t_*,x),w(t_*,x)) = (v_0 (x_1)_+^{\kappa_0/\alpha}, w_0(x_1)_+^{1/\alpha}) + \text{h.o.t.}
\] 
Here $v_0\in \R$ is arbitrary, and the propagation speed $c_* = \kappa_0
w_0^\alpha/\alpha$ depends only on $w_0$ as for the PME. 

For $\kappa_0<1/2$ the front is driven by a combination of $v$ and $w$. The
conjectured expansion takes the form 
\[
\big(v(t_*,x),w(t_*,x) \big) = \Big(v_0 \:(x_1)_+^{1/(2\alpha)}\:,\:  
\frac{v_0^2}{2(1{-}2\kappa_0)}\: (x_1)_+^{1/\alpha} \Big) + \text{h.o.t.}
\]
We observe that the leading terms of $v$ and $w$ are  coupled together
by the relation $w \approx v^2/(2{-}4\kappa_0)$. Moreover, the limit
$\kappa_0\to 0^+$ is consistent with the equipartition in the case $\alpha <
\beta$. The propagation speed is given by 
$c_* = w_0^\alpha/(2\alpha)$ which is different from $c_* =
\kappa_0w_0^\alpha/(2\alpha)$ in the case $\kappa_0>1/2$. Thus, the interaction
with the $v$ component prevents the deterioration of the wave speed for the
limit $\kappa_0\to 0^+$.  

In both subcases, we see that $\nabla v^\gamma \in \rmL^p(B_r((t_*,0)))$ for
$\gamma>2\alpha\,(1{-}1/p) $ and $\nabla w^\gamma \in
\rmL^p(B_r((t_*,0)))$ for $\gamma>\alpha\,(1{-}1/p)$.\bigskip

In summary, we find that the behavior of $w$ near the boundary of the support
is given by $w(t,x)= w_0 (x_1)_+^\gamma$ with $\gamma = \max\{
1/\alpha, 1/\beta\}$, which clearly shows that the front is driven by the
$v$-diffusion in case of $\beta >\alpha$. In the critical case $alpha=\beta$
the switch between the two regimes occurs for $\eta_0=2\kappa_0$.

\section{Weak and very weak solutions}
\label{se:WeakVeryWeakSol}

In general, we cannot expect to have strong solutions for our degenerate
coupled parabolic system. Hence, we define a suitable notions of weak and very
solutions. The problem is that the degeneracies of the viscosity $\eta$ and the
energy-transport coefficient $\kappa$, do not allows us to use parabolic
regularity, which is most easily seen for the trivial solutions
$(v(t,x),w(t,x))=(\INIv(x),0)$ that do not regularize at all. Hence, we provide a
proper definition of weak and very weak solutions in Section
\ref{su:DefWeakVeryW}, then discuss a compactly supported explicit solution in
Section \ref{su:NontrivSol}, and finally show non-uniqueness of very weak
solutions in Section \ref{su:NonuniqueVWS}.  

\subsection{Definition of weak and very weak solutions} 
\label{su:DefWeakVeryW}
 
Moreover, there is an intrinsic problem in passing to the limit in the the
``$\rmL^1$'' source term $\eta(w)|\nabla v|^2$, which typically generates a
nonnegative defect measure. This is particularly difficult because of the
degeneracies $\eta(0)=\kappa(0)$ in the viscosity $\eta$ and energy-transport
coefficient $\kappa$.  To avoid this problem, we use the strategy introduced by
Feireisl and M\'alek in \cite{FeiMal06NSET, BuFeMa09NSFS}. This means we
replace the ``partial energy equation'' \eqref{eq:SM01.b} by the equation for
the total energy $e=\frac12v^2 +w$ as given in \eqref{eq:ConsLaw.e}. Thus, we
are studying suitably defined weak solutions of the coupled system
\begin{equation}
  \label{eq:CoSy:v.e}
\begin{aligned}
\dot v= \Div\!\big( \eta(w) \nabla v\big), \quad \pl_t\big( \frac12v^2{+}w\big)
= \Div\!\big(\kappa(w) \nabla w + \eta(w) v \nabla v \big) 
\end{aligned}
\end{equation}
completed by no-flux conditions for $v$ and $w$ at the boundary $\pl\Omega$. 

Below we will give two different solution concepts, the first being a classical
weak solution. However, since $\nabla v$ only occurs together with $\eta(w)$ it
is difficult to obtain  
good a priori estimates guaranteeing that a limit $\INIv$ obtained from
$v_\eps\in \rmL^2(0,T;\rmH^1(\Omega))$ remains in that space. Hence, we define
a second concept called  very weak solutions, which can be defined for $v \in
\rmL^p(Q_T)$. For the latter we use the notion of a \emph{weak weighted
  gradient} generalizing terms of the form $a \nabla v$, where
$a$ takes the role of $\eta(w)$ and may degenerate. We 
will generalize to a mapping $G_a v$ that is valid under weak assumptions on
$v$ if $a$ is sufficiently well behaved. 

\begin{definition}[Weak weighted gradient]\label{de:WeightedGrad}
Let $q\in {]1,\infty[}$, $a \in \rmW^{1,q}(\Omega)$ and $v \in
\rmL^{q^*}(\Omega)$ with $1/q+1/q^*=1$. We say that $g \in \rmL^p(\Omega)$ with
$p\in [1,\infty]$ is the
$a$-weighted weak gradient of $v$ and write $g=G_a v$ if
\begin{equation}
  \label{eq:WeWeGrad}
  \forall\, \Psi\in \rmC^1_\rmN(\ol\Omega;\R^d):\quad 
\int_\Omega g\cdot \Psi \dd x = - \int_\Omega v\;\big( \Psi\,{\cdot} \nabla a + a
\Div \Psi\big) \dd x ,
\end{equation}
where $\rmC^1_\rmN(\ol\Omega;\R^d) := \bigset{\Psi\in \rmC^1(\ol\Omega;\R^d) }{
  \Psi\cdot \rmn=0 \text{ on }\pl\Omega }$. 
\end{definition}

As in the classical definition of weak derivatives, we see that $g$ is uniquely
defined by the pair $(a,v) \in \rmW^{1,q}(\Omega) \ti \rmL^{q^*}(\Omega)$.
Moreover, for $v \in \rmW^{1,q^*}(\Omega)$ we obviously have $G_a v= a\nabla v$
by applying Gau\ss' divergence theorem and using $\Psi\cdot\rmn=0$. However,
function for which $\nabla v$ does not exist may have a weighted gradient if
$a$ is canceling the singularity. For instance, on $\Omega ={]-1,1[}$ we may choose
$v:x\mapsto \mafo{sign}(x)$ and $a(x)=|x|^\alpha$ for $\alpha>0$, then $G_a v$
exists and equals $0$.

The next result shows that the notion of weighted gradients is stable
under limit passages, which will be crucial for constructing weak and very weak
solutions. 

\begin{lemma}[A closedness result for weak weighted gradients]\label{le:ClosWWG}
 Let $p,q\in {]1,\infty[}$ and consider $a_\eps \in
  \rmW^{1,q}(\Omega)$ and $v_\eps \in \rmL^{q^*}(\Omega)$ such that 
\[
a_\eps \weak a_0 \ \text{ in } \rmW^{1,q}(\Omega), \quad 
v_\eps \to v_0 \ \text{ in } \rmL^{q^*}(\Omega), \quad \text{and} \ 
g_\eps =G_{a_\eps} v_\eps \weak g_0 \in \rmL^p(\Omega).
\]
Then, we have $g_0=G_{a_0}v_0$.
\end{lemma}
\begin{proof}
We simply consider the defining identity \eqref{eq:WeWeGrad} for $\eps>0$ and
observe that we can pass to the limit $\eps\to 0^+$ in all three terms. On the
right-hand side it is crucial to have strong convergence of $v_\eps$. 
\end{proof}

Indeed, the notion of weak weighted gradients is implicitly defined in
\cite[Sec.\,3]{GLLMT03TSUE} and also used in \cite[Prop.\,A]{Naum13DPPT}.

We are now ready to give our two notions of solutions, where in the second
definition we have written out the definition of weak weighted gradients
explicitly to emphasize that the definition does not involve derivatives of $v$
and that the test function must be smoother. Moreover, as common in the PME we
use the pressure function
\[
   \Pi(w):= \int_0^w \kappa(s) \dd s \quad \text{ for } s\geq 0. 
\]

\begin{definition}[Weak and very weak solutions]
\label{de:WeakVeryWeak} 
Given  $T\in {]0,\infty]}$ and initial conditions
$(\INIv,\INIw)\in \rmL^2(\Omega)\ti \rmL^1_\geq (\Omega)$ we call a pair
$(v,w) \in \rmL^\infty([0,T];\rmL^2(\Omega)) \ti
\rmL^\infty([0,T];\rmL^1_\geq(\Omega))$ a \emph{weak solution} of system
\eqref{eq:SM01} if the following holds:
\begin{subequations}
  \label{eq:DefWeakSol}
\begin{align}
  \label{eq:DefWeakS.a}
& \nabla v, \ \eta(w)\nabla v,\ \eta(w) v\nabla v, \  \nabla \Pi(w) 
   \in \rmL^1(Q_T;\R^d),  \hspace{-9em}
\\[0.4em]
  \nonumber
& -\int_\Omega \INIv\varphi(0,\cdot)\dd x - \iint_{Q_T} v\pl_t\varphi \dd x \dd t
\\
\label{eq:DefWeakS.b}
&\quad = -\iint_{Q_T} \eta(w)\nabla v\cdot \nabla \varphi \dd x \dd t 
 &&\text{for all }\varphi\in \rmC^1_\rmc([0,T[\ti \ol\Omega),
\\[0.4em]
  \nonumber
& -\int_\Omega \big(\frac12(\INIv)^2{+}\INIw\big)\xi(0,\cdot)\dd x - \iint_{Q_T}\!
  \big(\frac12v^2{+}w\big)\pl_t\xi  \dd x \dd t\hspace*{-5em}
\\
\label{eq:DefWeakS.d}
&\quad= -\iint_{Q_T}\!\! \Big( \nabla\big(\Pi( w)\big)\cdot \nabla \xi 
 +\eta(w) v \nabla v\cdot \nabla \xi \Big) \dd x \dd t\hspace*{-0.5em}  && \text{for all
 }\xi\in \rmC^1_\rmc([0,T[\ti \ol\Omega).
\end{align}
\end{subequations}

A pair
$(v,w) \in \rmL^\infty([0,T];\rmL^2(\Omega)) \ti
\rmL^\infty([0,T];\rmL^1_\geq(\Omega))$ is called \emph{very weak solution} of
system \eqref{eq:SM01} if the following holds:
\begin{subequations}
  \label{eq:ADefWeakSol}
\begin{align}
  \label{eq:ADefWeakS.a}
& (v{+}v^2)\eta(w) \in \rmL^1(Q_T),\ \  (v{+}v^2)\nabla(\eta(w)),\ \nabla \Pi(w)
\in \rmL^1(Q_T;\R^d), 
\\[0.4em]
  \nonumber
& -\int_\Omega \INIv\varphi(0,\cdot)\dd x - \iint_{Q_T} v\pl_t\varphi \dd x \dd t
= \iint_{Q_T}\!\!v\,\big( \eta(w) \Delta \varphi + \nabla(\eta(w))\cdot \nabla
\varphi\big) \dd x \dd t 
\\
\label{eq:ADefWeakS.b}
&\hspace*{15em}
 \text{for all }\varphi\in \rmC^2_\rmc([0,T[\ti \ol\Omega)\text{ with }\nabla
 \varphi\cdot \rmn=0,
\\[0.4em]
  \nonumber
& -\int_\Omega \big(\frac12(\INIv)^2{+}\INIw\big)\xi(0,\cdot)\dd x - \iint_{Q_T}\!
  \big(\frac12v^2{+}w\big)\pl_t\xi  \dd x \dd t\hspace*{-5em}
\\
\label{eq:ADefWeakS.d}
&= \iint_{Q_T}\!\! \Big(\! {-}\nabla(\Pi(w))\cdot \nabla \xi 
 +\frac{v^2}2\big(\nabla (\eta(w))\cdot\nabla \xi  +\eta(w) \Delta \xi\big)
 \Big) \dd x \dd t\hspace*{-5em}  \\
\nonumber
&\hspace*{15em} \text{for all
 }\xi\in \rmC^2_\rmc([0,T[\ti \ol\Omega) \text{ with }\nabla \xi\cdot \rmn=0.
\end{align}
\end{subequations} 
\end{definition}

We first observe that functions of the form $(v(t,x),w(t,x))=(\INIv(x),0)$ with
$\INIv\in \rmL^2(\Omega)$ are very weak solutions. Because of $w\equiv 0$ it is
trivial to see that \eqref{eq:ADefWeakSol} is satisfied.  

We remark that weak solutions are not necessarily very weak solutions, because
the degeneracies do not allow us to transfer the necessary integrabilities
easily. 

For both notions of solutions,  we have conservation of momentum and energy. To
see this, we simply consider spatially constant test functions
$(t,x)=\varphi(t)$ and $(t,x)=\xi(t)$. As the spatial gradients of the test
functions vanish, we obtain 
\[
\calV(\INIv)\varphi(0)=\int_0^t \calV(v(t)) \dot\varphi(t) \dd t 
\quad \text{and} \quad 
\calE(\INIv,\INIw)\xi(0)=\int_0^t \calE(v(t),w(t)) \dot\xi(t) \dd t. 
\]
By the lemma of Du Bois-Reymond we conclude $ \calV(v(t))= \calV(\INIv)$ and
$\calE(v(t),w(t))  =  \calE(\INIv,\INIw)$ for a.a.\ $t\in [0,T]$. 

Of course, weak or very weak solutions that are sufficiently regular are even
strong solutions. 

\subsection{An explicit compactly supported solution}
\label{su:NontrivSol}

We consider the case $\Omega=\R$ and $\eta(w)=\kappa(w)=w$ and provide a
nontrivial solution with compact support that grows in time. The solution is
obtained by combining two self-similar solutions as discussed in Example
\ref{ex:SimiExplicit} in a suitable way. 

We choose real positive parameters $B$, $x_*$, and $t_*$ such that $t_* <
x_*^2/(4B^2)$ and set $T:=x_*^2/(4B^2) - t_*>0$. For $(t,x)\in [0,T]\ti \R$ we
define the functions $(v^*,w^*)$ via 
\begin{equation}
\label{eq:ExplSol*}
\begin{aligned}
v^*(t,x)&= \left\{ \ba{cl}
   \frac{x{+}x_*{+}2B\sqrt{t{+}t_*}}{\sqrt2\,\sqrt{t{+}t_*} }  
   & \text{for } |x{+}x_*| \leq 2B\sqrt{t{+}t_*} ,
 \\[0.4em]
 2 \sqrt 2\: B & \text{for }|x| \leq x_*-2B\sqrt{t{+}t_*},
\\[0.4em] 
\frac{x_*{+}2B\sqrt{t{+}t_*}\,{-}x}{\sqrt2\,\sqrt{t{+}t_*} } 
  & \text{for } |x{-}x_*| \leq 2B\sqrt{t{+}t_*}
,\\[0.4em] 0 & \text{otherwise}.  \ea \right.
\\ 
w^*(t,x)&= \left\{ \ba{cl} B^2-\frac{(x{+}x_*)^2}{4(t{+}t_*)}&
\text{for } |x{+}x_*| \leq 2B\sqrt{t{+}t_*} ,\\
B^2-\frac{(x{-}x_*)^2}{4(t{+}t_*)}& \text{for } |x{-}x_*| \leq 2B\sqrt{t{+}t_*}
,\\ 0 & \text{otherwise}.  \ea \right.
\end{aligned}
\end{equation}
\begin{figure}
\centering
\begin{tikzpicture}
\draw[thick,->] (-4,0)--(4.5,0) node[below]{$x$};
\draw[thick,->] (0,-0.3)--(0,1.7);
\draw[very thick, green!70!blue] (-3.8,0.06)--(-2.5,0.06)--(-1.5,1.414)
   --node[pos=0.8, above]{$v^*(t,x)$}(1.5,1.414)--(2.5,0.06)--(4.0,0.06);
\draw[very thick, red!80!blue] (-3.8,0.02)--(-2.5,0.02);
\draw[red!80!blue, very thick]   plot[smooth,domain=-1.0:1.0] (-2.0-\x/2.0, {1.0-(\x)^2});
\draw[very thick, red!80!blue] (4.0,0.02)--(2.5,0.02);
\draw[red!80!blue, very thick]   plot[smooth,domain=-1.0:1.0] (2.0-\x/2.0, {1.0-(\x)^2});
\draw[very thick, red!80!blue] (-1.5,0.02)--node[pos=0.77, above]{$w^*(t,x)$}(1.5,0.02);
\draw[thick] (-2,0.1)--(-2,-0.1) node[below]{$-x_*$};
\draw[thick] (2,0.1)--(2,-0.1) node[below]{$x_*$};
\end{tikzpicture}
\begin{minipage}[b]{0.3\textwidth}
\caption{Graph of the explicit solution $(v^*,w^*)$ from \eqref{eq:ExplSol*} 
having growing support
  $[-\xi(t),\xi(t)]$ where $\xi(t)=x_*+2B\sqrt{t{+}t_*}$.}
\label{fig:ExplSol}
\end{minipage}
\end{figure}
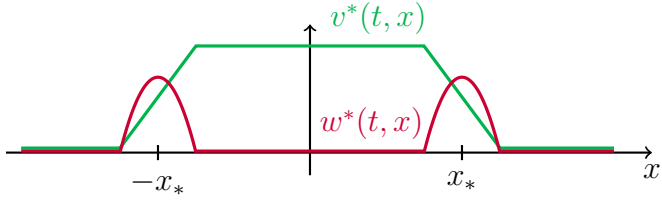%
A direct calculation shows that $(v^*,w^*)$ is a strong as well as a weak
solution.  Moreover, we have $e^*:=\frac12 (v^*)^2 + w^*= \sqrt 2
B\,v^*$, which is consistent with the fact that $\eta\equiv \kappa$ implies that
$e^*$ and $v^*$ satisfy the same equation, namely $\dot e^* = (w^* e^*_x)_x$ and $\dot v^*
=(w^*v^*_x)_x$.  

A simple calculation using the piecewise linear structure of $v^*(t,\cdot)$ and
the piecewise parabolic structure of $w^*(t,\cdot)$ gives the relations 
\[
\int_\R\!\! v^*(t,x)\dd x = \sqrt{32}\, B x_*, \ \ \int_\R\!\frac{(v^*)^2}2\dd x = 8B^2 x_*
{-} \tfrac{16}3B^3\sqrt{t{+}t_*} , \ \ \int_\R\!\!  w^*(t,x) \dd x = \tfrac{16}3 B^3
\sqrt{t{+}t_*}\:. 
\]  
This confirms the conservation of the total momentum and the total energy. 

We also note that the source term $\eta(w^*)(v^*_x)^2 $ reduces here to the simple
expression $w^*/(2t{+}2t_*)$, which vanishes at the boundary of
$\mafo{sppt}(w^*(t,\cdot))$. Hence, in this case the source term does not
contribute to the growth of the support.

\subsection{Nonuniqueness of very weak solutions}
\label{su:NonuniqueVWS}

We now consider the pair $(v^0,w^0)$ that is obtained from $(v^*,w^*)$ by
keeping $B$ and $x_*$ fixed but taking the limit $t_*\to 0^+$.  Then,
$(v^0,w^0)$ has the initial values $(\INIv_0,\INIw_0)= ( 2\sqrt 2\, B\, 
\bm 1_{[-x_* ,x_*]},0) \in \rmL^2(\Omega)\ti \rmL^1_\geq(\Omega)$.  
Next we observe that 
$\nabla v=v_x$ is piecewise constant with values $\pm 1/\sqrt{2t}$ in the
intervals $|x{\pm} x_*|\leq 2B\sqrt{t}$ and $0$ otherwise. Hence, we have
$\nabla v\in \rmL^p(Q_T)$ for $p\in {[1,3[}$. With $\Pi(w)=w^2/2$ we find 
$\nabla \Pi(w) = w w_x \in \rmL^q(Q_T)$ for all $q\in {[1,2[}$. Using
$(v^0,w^0)\in \rmL^\infty(Q_t)$ we have the conditions \eqref{eq:DefWeakS.a}
and \eqref{eq:ADefWeakS.a}. Moreover, inserting $(v^0,w^0)$ into the weak form 
\eqref{eq:DefWeakS.b}+\eqref{eq:DefWeakS.d} of the very weak form
\eqref{eq:ADefWeakS.b}+\eqref{eq:ADefWeakS.d} we can use the explicit formula
for $(v^0,w^0)$ to undo the integrations by parts and see that $(v^0,w^0)$ is
indeed a weak solution as well as a very weak solution. 

However, there is the trivial second very weak solution, namely 
$(v(t,\cdot),w(t,\cdot))=(\INIv_0,0)= ( 2\sqrt 2\, B\,\bm 1_{[-x_*
, x_*]},0)$.  Thus, we definitely have nonuniqueness in the class of very weak
 solutions. 

Indeed, we have a two-parameter family of very weak solutions for the initial
conditions $(\INIv_0,0)= ( 2\sqrt 2\, B\,\bm 1_{[-x_*
  ,x_*]},0)$. The point is that we may keep the solution constant in time for
an arbitrary $t_+>0 $ for $x\geq 0$ and then start with a delayed version of 
$(v^0,w^0)$. Moreover, we may choose $t_->0$ for starting starting a
delayed version of $(v^0,w^0)$ for $x\leq 0$. More precisely, we choose $t_+,t_-\in
[0,T]$ and set 
\[
(\wt v(t,x),\wt w(t,x)) = \left\{\ba{cl} 
(\INIv_0(x),0) & \text{for } x\geq 0 \text{ and } t\in [0,t_+], \\ 
(v^0(t{-}t_+,x),w^0(t{-}t_+,x)) & \text{for } x\geq 0 \text{ and } t\in [t_+,T], \\ 
(\INIv_0(x),0) & \text{for } x\leq 0 \text{ and } t\in [0,t_-], \\ 
(v^0(t{-}t_-,x),w^0(t{-}t_-,x)) & \text{for } x\leq 0 \text{ and } t\in
[t_-,T]. 
\ea  \right.
\]
We emphasize that the different delays do not produce any nonsmoothness,
because we have $(v^0(t,0),w^0(t,0))=(2\sqrt2 B,0)$ for all $t\in [0,T]$.  A
direct calculation shows that $(\wt v,\wt w)$ is a very weak solution for all
the choices of $t_+,t_-\in [0,T]$.

\subsection{Case $\eta\equiv \kappa$: families of solutions with growing
  support}\label{su:eta=kappa}   

In the case $\kappa \equiv \eta$ we have the additional simple equation for $e= \frac12
v^2 +w$, namely
\[
\dot e = \DIV\big(\kappa(w)\nabla e), \quad \kappa(w)\nabla e \cdot \rmn=0.
\]
Thus, we obtain exactly the same equation as for $v$ and may restrict to
a solution class defined via the relation $e=Bv$ for some fixed constant $B>0$.
This leads to the relation  $w=e-\frac12 v^2 = Bv {-}\frac12 v^2 $, where we now
have the restrictions $v \in [0,2B]$ and $w\in [0,B^2/2]$. With this the pair
$(v,w)=(v, Bv {-}\frac12 v^2)$ solves the coupled system \eqref{eq:SM01} if and
only if $v$ solves the scalar equation 
\begin{equation}
  \label{eq:v.scalar}
  \dot v = \DIV\!\big( \wt\kappa_B(v) \nabla v\big), \quad \wt\kappa_B(v)
  \nabla v\cdot \rmn=0 , \quad \text{where }\wt\kappa_B(v):=\kappa \big(Bv
  {-}\tfrac12 v^2\big).  
\end{equation}
For this scalar equation the general existence theory for the PME (cf.\
\cite{Vazq07PMEM}) can be applied and a huge set of solutions with
compact and growing support for $v$ are known to exist. 

For example, we may choose $\eta \equiv \kappa$ in the form 
\begin{equation}
  \label{eq:mu=nu}
  \kappa(w)=\eta(w)= \left\{ \ba{cl} B - \sqrt{B^2{-}2w} &\text{for }w\in [0,3B^2/8], \\
 2w/B - B/4& \text{for }w\geq 3B^3/8. \ea \right. 
\end{equation} 
Then, for $v\in [0,B/2]$ we have $\wt \kappa_B(v)=\kappa\big(Bv{-}v^2/2\big)=v$
and arrive at the classical PME $\dot v=\Div(v\nabla v)=\Delta(v^2/2)$, which
has explicit similarity solutions satisfying $v\in [0,B/2]$ for sufficiently
large $t$, see \eqref{eq:PME.SimiSol}.

\section{General a priori estimates}
\label{se:GenApriooriEst}

To provide a first existence theory for the coupled system \eqref{eq:SM01}, 
we restrict now to
the case that $\Omega \subset \R^d$ is a smooth and bounded domain. This will
simplify certain compactness arguments. Moreover, throughout this section we
will assume that the solutions $(v,w)$ are classical solutions. 

To derive general a priori estimate we consider a smooth
function $(v,w)\mapsto \varphi(v,w)$ and find for solutions $(v,w)$ of
\eqref{eq:SM01} the relation  
\begin{subequations}
  \label{eq:varphi.PDE}
\begin{equation}
  \label{eq:varphi.PDE.a}
\pl_t \varphi(v,w)= \DIV\!\big( \eta(w) \varphi_v(v,w) \nabla v + \kappa (w) 
\varphi_w(v,w) \nabla w) \big) + R_\varphi(v,w,\nabla v,\nabla w), 
\end{equation}
where the remainder $R$ is given explicitly via 
\begin{align}
\label{eq:varphi.PDE.b}
R_\varphi(v,w,\nabla v,\nabla w) &= \eta(w)\,(\varphi_w(v,w) {-}
\varphi_{vv}(v,w))\big| \nabla v \big|^2\\
& \nonumber  \quad  -
(\eta(w){+}\kappa(w))\, \varphi_{vw}(v,w)\,  \nabla v {\cdot}\nabla w - \kappa(w)\,
\varphi_{ww}(v,w)\big|\nabla w \big| ^2.  
\end{align}
\end{subequations}
Integrating \eqref{eq:varphi.PDE.a} and using the boundary condition we find
along solutions 
\begin{equation}
  \label{eq:varphi.R}
  \frac\rmd{\rmd t} \int_\Omega \varphi(v,w) \dd x = \int_\Omega
R_\varphi(v,w,\nabla v,\nabla w) \dd x.
\end{equation}

\subsection{Estimates for the $\rmL^p$ norms}
\label{suu:IntegralEstim}

Clearly, choosing $\varphi(v,w)=v$ or $\varphi(v,w)=\frac12v^2 + w$ gives
$R_\varphi \equiv 0$, which is the conservation of momentum $\calV$ and energy
$\calE$ as discussed in Section \ref{su:ConservLaws}.

Moreover, choosing $\varphi(v,w)=\phi(v)$ we obtain 
\[
\frac\rmd{\rmd t} \int_\Omega \phi(v(t,x))\dd x = - \int_\Omega \phi''(v) 
\eta(w) \big|\nabla v \big|^2 \dd x 
\]
Hence, for all convex functions $\phi$ we obtain that $\int_\Omega
\phi(v(t,x))\dd x$ is nondecreasing in $t$. This implies the decay of all
$\rmL^p$ norms, namely 
\begin{equation}
  \label{eq:LpNormDecayV}
  \forall\, p \in [1,\infty]\  \forall\, t\geq 0: \quad 
\| v(t,\cdot)\|_{\rmL^p(\Omega)} \leq \| v(0,\cdot)\|_{\rmL^p(\Omega)}.
\end{equation}

For the $w\geq 0$ we obviously have an a priori bound in $\rmL^1(\Omega)$ via
\[
\int_\Omega w(t) \dd x \leq \calE(v(t),w(t)) =  \calE(\INIv,\INIw). 
\]
However, because of the $\rmL^1$-right-hand side $\eta(w)|\nabla v|^2$ it is
difficult to derive a priori bounds for high $\rmL^p$ norms. 

The class $\varphi(v,w)=\Phi(w)$ is also important and leads to the relation
\begin{equation}
  \label{eq:Phi(k)}
\frac\rmd{\rmd t} \int_\Omega \Phi(w(t)) \dd x = \int_\Omega \big\{
\eta(w)\Phi'(w)|\nabla v|^2 -\kappa(w)\Phi''(w)|\nabla w|^2  \big\}
\dd x .
\end{equation} 
Thus, we have growth of the integral functional if $\Phi'\geq 0$ and $\Phi''\leq
0$, e.g.\ for $\Phi(k)=-k^\alpha$ with $\alpha\in [0,1]$. Such functionals
include the physically relevant entropies discussed in Section \ref{su:Thermodyn}.

\begin{remark}[$\rmL^p$ norms for $\eta \equiv \kappa$]
\label{re:LpNormsE=K}
For the case $\eta \equiv \kappa$ we are in a special situation, where we can use
$\varphi(u,k)=\phi\big(\frac12u^2{+}k)$ to obtain 
\begin{equation}
  \label{eq:nu=mu.sigma}
  \frac\rmd{\rmd t} \int_\Omega \phi\big(\frac12u^2{+}k\big) \dd x =
  -\int_\Omega \eta(k)\, \phi''\big(\frac12u^2{+}k\big)\: \big|u\nabla u {+} 
  \nabla k \big|^2 \dd x.  
\end{equation}
Thus,  we have decay for all convex $\phi$, whereas concave $\phi$ leads to 
growth. This estimate is also easily derived from the simple equation $\dot e =
\DIV\!\big(\eta(w)\nabla e\big)$, which holds for $\eta \equiv \kappa$.  
Together with \eqref{eq:LpNormDecayV} we find
\begin{align}
  \label{eq:LpNormDecayW}
  \forall\, p \in [1,\infty]\  \forall\, t\geq 0: \quad 
\| w(t)\|_{\rmL^p(\Omega)}& 
\leq \| e(t)\|_{\rmL^p(\Omega)} 
\\
& \nonumber 
\leq \| e_0\|_{\rmL^p(\Omega)}
\leq \| \INIw\|_{\rmL^p(\Omega)} + \frac12 \| \INIv\|^2_{\rmL^{2p}(\Omega)}.
\end{align}
\end{remark}

To see that $\rmL^p$ bounds for $w$ can also be derived for cases without
$\eta\equiv \kappa$ we consider now the situation $\kappa(w)=\kappa_0\eta(w)$. 
For this case we can write $R_\varphi$ as a multiple of $\eta(w) $ in the form 
\[
R_\varphi(v,w,\nabla v,\nabla w)= - \eta(w) \binom{\nabla v}{\nabla w} \cdot
\bbA^{\kappa_0}_\varphi \binom{\nabla v}{\nabla w} \ \text{ with } 
\bbA^{\kappa_0}_\varphi = \bma{cc} \varphi_{vv}-\varphi_w& \frac{1+\kappa_0}2
\varphi_{vw} \\ \frac{1+\kappa_0}2
\varphi_{vw} & \kappa_0 \varphi_{ww}  \ema .
\] 
Thus,  it suffices to show that $\bbA^{\kappa_0}_\varphi(v,w)$ is positive
semidefinite for all arguments to obtain decay estimates for $\int_\Omega
\varphi(v,w)\dd x$.  

\begin{lemma}[Higher norms in Spec.\,II]
\label{le:HighNorms} For all integers $m\in \N$ and all $\kappa_0>0$ 
there exists coefficients $a_j\geq 0$ for $j=1,...,m$ such that 
$\varphi_m(v,w)= w^m+\sum_{j=1}^m a_j w^{m-j}v^{2j}$ satisfies
$\bbA^{\kappa_0}_{\varphi_m}(v,w)\geq 0$ for all $(v,w)\in \R\ti {[0,\infty[}$.   
\end{lemma}
\begin{proof} The case $m=1$ is solved by $\varphi_1(v,w)=w + \frac12 v^2$
  independently of $\kappa_0$. For case $m=2$ a direct calculation for
  $\varphi(v,w)= w^2+a_1 wv^2 + a_2v^4$ yields
\[ 
\bbA^{\kappa_0}_{\varphi_2} = \bma{cc} (12a_2{-}a_1)v^2 + (2a_1{-}2)w &
(1{+}\kappa_0) a_1 v\\ (1{+}\kappa_0) a_1 v& 2\kappa_0 \ema ,
\]
which is positive semi-definite for all $(v,w)$ if and only if $a_1\geq 1$ and
$2\kappa_0(12a_2{-}a_1)\geq (1{+}\kappa_0)^2 a_1^2$. Clearly, it suffices to choose
$a_1=1$ and then $a_2\geq 1/12 + (1{+}\kappa_0)^2/(24\kappa_0)$ to fulfill all
requirements in the case $m=2$.

For general $m\geq 3$ we set $e=w+\frac12v^2$ and define
$\phi^m_j(v,w):=e^{m-j} v^{2j}$, 
\[
\varphi_m(v,w)= \sum_{j=0}^m b_j \phi_j^m(v,w) , 
 \quad \bbV^{\kappa_0}:=\bma{@{}cc@{}}v^2 &\frac{1+\kappa_0}2 v\\ 
 \frac{1+\kappa_0}2 v &\kappa_0 \ema, \quad 
 \bbE:=\bma{@{}cc@{}}1&0\\0&0\ema.   
\]
With $\kappa_*=(1{-}\kappa_0)^2/(4\kappa_0)$ we obtain $\bbV^{\kappa_0} \geq -\kappa_* v^2
\bbE$ in the sense of positive definite matrices. Since $\bbA^{\kappa_0}_\phi$
depends linearly, we first calculate the individual terms for $\phi^m_j$:
\begin{align*}
\bbA^{\kappa_0}_{\phi^m_j} &= (m{-}j)(m{-}j{-}1)e^{m-j-2}v^{2j} \bbV^{\kappa_0} +
  \big((m{-}j)4j e^{m-j-1}v^{2j} +2j(2j{-}1) e^{m-j}v^{2j-2}\big) \bbE 
 \\
& \geq -\kappa_*(m{-}j)(m{-}j{-}1)e^{m-j-2}v^{2j+2}\bbE + j(4m{-}2j{-}1)e^{m-j-1}v^{2j} \bbE,
\end{align*} 
where we used $e\geq \frac12 v^2$ for the last term in the upper line. 
Summing the estimates over $j$ and using $b_j\geq 0$ we obtain the estimate  
$\bbA^{\kappa_0}_{\varphi_m} \geq b_*(v,e) \bbE $ with 
\begin{align*}
b_*(v,e)&= \sum_{j=0}^m b_j \Big(j(4m{-}2j{-}1)e^{m-j-1}v^{2j} -
-\kappa_*(m{-}j)(m{-}j{-}1)e^{m-j-2}v^{2j+2} \Big)  
\\
&= \sum_{i=1}^{m} \big( i(4m{-}2i{-}1) b_i - \kappa_*(m{-}i{+}1)(m{-}i)b_{i-1}
\big)  e^{m-i-1}v^{2i}.  
\end{align*} 
Clearly, $\bbA^{\kappa_0}_{\varphi_m}(v,w)\geq 0$ for all $(v,w)$, if
$b_*(v,e)\geq 0$ for all $(v,e)$, and this is easily reach by starting with
$b_0=1$ and then choosing iteratively $b_i$ via $ i(4m{-}2i{-}1) b_i =
\kappa_*(m{-}i{+}1)(m{-}i)b_{i-1}$ for $i=1,\ldots, m$. 
\end{proof} 

The following result follows easily by applying the previous lemma with
\eqref{eq:varphi.R} and the fact that $\varphi_m$ can be estimated from above
and below by $v^{2m}+w^m$.

\begin{proposition}[$\rmL^p$ bounds for $(v,w)$]
\label{pr:LpBounds}
For the case $\kappa \equiv\kappa_0\eta$ and $p=m\in \N$, there exists
a constant $C(\kappa_0,p)$ such that all smooth solutions $(v,w)$ of
\eqref{eq:SM01} satisfy 
\[
\big\|v(t) \big\|_{\rmL^{2p}(\Omega)}^2 + \big\| w(t) \big\|_{\rmL^p(\Omega)} 
\leq  C(\kappa_0,p) \big( \big\|\INIv \big\|_{\rmL^{2p}(\Omega)}^2 + \big\| \INIw
\big\|_{\rmL^p(\Omega)} \big)  \text{ for all } t \in [0,T].
\]
\end{proposition} 

\subsection{Estimates based on comparison principles}
\label{su:ComparisonEst}

As our system is given in terms of two scalar diffusion equations, we can apply
comparison princples when taking care of the interaction between the two
equations.

\STEP{(C1)} We first observe that $\INIw(x)\geq 0$ immediately implies
$w(t,x)\geq 0$ for all $t\geq 0$ and $x \in \Omega$. Of course this similarly
holds for $v$, but we do not need a sign condition for $u$.\smallskip 

\STEP{(C2)} Moreover, if $(v,w)$ is a smooth solution of our
system and $\underline w$ solves the scalar PME 
\[
\dot {\underline w} = \DIV\big(\kappa(\underline w) \nabla \underline w\big),
\quad \underline w(0,\cdot)=\INIw,
\]
then we have $w(t,x)\geq \underline w(t,x)$ for all $t\geq 0$ and $x \in
\Omega$, see \cite{Vazq07PMEM} for a proof. In particular, if $\INIw\geq c>0$
then $w(t,x)\geq c$ for all $(t,x)$.\smallskip  

\STEP{(C3)} If $v_* \leq \INIv(x) \leq v^*$ for all $x\in \Omega$, then $v(t,x)
\in [v_*,v^*]$ for  all $t\geq 0$ and $x \in \Omega$.\medskip

The next result is more advanced and truly uses the interaction of the two
equations. However,  it is restricted to the case $\kappa \equiv \eta$. Under this
assumption, equation \eqref{eq:varphi.PDE} gives 
\[
\pl_t\varphi = \DIV(\eta \nabla \varphi) - \eta \wh R \quad \text{with } \wh R= 
(\varphi_{vv}{-}\varphi_w)|\nabla v|^2+2\varphi_{vw}\nabla v{\cdot} \nabla w +
\varphi_{ww} |\nabla w|^2. 
\]
We are interested in the case $\wh R\geq 0$, for which initial conditions
 with $\varphi(\INIv,\INIw)\leq 0$ immediately implies $\varphi(v,w)\leq 0$ for
all $t \in [0,T]$ and $x \in \Omega$.

\begin{proposition}[Comparison for $\eta\equiv \kappa$] 
\label{pr:Com.eta=kappa}
Consider $a,b,c\in \R$ with $c\leq 2b$. If $(v,w)$ is a classical solution of
our coupled system \eqref{eq:SM01}, then we have
\[
a\INIv+b(\INIv)^2+c\INIw\leq 0 \text{ on }\Omega \quad \Longrightarrow \quad
av+bv^2+cw \leq 0 \text{ on }[0,T]\ti\Omega.
\]
In particular, $|\INIv|\leq M_* \INIw$ on $\Omega$ implies $|v(t,x)|\leq
M_*w(t,x)$ for $(t,x)\in [0,T]\ti \Omega$. 
\end{proposition}
\begin{proof}  We simply observe that the function $\varphi(v,w)=a v + bv^2
  +cw$ gives $\wh R(v,w)=(2b{-}c)|\nabla v|^2$ which is non-negative because of 
$c\leq 2b$. Then, the maximum principle applied to
$\zeta(t,x)=\varphi(v(t,x),w(t,x))$ shows that $\zeta_0\leq 0$ implies
$\zeta(t,x)\leq 0$, which is the desired result. 
 
The last assertion follows by choosing $(a,b,c)=(1,0,{-}M_*)$ for finding
$v\leq M_* w$ and by choosing $(a,b,c)=(-1,0,{-}M_*)$ for finding
$-v\leq M_*w$.
\end{proof}

\subsection{Dissipation estimates}
\label{su:DissipEstim}

Here we provide space-time estimates for various quantities, which will allow
us to derive suitable compactness for approximating sequences. 

We start with the estimate obtained in Section \ref{suu:IntegralEstim} for the
choice $\varphi(v,w)=\frac1p |v|^p$ with $p>1$. Integration in time gives the
estimate 
\begin{equation}
  \label{eq:DissEst.v.Lp}
  \iint_{Q_T}  \eta(w) (p{-}1) |v|^{p-2} \big| \nabla v\big|^2 \dd x \leq
\int_\Omega \frac1p |\INIv|^p\dd x .
\end{equation}
The special case $p=2$ shows that it is sufficient to use the initial energy,
namely 
\begin{equation}
  \label{eq:Diss.u.est}
\iint_{Q_T}  \eta(w)  \big| \nabla v\big|^2 \dd x \leq
\int_\Omega \frac12 |\INIv|^p\dd x \leq \calE(\INIv,\INIw).
\end{equation}
This estimate shows that the right-hand side in the $w$-equation in
\eqref{eq:SM01} is always in $\rmL^1(Q_T)$, but it will be difficult to
obtain higher integrability. \medskip

Choosing $\varphi(v,w)=\Phi(w)$ leads to the dissipation relation
\begin{align}
  \label{eq:Phi.DissRel}
\iint_{Q_T}\!\!\big(\kappa(w)\Phi''(w)|\nabla w|^2 - \eta(w) \Phi'(w) |\nabla
v|^2\big) \dd x \dd t 
  =   \int_\Omega\!\big\{\Phi(\INIw){-}\Phi(w(T))  \big\} \dd  x,
\end{align}
which is obtained by integrating \eqref{eq:Phi(k)}. Here we additionally have
to impose $w(t,x)>0$ in the case that $\Phi''(w)\to \infty$ for $w\to 0^+$. 

The last relation can be used in to two different ways to obtain a bound (i) on
$\nabla  w$ and (ii) on $\nabla v$.\medskip 

\STEP{(i) Estimates for $\nabla w$:} we observe that we may use
\eqref{eq:Diss.u.est} if $\Phi$ satisfies  
\begin{equation}
  \label{eq:Cond.Phi}
  0=\Phi(0) \leq \Phi(w) \leq C_\Phi w,\quad 0\leq \Phi'(w) \leq C_\Phi, \quad
\Phi''(k)\geq 0 
\end{equation}
for some $C_\Phi>0$. Then, we find 
\begin{align}
  \nonumber
\iint_{Q_T} \!\!\kappa(w)\Phi''(w)|\nabla w|^2 \dd x \dd t \leq  
C_\Phi\iint_{Q_T} \!\! \eta(w)|\nabla v|^2 \dd x \dd t +
  C_\Phi \int_\Omega \INIw \dd   x  \leq 2C_\Phi \calE(\INIv,\INIw).
\end{align}
In particular, one may want to make $\Phi''(w)$ big for $w\approx 0$ such that
$\kappa(w)\Phi''(w)$ is large there. Of course, we need 
integrability of $\Phi''$ on ${]0,\infty[}$ to allow for $0\leq \Phi'(w)\leq C_*$.  
Hence, generalizing  \cite{BocGal89NEPE,BoDaOr97NPEM,MieNau18?EGTW} we will use
the family  
\begin{equation}
  \label{eq:Phi.delta}
  \Phi(w)=\frac1{w^\delta(1{+}w)} \ \text{ for } w>0, \quad \text{where }
  \delta \in {]0,1[}. 
\end{equation}
The case $0<\delta \ll 1$ provides good results for large $w$, while
$0<1{-}\delta \ll 1$ produces good results for small $w$:  
Thus, for all $\delta \in
{]0,1[}$ there exists $C(\delta)>0$ such that for all classical solutions
$(v,w)$ with $w(t,x)>0$ on $[0,T]\ti \ol\Omega$ we obtain the estimate
\begin{equation}
  \label{eq:EstNablaW.delta}
\forall\,\delta \in {]0,1[}\ \exists\, C_\delta>0:\quad 
\iint_{Q_T} \frac{\kappa(w)}{w^{\delta}(1{+}w)}
\,\big|\nabla w \big|^2 \dd x \dd t \leq  
C_\delta\,  \calE(\INIv,\INIw).  \medskip 
\end{equation}

Alternatively, we may also consider a function $\varphi(v,w)=\Psi(w)$ such that 
\begin{equation}
  \label{eq:Cond.Psi}
  0 \geq \Psi(w) \geq -C_\Psi(1{+}w), \quad \Psi'(w)\leq 0, \quad \Psi''(w)>0
\quad \text{for all }w>0,
\end{equation}
where $\Psi(w)=- w^\gamma$ with $\gamma \in {]0,1[}$ is a typical candidate. We
can then drop the term $-\eta(w) \Psi'(w)|\nabla v|^2\geq 0$ and take advantage
of the fact that $\Psi''(w)$ can be much bigger than $\Phi'(w)$ for $\Phi$
satisfying \eqref{eq:Cond.Phi}. 
The estimates for $\Psi$ and the identity
\eqref{eq:Phi.DissRel} show that that positive classical solutions satisfy
\begin{align}
  \nonumber
   & \iint_{Q_T} \!\! \Big( \Psi''(w)\kappa(w)\big|\nabla w \big|^2
  -\Psi'(w)\eta(w)\big|\nabla v \big|^2  \Big) \dd x \dd t = 
 \int_\Omega\!\big\{\Psi(\INIw) - \Psi(w(T)) \big\} \dd x
\\
  \label{eq:Diss.Psi}
& \qquad \leq \int_\Omega C_\Psi(1+w(T)) \dd x \leq C_\Psi \big(|\Omega| +
\calE(\INIv,\INIw)\big),  \medskip 
\end{align}
where we used that $\Omega\subset \R^d$ has a finite Lebesgue volume.

\STEP{(ii) Estimates for $\nabla v$:}
The best estimate for $\nabla v$ can be obtained if the function $w \mapsto
1/\eta(w)$ is integrable near $w=0$, which is certainly the case for
$\eta(w)=w^\alpha$ with $\alpha\in {]0,1[}$. In this case, we define the
negative function $\Psi_\eta$ via 
\[
\Psi_\eta(w):=-\int_0^w \frac1{\eta(\varkappa)} \dd \varkappa .
\]
The choice leads to the relation $\Psi_\eta'(w)=-1/\eta(w)$ and $\Psi''_\eta(w)
= \eta'(w)/\eta(w)^2\geq 0$, such that the conditions \eqref{eq:Cond.Psi} hold
and \eqref{eq:Diss.Psi} leads to the special case 
\begin{align}
    \label{eq:regul.u.1}  
  \iint_{Q_T} \!\! \Big(\big|\nabla v \big|^2 + 
  \frac{\eta'(w)\kappa(w)}{\eta(w)^2}\big|\nabla w \big|^2\Big) \dd x \dd t 
  \leq C_{\Psi_\eta} \big(|\Omega| + \calE(\INIv,\INIw)\big),  \medskip 
\end{align}
see \cite[Eqn.\,(4.6)]{Naum13DPPT} for an earlier occurrence for the case
$\alpha=1/2$.  

Another way of estimating $\nabla u$ can be derived from \eqref{eq:Diss.u.est}
if we additionally have a pointwise estimate of the form 
$|v(t,x)| \leq M_* w(t,x)$, which was derived in Proposition
\ref{pr:Com.eta=kappa} for the special case $\kappa \equiv \eta$. Using the
monotonicity $\eta'(w)\geq 0$ we have $\eta(|v|/M_*)\leq \eta(w)$ and conclude
\begin{equation}
  \label{eq:regul.u.II}
\begin{aligned}  \iint_{Q_T} \!\! \eta\big( \frac1{M_*} |v|\big) |\nabla
  v|^2 \dd x \dd t \leq \calE(\INIv,\INIw) .
\end{aligned}
\end{equation}

We summarize the dissipation estimate by restricting to the homogeneous case
$\eta(w)=\eta_0w^\alpha$ and $\kappa(w)=\kappa_0 w^\beta$. More general cases 
can easily be deduced in the same fashion. We explicitly show the dependence on
$|\Omega|$ such that unbounded domains like $\Omega=\R^d$ can be treated as
well in cases where no dependence on $|\Omega|$ is indicated. 

\begin{proposition}[Dissipation estimates] 
\label{pr:DissipEstim} Consider the case $\eta(w)=\eta_0w^\alpha$ and
$\kappa(w)=\kappa_0 w^\beta$ with $\alpha,\beta, \eta_0,\kappa_0>0$ and bounded
$\Omega\subset \R^d$. Then for
all $\delta \in {]0,1[} $ there exists a constant $C_*>0$ such that smooth solutions
$(v,w)$ of \eqref{eq:SM01} with $w(t,x)\geq \underline w>0$ satisfy 
\begin{subequations}
  \label{eq:DissEstAll}
  \begin{align}
      \label{eq:DissEst.a}
 \beta>0:\quad  &\iint_{Q_T} \frac1{1{+}w}\Big|\nabla\big( 
   \,w^{1+(\beta-\delta)/2} \big)  \Big|^2 \dd x \dd t \leq C_* \calE(\INIv,\INIw)   , 
      \\
      \label{eq:DissEst.b}
 \beta>0:\quad  &\iint_{Q_T} \Big|\nabla\big( 
   w^{(\beta+\delta)/2} \big)  \Big|^2 \dd x \dd t 
      \leq C_*\big(|\Omega|+ \calE(\INIv,\INIw) \big)  ,
      \\
      \label{eq:DissEst.c}
 \alpha \in {]0,1[}:\quad  &\iint_{Q_T} \big|\nabla v   
   \big|^2 \dd x \dd t  \leq C_*\big(|\Omega|+ \calE(\INIv,\INIw) \big) .
  \end{align}
If we additionally have $\kappa \equiv \eta$ (i.e.\ $\alpha=\beta$ and
$\kappa_0=\eta_0$) and $|\INIv(x)| \leq M_* \INIw(x)$, then for $C_*$ now depending
on $\delta$ and $M_*$ only, we further have 
\begin{align}
      \label{eq:DissEst.d}
 \alpha=\beta:\quad & \iint_{Q_T} \Big|\nabla\big( 
   v^{1+\alpha/2} \big)  \Big|^2 \dd x \dd t \leq C_* \calE(\INIv,\INIw), 
\\   
 \label{eq:DissEst.e}
 \alpha=\beta \geq 1{-}\delta:\quad & \iint_{Q_T} \Big|\nabla\big( 
   v^{(1+\alpha+\delta)/2} \big)  \Big|^2 \dd x \dd t \leq C_* \big(|\Omega|+
   \calE(\INIv,\INIw) \big) .  
\end{align}
\end{subequations}
Throughout, the constant $C_*$ is independent of the lower bound $\ol w$. 
\end{proposition}  
\begin{proof} The results are obtained by applying the above estimate for
  suitable functions $\Phi$ and $\Psi$. Throughout use $w^\gamma |\nabla
  w|^2 = c_\gamma \big| \nabla(w^{1+\gamma/2})\big|^2$ for the smooth solutions. 

Estimate \eqref{eq:DissEst.a} follows directly from
\eqref{eq:EstNablaW.delta}. 

For \eqref{eq:DissEst.b} we exploit
\eqref{eq:Diss.Psi} with $\Psi(w)=-w^\delta$. 

Estimate  \eqref{eq:DissEst.c} is a consequence of \eqref{eq:regul.u.1},
because $1/w^\alpha$ is integrable for $\alpha <1$. 

For \eqref{eq:DissEst.d} we simply use \eqref{eq:regul.u.II}.  

Finally, \eqref{eq:DissEst.e} follows by combining \eqref{eq:Diss.Psi} with
$\Psi(w)=-w^\delta $, the estimate $|v|\leq M_* w$  from Proposition
\ref{pr:Com.eta=kappa}, and the monotonicity of $w \mapsto |\Psi'(w)| \eta(w)=
c w^{\alpha -1+\delta}$. 
\end{proof}

Based on the conjectures concerning the typical front behavior, which were
discussed in Section \ref{su:Conj.GrowSppt}, we see that $\nabla v \in
\rmL^2(Q_T)$ can be expected only in the case $\alpha<1$, i.e.\ the restriction
$\alpha \in {]0,1[}$ in \eqref{eq:DissEst.c} seems to be sharp. Similarly,
Section \ref{su:Conj.GrowSppt} show that $\nabla w^\gamma \in \rmL^2(Q_T)$
implies $\gamma >\beta/2$, which corresponds to the restriction $\delta>0$ in
\eqref{eq:DissEst.b}. 

\begin{remark}[$\int_\Omega \log \INIw \dd x$ finite]
  \label{re:BuliMalek} In their treatment of Kolmogorov's two-equation model in
  \cite{BulMal19LDAK} the authors consider the the situation corresponding to
  $\eta(w)=\eta_0 w$ and $\kappa(w)=\kappa_0 w$, i.e.\
  $\alpha=\beta=1$. Moreover, they assume that $\INIw$ is positive almost
  everywhere such that $L_0:= -\int_\Omega \log \INIw \dd x <\infty$. Choosing
  the convex functional $\Psi(w)=- \log w$ the dissipation relation
  \eqref{eq:Phi.DissRel} leads to
\begin{align*}
\iint_{Q_T}\!\! \Big(\frac{\kappa_0}{w}|\nabla w|^2 + \eta_0 |\nabla v|^2 \Big) \dd
x \dd t &= \int_\Omega \! \big(\log w(T)- \log \INIw\big) \dd x 
\\
&  \leq \int_\Omega w(T)\dd x  + L_0  \leq  \calE(\INIv,\INIw) + L_0. 
\end{align*}
Of course, the assumption $L_0<\infty$ does not allow to study solutions
with nontrivial support, i.e.\ it implies $\mafo{sppt}(w(t))=\ol\Omega$ for all
$t\geq 0$. 
\end{remark} 

\section{Existence of solutions}
\label{se:ExistTheory}

Through this section we use the following standard assumptions:
\begin{subequations}
  \label{eq:StandAssum}
\begin{align}
& \Omega \subset \R^d \text{ is bounded with $\rmC^2$ boundary}, 
\\
&\exists\, \alpha,\beta, \eta_0,\kappa_0>0 \ \forall\, w\geq 0:\quad 
\eta(w)= \eta_0 w^\alpha \text{ and }\kappa (w) = \kappa_0 w^\beta .
\end{align} 
\end{subequations}

\subsection{Three prototypical existence results}
\label{su:Prepare}

For given initial conditions $(\INIv,\INIw) \in \rmL^2(\Omega)\ti \rmL^1_\geq
(\Omega)$ we choose a sequence of smooth initial 
data $(\INIv_\eps,\INIw_\eps)\in \rmC^\infty(\ol\Omega)^2$ such that
\begin{subequations}
  \label{eq:IniConds.eps}
\begin{align}
  \label{eq:IniConds.eps.a}
  &\INIw_\eps(x)\geq \eps>0 \text{ on }\Omega,  
\\
  \label{eq:IniConds.eps.b}
 & (\INIv_\eps,\INIw_\eps)
  \to (\INIv,\INIw) \in \rmL^2(\Omega)\ti \rmL^1 (\Omega),
\\
  \label{eq:IniConds.eps.c}
 &\INIv \in \rmL^{p_0}(\Omega) \quad \text{and} \quad \INIv_\eps \to \INIv
 \text{ in } \rmL^{p_0}(\Omega). 
\end{align}
\end{subequations}
Condition \eqref{eq:IniConds.eps.b} implies the convergence of the conservation laws 
$\calV(\INIv_\eps)\to \calV(\INIv)$ and $\calE(\INIv_\eps,\INIw_\eps)\to
\calE(\INIv, \INIw)$. Condition \eqref{eq:IniConds.eps.c} with $p_0 \in [1,2]$  
follows from \eqref{eq:IniConds.eps.b}, hence, it will only be an extra
condition for $p_0>2$.

With this we have a classical parabolic system with positive viscosity such
that local existence of solutions $(v_\eps,w_\eps)$ is classical, see e.g.\
\cite{Aman89DTQP3, Wieg92GSCS} or \cite[Cha.\,8]{Luna95ASOR}. Of course all
these solutions satisfy our a priori bounds, in particular, we have
$w_\eps(t,x)\geq \eps$ as long as the solution exists, cf.\ Section
\ref{su:ComparisonEst} (C2).  This implies that the solution stays as
smooth. Moreover, the structure of the equation for $v$ implies the global
$\rmL^\infty$ bound $\|v_\eps(t)\|_{\rmL^\infty(\Omega)} \leq
\|\INIv_\eps\|_{\rmL^\infty(\Omega)}$, see \eqref{eq:LpNormDecayV}. Finally,
the energy conservation $\calE(v_\eps(t),w_\eps(t)) = 
\calE(\INIv_\eps,\INIw_\eps)$ which also implies
$\iint_{Q_T} \eta(w_\eps)|\nabla v_\eps|^2 \dd x \dd t \leq
\calE(\INIv_\eps,\INIw_\eps) $. Thus, blow-up is not possible and the classical
solutions exist for all time.
  
The aim is now to show that the solutions $(v_\eps,w_\eps)$, or better a
suitable subsequence thereof, converge to a limit
$(v,w)$ that is a weak or very weak solution of our coupled system. The problem
here is that the limit $w$ may have a nontrivial support strictly contained in
$\Omega$. As a consequence, an integral bound 
\begin{equation}
  \label{eq:GenEstimate}
   \iint_{Q_T} \big(
w_\eps^{\gamma} |\nabla w_\eps^\alpha|^2 + w_\eps^{\delta} |\nabla v_\eps|^2\big)
\dd x \dd t \leq C_*
\end{equation} 
does not necessarily imply spatial compactness for
$v_\eps$.  Thus, we will provide two different existence results, in
the case $\eta(w)=\eta_0w^\alpha$ with $\alpha<1$ one may obtain the exponent
$\delta=0$ and a bound for $\nabla v_\eps $ in $\rmL^2(Q_T)$ follows. 
In the case, $\alpha\geq 1$, we can only treat the case $\eta\equiv \kappa$
under the additional assumption $|\INIv(x)|\leq M_* \INIw(x)$, which allows us
to obtain a bound on $\nabla (v_\eps)^{1+\delta/2}$  in $\rmL^2(Q_T)$.   

Another way of obtaining weak solutions with a positive lower bound $\eps>0$
would be to adapt the methods in \cite{Naum13DPPT, BulMal19LDAK,
  MieNau18?EGTW} to the present case. Our approach is different is similar to
the one taken in \cite{BerKam90SDPE, DalGia99WSSC}; however, our goal is more
restrictive.  There, for the solutions
$(\rho,\theta)$ of the plasma model \eqref{eq:Rosenau} one only asks for
$\rho^{\delta/2}\nabla \theta \in \rmL^2(\calP)$ where $\calP:=\bigset{(t,x)\in
[0,T]\ti \Omega}{ \rho(t,x)>0}$.\medskip

Throughout the rest of this section we use the following short-hand notations:
\begin{align*}
\text{(a) \ }&u_\eps \bbin X \quad \Longleftrightarrow \quad 
\exists\, C>0 \ \forall \, \eps\in{]0,1[}: \quad u_\eps \in X \text{ and } \|
u_\eps\|_X \leq C,
\\[0.5em]
\text{(b) \ }& \rmL^p\rmL^q := \rmL^p\big([0,T];\rmL^q(\Omega)\big) \quad \text{and similarly }
\rmL^q\rmH^1  \text{  or } \rmL^p\rmW^{s,q}. 
\end{align*}
We will use the following standard interpolation in $\rmL^s\rmL^p$ (cf.\ 
\cite[Lem.\,4.2]{MieNau18?EGTW}):
\begin{equation}
  \label{eq:Interpol.LsLp}
\begin{aligned}
&\forall\, \theta\in [0,1],\ s,s_1,s_2,p,p_1,p_2\in [1,\infty] \text{ with }
\textstyle \frac1s= \frac{1{-}\theta}{s_1} + \frac\theta{s_2} \text{ and } 
\frac1p= \frac{1{-}\theta}{p_1} + \frac\theta{p_2}  
\\ 
& \exists \,C>0\ \forall\, u\in \rmL^{s_1}\rmL^{p_1} \cap \rmL^{s_2}\rmL^{p_2}: 
\quad
\| u\|_{\rmL^s \rmL^p} \leq C \| u\|_{\rmL^{s_1} \rmL^{p_1}}^{1-\theta} 
 \| u\|_{\rmL^{s_2} \rmL^{p_2}}^\theta.  
\end{aligned}  
\end{equation}
Moreover, we define $2_d^*$ as the optimal exponent in the embedding
$\rmH^1(\Omega) \subset \rmL^{2_d^*}(\Omega) $, i.e.\ $2_d^*=2d/(d{-}2)$ for
$d\geq 3$, $2_1^*=\infty$, and  $2_2^*<\infty$.

We are now ready to state three different existence results. The first result
provides weak solutions and relies on the restriction \eqref{eq:p0.cond}
concerning the exponents 
$\alpha$ and $\beta$ as well as the integrability power $p_0$ for $\INIv$. We
do not expect that condition \eqref{eq:p0.cond} is sharp.

\begin{theorem}[Existence of weak solutions for $\alpha<1$] 
\label{th:ExistWeakSol} Assume that \eqref{eq:StandAssum} holds with $\alpha
\in {]0,1[}$.  
Moreover, consider $p_0\geq 2$ such that
\begin{equation}
  \label{eq:p0.cond}
  \frac{\alpha}{\beta + 2 - 2/2_d^*} + \frac1{p_0+2- 2p_0/2_d^*} < \frac12.
\end{equation}
Then, for all initial data $(\INIv,\INIw) \in
\rmL^{p_0}(\Omega) \ti \rmL^1_\geq (\Omega)$ there exists a weak solution $(v,w)$
in the sense of Definition \ref{de:WeakVeryWeak}. Moreover, this solution
satisfies 
\[
\nabla v\in \rmL^2(Q_T) \quad \text{ and } \quad w \in
\rmL^{\beta+2-2/2_d^*}(Q_T). 
\]
\end{theorem}

Condition \eqref{eq:p0.cond} is somehow restrictive and shows that small
$\alpha$ and large $\beta$ are desirable. Even assuming $\INIv \in
\rmC(\ol\Omega)\subset \rmL^\infty(\Omega) $, i.e.\ $p_0=\infty$, we still have
to satisfy 
\[
2\alpha < \beta + 2 - \max\{ 0, 1{-}2/d\}. 
\]   
Because of $\alpha\in {]0,1[}$ this is always satisfied for $d\leq 2$, but may
provide a nontrivial lower bound for $d\geq 3$ and $\alpha >1/2$. However, the
important case $\alpha=\beta\in {]0,1[}$ is always possible, even for finite
$p_0\in {[2,\infty[}$. 

For example, assuming $\alpha=\beta=1/2$ (as in Kolmogorov's one-equation or
Prandtl's model discussed in Section \ref{suu:Kolmo1eqn}) we obtain, after some
computation, the restriction
\[
p_0 \geq 2 \quad \text{ and } \quad p_0> \frac{2d^2}{d{+}4}. 
\]

For the next result we stay in the case $\alpha\in {]0,1[}$ and show that very
weak solutions can be obtained in a previously inaccessible regime. However,
the usage of the weighted gradients as introduced in Definition
\ref{de:WeightedGrad} forces us to control $ \nabla w_\eps^\alpha$ in $\rmL^1
(Q_T)$ whereas the equation provides control on $\nabla w_\eps^{\beta/2}$
only. Thus, the restriction $\beta<2\alpha$ will be needed.  

\begin{theorem}[Existence of very weak solutions for $\alpha<1$] 
\label{th:ExistVWS.al} Assume that \eqref{eq:StandAssum} holds with 
$0 < \beta < 2\alpha <2$.  
Moreover, consider $p_0\geq 2$ such that
\begin{equation}
  \label{eq:p0.cond.VWS}
  \frac{2\alpha+1 -2/2_d^*}{2(\beta {+} 2 {-} 2/2_d^*)} + \frac1{p_0+2- 2p_0/2_d^*} < 1.
\end{equation}
Then, for all initial data $(\INIv,\INIw) \in
\rmL^{p_0}(\Omega) \ti \rmL^1_\geq (\Omega)$ there exists a very weak solution $(v,w)$
in the sense of Definition \ref{de:WeakVeryWeak}. Moreover, this solution
satisfies 
\[
\nabla v\in \rmL^2(Q_T) \quad \text{ and } \quad w \in
\rmL^{\beta+2-2/2_d^*}(Q_T). 
\]
\end{theorem}

We emphasize that condition \eqref{eq:p0.cond.VWS} is weaker than
\eqref{eq:p0.cond}. Indeed, consider the case $p_0=\infty$ for simplicity, then 
\eqref{eq:p0.cond.VWS} is equivalent to $2\alpha < 2\beta + 3-2/2_d^*$, which
is automatically satisfied for $0\leq \beta < 2\alpha <2$, because of $2\alpha
\leq \alpha + 1 < 2\beta +1 < 2\beta +3-2/2_d^*$. 

However, for $p_0=\infty$
\eqref{eq:p0.cond} reduces to $2\alpha < \beta +2-2/2_d^*$ which is violated
for $d\geq 3$ and suitable $\alpha $ and $\beta$, e.g.\ for
$(d,\alpha,\beta)=(4,7/8,1/8)$ we have 
\[
0< \beta=1/8 < 2\alpha=7/4 < 2, \quad \text{and} \quad 2\alpha = 7/4 \not< 13/8
= 2 {+} 1/8 {-} 1/2 = 2 {+} \beta  {-} 2/2_d^*.
\]  
Thus, there are cases where Theorem \ref{th:ExistVWS.al} can be applied but
Theorem \ref{th:ExistWeakSol} cannot. 

Finally we treat the case $\eta \equiv \kappa$ which is special in two ways:
first the energy density $e= \frac12v^2 +w$ satisfies the simple equation $\dot
e = \DIV\!\big( \eta(w)\nabla e\big)$, and second we can pointwise bound $v$ by
$w$, see Proposition \ref{pr:Com.eta=kappa}. The case $\kappa \equiv  \eta=
w^\alpha$ is very special, but it is presently the only case where $\alpha\geq
1$ can be handled. 

\begin{theorem}[Existence of very weak solutions for $\eta\equiv \kappa$]
\label{th:Exi.eta=kappa}
Assume that \eqref{eq:StandAssum} holds with $\eta(w)= \kappa(w)=w^\alpha$ with
$\alpha >0$.  Moreover, consider $p_0\geq 2$ such that
\begin{equation}
  \label{eq:p0.cond.e=k}
  p_0 \Big(1- \frac2{2_d^*} \Big) > \sqrt{8+4\alpha + \alpha^2/4}\: -\frac\alpha2.
\end{equation}
Then, for all $M_*>0$ and all initial data $(\INIv,\INIw) \in
\rmL^{p_0}(\Omega) \ti \rmL^{p_0/2}_\geq (\Omega)$ satisfying 
\begin{equation}
  \label{eq:v.leq.Mw}
  |\INIv(x)| \leq M_* \INIw(x) \ \text{ for a.a. } x\in \Omega
\end{equation}
there exists a very weak solution $(v,w)$
in the sense of Definition \ref{de:WeakVeryWeak}. Moreover, this solution
satisfies 
\[
 |v(t,x)| \leq M_* w(t,x) \ \text{ for a.a. } (t,x)\in Q_T.
\]
\end{theorem} 

The right-hand side in \eqref{eq:p0.cond.e=k} is strictly increasing with
range ${]2,4[}$.  Hence, the assumption $p_0\geq \max\{4,2d\}$ is sufficient for all
$\alpha>0$.
  
The proofs of these three results are the contents of the following three
subsections.

\subsection{Limit passage to weak solutions in the case $\alpha<1$} 
\label{su:LimPasAlpha1}

Here we provide the proof of existence of weak solutions first, and then show
what has to be changed for obtaining very weak solutions.\medskip

\noindent
\begin{proof}[Proof of Theorem \ref{th:ExistWeakSol}] \!\!We proceed in three
  steps: 1. a priori estimates, 2. compactness, and 3. identification of nonlinear
  limits. 

\STEP{Step 1: A priori estimates.}\\ Because of $\calE(v_\eps(t),w_\eps(t))=E_0$
and \eqref{eq:DissEst.c} we obviously have $v_\eps \bbin \rmL^2\rmH^1$.

For $w_\eps$ we use \eqref{eq:DissEst.b} to obtain
$ \nabla w_\eps^{(\beta+\delta)/2} \bbin \rmL^2\rmL^2$ for all
$\delta \in {]0,1[}$.  Since energy conservation implies
$w_\eps \bbin \rmL^\infty\rmL^1$, we conclude
$w_\eps^{(\beta+\delta)/2} \bbin \rmL^2\rmH^1$. For bounded Lipschitz domains
$\Omega \subset \R^d$ we have the embedding
$\rmH^1(\Omega) \subset \rmL^{2_d^*}(\Omega) $ with $2_1^*=\infty$,
$2_2^*<\infty$, and $2_d^*=2d/(d{-}2)$ for $d\geq 3$. We obtain
$w_\eps^{(\beta+\delta)/2} \bbin \rmL^2\rmL^{2_d^*}$, or equivalently
$w_\eps \bbin \rmL^{\beta+\delta}\rmL^q$ for $q= (\beta{+}\delta)2_d^*/2$.
Using the interpolation \eqref{eq:Interpol.LsLp} with
$w_\eps \bbin \rmL^\infty\rmL^1$ gives
\begin{equation}
  \label{eq:Interpol.weps}
  w_\eps \bbin \rmL^p(Q_T)=\rmL^p\rmL^p \text{ with } p = p(\delta):= \beta
  +\delta +1 -2/2_d^*.   
\end{equation}

For controlling the gradient of $\Pi(w)=\frac1{\beta+1} w^{\beta+1}$ in
$\rmL^{1+\eta}(Q_T)=\rmL^{1+\eta}\rmL^{1+\eta}$ for some $\eta>0$ we use the
simple identity
\[
\nabla (\Pi(w))=w^\beta \nabla w = \frac1{\gamma_2} w^{\gamma_1} \nabla \big(
w^{\gamma_2}\big) \quad \text{for }\gamma_1+\gamma_2=\beta +1.
\]
Choosing $\gamma_2=(\beta{+}\delta)/2$ and using $ \nabla
w_\eps^{(\beta+\delta)/2}  \bbin \rmL^2\rmL^2$ we find 
$\nabla \Pi(w_\eps) \bbin \rmL^{1+\eta}(Q_T)$ if we have  $w_\eps^{\gamma_1}
\bbin \rmL^q\rmL^q$ or $w_\eps \bbin 
\rmL^{\gamma_1 q}\rmL^{\gamma_1 q}$ with $1/2 +1/q <1$. Hence, it suffices to
show $p(\delta)>2\gamma_1(\delta)=:\beta +2 -\delta$, which holds for  for some
$\delta \in {]0,1[}$ because of the continuity of 
$\delta\mapsto (\gamma_2(\delta),p(\delta))$ on $[0,1]$ and $p(1)= \beta
+2-2/2_d^*> 2\gamma_2(1)= \beta+1$. In summary, we
have established 
\begin{equation}
  \label{eq:nabla.pi}
  \exists\, \eta >0 :\qquad \nabla \big(\Pi(w_\eps)\big)\ \bbin \
  \rmL^{1+\eta}(Q_T). 
\end{equation}

Clearly, the classical dissipation estimate \eqref{eq:Diss.u.est} for the
$v$-equation gives
\begin{align}
\label{eq:SQRTeta.n.veps}
U_\eps:= \eta(w_\eps)^{1/2} \nabla v_\eps = w_\eps^{\alpha/2} \nabla v_\eps 
\ \bbin \ \rmL^2(Q_T) 
\end{align}
Similarly, using  $w_\eps^{\alpha/2} \bbin \rmL^\infty
\rmL^{2/\alpha} \subset \rmL^{2/\alpha}(Q_T)$ and $\eta(w_\eps)\nabla v_\eps =
w_\eps^{\alpha/2} \,w_\eps^{\alpha/2}\nabla v_\eps$ yields 
\begin{align}
\label{eq:eta.nabla.veps}
 \eta(w_\eps) \nabla v_\eps = w_\eps^\alpha \nabla v_\eps \ \bbin \
 \rmL^{2/(1{+}\alpha)}(Q_T).  
\end{align}

%%% sssssssssssssssssssssssssssssssssssssssss

\STEP{Step 2: Compactness.}\\
 As usual we apply an Aubin-Lions-Simon theorem, see 
\cite[Cor.\,4]{Simo87CSSL} , \cite[Thm.\,5.1]{Lion69QMRP},
\cite[Lem.\,7.7]{Roub13NPDE}  and for nonlinear versions see 
\cite[Thm.\,1]{Mous16VCAL}.  

For this we first create \emph{spatial compactness},
which is trivial for $v_\eps$ because of $v_\eps \bbin \rmL^2\rmH^1$ with
$\rmH^1=\rmH^1(\Omega)$ compactly embedded into $\rmL^2(\Omega)$.  

For $w_\eps$ we have the difficulty that we only control the derivative of
power of $w$ in \eqref{eq:nabla.pi}, namely
$w_\eps^{1+\beta} \bbin \rmL^{1+\eta}\rmW^{1,1+\eta}$. However, employing the
fact that the Nemitskii operator $u\mapsto u^{\theta}$ maps
$\rmW^{s,p}(\Omega)$ into $\rmW^{\theta s, p/\theta}(\Omega)$ for all
$s\in {]0,1[}$ and $p \in {]1,\infty[}$, see Proposition
\ref{pr:FractPow}. This leads to
\[
w_\eps \bbin \rmL^{(1{+}\eta)(1+\beta)} \rmW^{\sigma, (1{+}\eta)(1+\beta)}
\quad\text{for all } \sigma \in {\big[ 0,\tfrac1{1{+}\beta}\big[}.
\]
Clearly, $\rmW^{\sigma, (1{+}\eta)(1+\beta)}(\Omega)$ is still compactly
embedded into $\rmL^{1{+}\eta}(\Omega)$. 

To derive temporal compactness we use the PDEs for $\dot v_\eps$ and  $\dot
w_\eps$, see \eqref{eq:SM01}. From \eqref{eq:SM01.a} and
\eqref{eq:eta.nabla.veps} we obtain 
\begin{equation}
  \label{eq:dot.veps}
  \dot v_\eps \in \rmL^{2/(1{+}\alpha)}
  \rmW_0^{-1,2/(1{+}\alpha)}=\rmL^{2/(1{+}\alpha)} \big(
  \rmW^{1,2/(1{-}\alpha)}(\Omega)\big)^*  . 
\end{equation}
Starting from $\nabla \Pi(w_\eps) \bbin
\rmL^{1+\eta}(Q_T)$ and $\eta(w_\eps)|\nabla v_\eps|^2 \bbin \rmL^1(Q_T)$ 
the PDE \eqref{eq:SM01.b} for $w_\eps$ gives 
\begin{equation}
  \label{eq:dot.weps}
  \dot w_\eps \in \rmL^{1+\eta} \rmW_0^{-1,1+\eta} + \rmL^1 \rmL^1  \ \subset \
  \rmL^1 \big(\rmW^{1,1+1/\eta}(\Omega)\big)^*. 
\end{equation}
Now applying Banach's selection principle for weak convergence and 
Aubin-Lions-Simon theorem strong convergence we obtain a limit pair $(v,w) $ with 
\[
 v \in \rmC_\rmw([0,T];\rmL^2(\Omega)) \cap \rmL^2\rmH^1
  \quad \text{and} \quad  w \in \rmL^\infty\rmL^1 \cap  \rmL^{(1{+}\eta)(1+\beta)}
  \rmW^{\sigma, (1{+}\eta)(1+\beta)} 
\]
and, along a suitable subsequence, the convergences 
\begin{subequations}
  \label{eq:Conv.ve.we}
\begin{align}
  \label{eq:Conv.ve.we.a}
(v_\eps(x), w_\eps(x)) &\to (v(x),w(x))
    &\text{a.e.\ in }&  \Omega;
\\
  \label{eq:Conv.ve.we.b}
(v_\eps,w_\eps) &\to (v,w)  &\text{in } &\rmL^2(Q_T)\ti \rmL^{1+\eta}(Q_T);
\\
  \label{eq:Conv.ve.we.c}
v_\eps &\weak v & \text{ in } &\rmL^2 \rmH^1,
\\
  \label{eq:Conv.ve.we.d}
w_\eps &\weak w & \text{in }&\rmL^{1+\beta}\rmW^{1/(2+\beta), 1+\beta}.
\end{align}
\end{subequations}

\STEP{Step 3: Identification of nonlinear limits.}\\
It remains to pass to the limit in the nonlinear terms of the weak formulation,
namely in 
\[
\text{(i) } \eta(w_\eps)\nabla v_\eps, \qquad \text{(ii) } \nabla \Pi(w_\eps) ,
\qquad \text{(iii) } \eta(w_\eps)v_\eps \nabla v_\eps. 
\]
In the trivial nonlinear term $\iint_{Q_T} v_\eps^2\pl_t \xi \dd x \dd t $ in the
left-hand side of \eqref{eq:DefWeakS.d} , the limit passage follows easily from
the strong convergence of $v_\eps$ in $\rmL^2(Q_T)$. 

Concerning the term (i) we recall the decomposition
\[
w_\eps^\alpha \nabla v_\eps = w_\eps^{\alpha/2} U_\eps \quad \text{with }
U_\eps = w_\eps^{\alpha/2} \nabla v_\eps. 
\]
By \eqref{eq:Conv.ve.we.b} we have $ w_\eps^{\alpha/2} \to w^{\alpha/2}$ in
$\rmL^{2/\alpha} $. One the one hand, we  have $U_\eps \weak U$ in $\rmL^2(Q_T)$ because of
\eqref{eq:SQRTeta.n.veps}. On the other hand, the strong convergence of
$w_\eps^{\alpha/2}$ and $\nabla v_\eps \weak \nabla v$ in $\rmL^2(Q_T)$ implies 
$U_\eps \weak  w^{\alpha/2}\nabla v$ in $\rmL^{2/(1{+}\alpha)}(Q_T)$. Thus,
$U=w^{\alpha/2}\nabla v$ and, with the same argument we find 
\begin{align}
  \label{eq:Cvg.(i)}
w_\eps^\alpha \nabla v_\eps = w_\eps^{\alpha/2} U_\eps \weak w^{\alpha/2}U =
w^\alpha \nabla v \quad \text{in }\rmL^{2/(1+\alpha)}.  
\end{align}

For term (ii) we use $\Pi(w)=\frac1{1+\beta} w^{1+\beta}$ and the a priori
estimates from Step 2, such that along a further subsequence (not relabeled) we
have, for all $q\in {[1,1{+}\beta 2_d^*/2[}$, 
\[
 w_\eps \weak w \text{ in } \rmL^q(Q_T), \quad 
w_\eps^{1+\beta} \weak g \text{ in } \rmL^{1+\eta}(Q_T), \quad 
\nabla w_\eps^{1+\beta} \weak G \text{ in } \rmL^{1+\eta}(Q_T). 
\]
Together with the strong convergence \eqref{eq:Conv.ve.we.b} we obtain $w_\eps
\to w$ in $\rmL^q(\Q_T)$ for the same $q$ such that $
w_\eps^{1+\beta} \to w^{1+\beta}$  in $\rmL^1(Q_T)$, which implies
$g=w^{1+\beta}$ as desired, and $G=\nabla g= \nabla \Pi(w)$ follows.

For the term (iii) we may need the extra condition $\INIv \in
\rmL^{p_0}(\Omega)$, which implies $v_\eps \bbin \rmL^\infty
\rmL^{p_0}$. Interpolation with $v_\eps \bbin \rmL^2\rmH^1 \subset
\rmL^2\rmL^{2_d^*}$ via \eqref{eq:Interpol.LsLp} we obtain 
\begin{align*}
&v_\eps \bbin \rmL^{q_0}(Q_T) \text{ with } q_0:= p_0+2 -2p_0/2_d^*,
%\\ &
\quad \text{and}\quad
v_\eps \to v  \text{ in } \rmL^q(Q_T)  \text{ for } q\in {[1,q_0[},
\end{align*} 
where the second statement uses \eqref{eq:Conv.ve.we.b}. Similarly, 
\eqref{eq:Interpol.weps} and \eqref{eq:Conv.ve.we.b} yield
\[
w_\eps^\alpha \to w^\alpha \text{ in } \rmL^{p/\alpha }(Q_T) \text{ for } p \in
{[1, \beta {+} 2 {-} 2/2_d^*[}.
\] 

Combining these two strong convergences and exploiting the condition
\eqref{eq:p0.cond}, we conclude $ \eta(w_\eps) v_\eps \to \eta(w)v$ in
$\rmL^2(Q_T)$. Together  with the weak convergence $\nabla v_\eps \to \nabla v$
(cf.\ \eqref{eq:Conv.ve.we.c}) we obtain the desired weak convergence for term
(iii), i.e.\ $ \eta(w_\eps) v_\eps \nabla v_\eps \to \eta(w)v\nabla v$ in
$\rmL^1(Q_T)$. Thus, the limit passage for all terms on the weak formulation
\eqref{eq:DefWeakSol}, which shows that $(v,w)$ is indeed a weak solution. 

This finishes the proof of Theorem \ref{th:ExistWeakSol}. 
\end{proof}

\subsection{Limit passage to very weak solutions for $\alpha<1$} 
\label{su:LimPasVWSAlpha1}

\noindent
\begin{proof}[Proof of Theorem \ref{th:ExistVWS.al}] \!\!We follow the same
  three steps as in the proof of Theorem \ref{th:ExistWeakSol}.  Step 1 and
  Step 2 are in fact identical, and in Step 3 we have to exploit the idea of
  weighted gradients for $\eta(w)\nabla v$, which means that we have to prove
  weak convergence in $\rmL^1(Q_T)$ for the following five terms
\[
\text{(i)' } v_\eps \nabla\eta(w_\eps), \quad \text{(ii) } \nabla \Pi(w_\eps) ,
\quad \text{(iii)' } v_\eps^2\nabla \eta(w_\eps), \quad \text{(i)'' }
v_\eps\eta(w_\eps), \quad \text{(iii)'' } v_\eps^2 \eta(w_\eps).
\]
Clearly, (ii) works as before and (i)' and (iii)' look as (i) and (iii) above,
but the gradient is moved from $v_\eps$ to $\eta(w_\eps)$, while in (i)'' and
(iii)'' there are no gradients at all. 

It is easy to see that the most critical term (iii)'. If the convergences work
for this case, then they also work for all the other ones. 

The problem is now that we need to control $\nabla w_\eps^\alpha$ which is
nontrivial if $\alpha$ is small, because \eqref{eq:DissEst.a} and
\eqref{eq:DissEst.b} provide bounds for $\nabla w_\eps^\gamma$ for $\gamma$
defined in terms of $\beta$. We hence assume 
\[
\beta < 2\alpha <2.
\]
Using $\nabla w^\alpha = \frac{\alpha}\gamma w^{\alpha-\gamma} \nabla w^\gamma$
for $\gamma \in {]0,\alpha[}$ and choosing $\gamma=(\beta{+}\delta)/2$ for some
$\delta \in {]0,1[}$ we obtain from \eqref{eq:DissEst.b} and \eqref{eq:Interpol.weps} 
the boundedness
\begin{align*}
\alpha \in {\big]\tfrac\beta2, \tfrac{\beta+1}2 \big[}:\quad 
 &\nabla w_\eps^\alpha \ \bbin \ \rmL^2(Q_T) 
\\
\alpha \in {\big[ \tfrac{\beta+1}2 ,1\big[}:\quad 
 &\nabla w_\eps^\alpha \ \bbin \ \rmL^r(Q_T) \ \text{ for all } r \in {[1,r_0[} 
\text{ with } r_0= \frac{2(\beta{+}2{-}2/2_d^*)}{ 2\alpha{+}1{-}2/2_d^*}. 
\end{align*}
Because of $r_0>1$ and the strong convergence \eqref{eq:Conv.ve.we.b} we 
obtain $\nabla w_\eps^\alpha \weak \nabla w^\alpha$ in $\rmL^r(Q_T)$. 

As before we have the strong convergence $v_\eps^2 \to v^2$ in
$\rmL^{q/2}(Q_T)$ for all $q\in {[1,q_0[}$ and now condition
\eqref{eq:p0.cond.VWS} is exactly made to allows us to conclude that $v_\eps^2
\nabla \eta(w_\eps) \weak v^2 \nabla \eta(w)$, as desired. 
This finishes the proof of Theorem \ref{th:ExistVWS.al}.   
\end{proof} 

\subsection{Limit passage in the case $\eta=\kappa$} 
\label{su:LimPasEta=Ka}

For $a\in \R$ and $\gamma>0$ we set $\{a\}^\gamma:= \mafo{sign}(a)\,|a|^\gamma$. 
\\[0.3em]
\begin{proof}[Proof of Theorem \ref{th:Exi.eta=kappa}] \!\!We proceed along the
same three steps as in Section \ref{su:LimPasAlpha1} but now use that it is
easy to construct smooth initial conditions $(\INIv_\eps,\INIw_\eps)$ such that
\eqref{eq:IniConds.eps} holds together with
$|\INIv_\eps(x) | \leq M_* \INIw_\eps(x)$ for $x\in \Omega$. For this one
simply mollifies the nonnegative functions $a^\pm :=M_*\INIw\pm \INIv$ to
obtain smooth and nonnegative $a^\pm_\eps$. Then, we set
$\INIv_\eps = (a^+_\eps{-}a^-_\eps)/2$ and
$ \INIw_\eps =(a^+_\eps{+}a^-_\eps)/(2M_*)$. With Proposition
\ref{pr:Com.eta=kappa} we obtain
\begin{equation}
    \label{eq:veps.leq.Mweps}
|v_\eps(t,x) | \leq M_*w_\eps(t,x) \quad \text{for all }(t,x) \in Q_T.  
\end{equation}

\STEP{Step 1: a priori bounds.} We first observe $e^0:=\frac12(\INIv)^2+\INIw
\in \rmL^{p_0/2}(\Omega)$, and $\dot e = \DIV(w^\alpha e)$ implies
$\|e(t)\|_{\rmL^{p_0/2}} \leq \|e^0\|_{\rmL^{p_0/2}}$. Thus, we have $
v_\eps  \bbin \rmL^\infty \rmL^{p_0}$ and $w_\eps  \bbin \rmL^\infty
\rmL^{p_0/2}$, and
 using \eqref{eq:DissEst.v.Lp} and \eqref{eq:veps.leq.Mweps} we obtain 
$\{v_\eps\}^{(p_0+\alpha)/2} \bbin \rmL^2 \rmH^1\leq \rmL^2\rmL^{2_d^*}$ and conclude 
$v_\eps \bbin \rmL^{p_0+\alpha}\rmL^{(p_0+\alpha) 2_d^*/2}$. Using $
v_\eps \bbin \rmL^\infty \rmL^{p_0}$ and the
interpolation \eqref{eq:Interpol.LsLp} we arrive at 
\begin{equation}
  \label{eq:v.est.e=k}
  v_\eps \bbin \rmL^{q_\rmv} (Q_T) \ \text{ with } q_\rmv := 
  \alpha + 2p_0\big(\tfrac12 -\tfrac1{2_d^*}\big).
\end{equation}
From \eqref{eq:DissEst.b} with $\alpha=\beta$ we find $w^{(\alpha+\delta)/2} \bbin
\rmL^2\rmH^1$ for all $\delta \in {]0,1[}$ and interpolation
with  $w_\eps \bbin \rmL^\infty \rmL^{p_0/2}$ yields
\begin{equation}
  \label{eq:w.est.e=k}
  w_\eps \bbin \rmL^{q} (Q_T) \ \text{for all }q\in {[1,q_\rmw[} 
\text{ with } q_\rmw := 1+\alpha+ p_0\big(\tfrac12 -\tfrac1{2_d^*}\big).
\end{equation}

\STEP{Step 2: compactness.} From $\{v_\eps\} ^{(p_0+\alpha)/2} \bbin \rmL^2 \rmH^1$
and $ w^{\alpha+\delta} \bbin \rmL^2\rmH^1$ we again obtain  
\[
v_\eps \bbin \rmL^{p_0+\alpha}\rmW^{2s/(p_0+\alpha), p_0+\alpha} \ 
\text{ and } \  w_\eps \bbin \rmL^{2(\alpha+\delta)}
\rmW^{s/(\alpha+\delta),2(\alpha+\delta)} \  \text{ for } s\in {]0,1[}
\text{ and } \delta \in {]1/2,1[},
\]
where we used Proposition \ref{pr:FractPow}.  
Thus, compactness works as in Step 2 of Section \ref{su:LimPasAlpha1}. In
particular, we have 
\begin{equation}
  \label{eq:stro.cvg.weps}
  v_\eps \to v \ \text{ in }\ \rmL^p(Q_T) \text{ for }p\in {[1,q_\rmv[}
\quad \text{and} \quad 
  w_\eps \to w \ \text{ in }\ \rmL^q(Q_T) \text{ for }q\in {[1,q_\rmw[} .
\end{equation}
Using the pointwise convergence a.e.\ in $Q_T$, the limits $(v,w)$ still
satisfies \eqref{eq:veps.leq.Mweps}. 

\STEP{Step 3: Identification of nonlinear limits.} We have to show that the
five terms 
\[
\text{(i)' } v_\eps \nabla\eta(w_\eps), \quad \text{(ii) } \nabla \Pi(w_\eps) ,
\quad \text{(iii)' } v_\eps^2\nabla \eta(w_\eps), \quad \text{(i)'' }
v_\eps\eta(w_\eps), \quad \text{(iii)'' } v_\eps^2 \eta(w_\eps)
\]
converge to their respective limits. 

Clearly, $w_\eps^{1+\alpha} \bbin \rmL^p(Q_T)$ for all $p\in
{[1,q_\rmw/(1{+}\alpha)[}$ and by compactness converges strongly to
$w^{1+\alpha}$. Moreover, $\nabla w_\eps^{(\alpha+\delta)/2} \bbin \rmL^2(Q_T)$
and \eqref{eq:w.est.e=k} imply
\[
\nabla w_\eps^{1{+}\alpha} = w_\eps^{1+(\alpha-\delta)/2} \nabla
w_\eps^{(\alpha+\delta)/2} \bbin \rmL^p(Q_T) \text{ for } 1 \leq p
 < p_\Pi:=\frac{2q_\rmw}{q_\rmw{+}1{+}\alpha} >1. 
\]
Here we need to choose $\delta \approx 1$. With this $\nabla w_\eps^{1+\alpha}
\weak \nabla w^{1+\alpha} $ in $\rmL^r(Q_T)$ follows, and the case (ii) is
settled. 

In exactly the same way with $\eta(w)=w^\alpha$ we obtain 
\begin{equation}
  \label{eq:nabla.w.alpha}
  \nabla \eta(w_\eps)\weak \nabla \eta(w) \ \text{ in } \rmL^p(Q_T)
   \text{ for all } p\in {[1,p_\eta[}    \text{ with } 
   p_\eta:= \frac{2q_\rmw}{q_\rmw {+}    \alpha {-} 1}> p_\Pi.   
\end{equation}
By \eqref{eq:stro.cvg.weps} we also have $\eta(w_\eps)\to \eta(w)$ in $
\rmL^p(Q_T $) for all $p\in {[1,q_\rmw/\alpha[}$. Thus, for the four weak convergences of
the terms in (i)', (iii)', (i)'', and (iii)'' in $\rmL^{1+\zeta}(Q_T)$ for some
$\zeta>0$, we need the four relations $1/q_\rmv + 1/p_\eta <1$, 
$2/q_\rmv + 1/p_\eta <1$, $1/q_\rmv+ \alpha/q_\rmw <1$, and $2/q_\rmv+
\alpha/q_\rmw<1$, respectively. Because of $q_\rmw/\alpha < p_\eta$ all four
inequalities follow if the second holds. Inserting the definitions of $q_\rmv$
and $p_\eta$ into $2/q_\rmv+ 1/p_\eta<1$ show that this condition is equivalent
to our assumption \eqref{eq:p0.cond.e=k}. 

Thus, the proof of Theorem \ref{th:Exi.eta=kappa} is finished. 
\end{proof}

\subsection{Fractional Sobolev spaces}
\label{su:FracSobol}

Here we provide a result that was used in the compactness arguments.
 
For $\gamma >0$ we set $\{u\}^\gamma: u\mapsto \mafo{sign}(u) |u|^\gamma$. Then,
we have the relations  
\begin{equation}
  \label{eq:alpha.Power}
\forall\, \gamma\geq 1\ \forall\, a,b\in \R: \quad 
\big|\{a\}^{1/\gamma} - \{b\}^{1/\gamma}\big|^\gamma \leq 
   2^{\gamma-1} \,\big| a - b \big|.
\end{equation}
With this we easily obtain the following statement.

\begin{proposition}[Fractional powers]\label{pr:FractPow}
  For all $\gamma>1$, $ s\in {]0,1[} $, $ p\in {]1,\infty[} $, and
  $ u\in \rmW^{s,p}(\Omega) $ we have
  $\{u\}^{1/\gamma} \in \rmW^{s/\gamma, p\gamma}(\Omega) $ with the estimate
\begin{equation}
  \label{eq:FractSobol}
  \big\| \{u\}^{1/\gamma} \big\|_{\rmW^{s/\gamma,p\gamma}(\Omega)} \ \leq \ 
  2^{(1-1/\gamma)/p} \, \big\| u \big\|_{\rmW^{s,p}(\Omega)}^{1/\gamma}.
\end{equation}
\end{proposition}
\begin{proof}
Clearly we have $\|\{u\}^{1/\gamma}\|_{\rmL^{p\gamma}}^{p\gamma} = 
\|u\|_{\rmL^p}^p $. For the Sobolev-Slobodeckij semi-norm 
$[\!\,[\,\cdot ]\!\,]_{s,p}$ we apply \eqref{eq:alpha.Power} and obtain
\begin{align*}
\big[\!\big[ \{u\}^{1/\gamma} \big]\!\big]_{s/\gamma,p\gamma}^{p\gamma}
&:=  \iint_{\Omega\ti\Omega}  \frac{\big|\{u(x)\}^{1/\gamma}{-}
  \{u(y)\}^{1/\gamma}\big|^{p\gamma}}{|x{-}y|^{d+(s/\gamma)p\gamma}}  
\dd x \dd y \ \leq \ 2^{\gamma-1} \big[\!\big[ u \big]\!\big]_{s,p}^p.
\end{align*}
This proves \eqref{eq:FractSobol}.
\end{proof}

\section{Conjectures on the self-similar behavior} 
\label{se:Scaling}

In this section we speculate about the longtime behavior of solutions on the
full space $\Omega=\R^d$. First we recall the similarity solution
$w_\mafo{PME}$ for the
classical PME in Section \ref{su:SimilPME}. In Section \ref{su:CONJ.al=be} we discuss
the case $\alpha=\beta$, where we have full scaling invariance. We
expect $\int_{\R^d} \frac12v(t,x)^2\dd x \to 0$ such that $\int_{\R^d} w(t,x)\dd x \to
\calE(\INIv,\INIw)$ for $t\to \infty$. Moreover, we conjecture that 
$w(t,\cdot)$ behaves like the
$w_\mafo{PME}$ for the corresponding mass $\calE(\INIv,\INIw)$, while $v(t)$
behaves like $t^\gamma( w_\mafo{PME})^{\kappa_0/\eta_0}$, which again shows
that $\eta_0\gg \kappa_0$ leads to singular behavior of $v$ at the boundary of
the support. Finally, Section
\ref{su:CONJ.a.neq.b} addresses the cases $\alpha >\beta>0$ and
$\beta>\alpha>0$. In the latter we still expect that $w(t)$ behaves like
$w_\mafo{PME}$ while $v(t)$ should behave like $ct^{-d\delta}
\bm1_{B_{t^\delta/b}(0)}(\cdot)$. In the former case, a similar result should
not be expected.

\subsection{Similarity solution for the classical PME}
\label{su:SimilPME}

To describe the exact similarity solutions we introduce the shape functions
\[
\bbW_\sigma(y)=\big(1{-}|y|^2\big)_+^\sigma.
\]
We have the obvious relations 
\begin{subequations}
  \label{eq:DIV.bbW}
\begin{align}
 \label{eq:DIV.bbW.a}
&\bbW_\sigma^\gamma =\bbW_{\sigma \gamma}, \quad
\bbW_\sigma\bbW_\gamma=\bbW_{\sigma+\gamma} , \quad
\Theta(\sigma,d):=\int_{\R^d}\bbW_\sigma(y) \dd y = \frac{\pi^{d/2} 
  \Gamma(1{+}\sigma)} {\Gamma(1{+}\sigma {+}d/2)},
\\[0.0em]
 \label{eq:DIV.bbW.b}
& \nabla \bbW_\sigma(y) = -2\sigma \bbW_{\sigma-1}(y) \,y, 
\\[0.2em]
 \label{eq:DIV.bbW.c}
&  \DIV\!\big(\bbW_\gamma \nabla \bbW_\sigma\big)(y) = - 2d\sigma
  \bbW_{\gamma+\sigma-1}(y) -2\sigma \,y{\cdot}\nabla \bbW_{\gamma+\sigma-1}(y) \quad
  \text{for } |y|\neq 1. 
\end{align}
\end{subequations}
With this, we easily obtain the following well-known, explicit form of the similarity
solution for the PME.

\begin{lemma}[Similarity solution for PME]\label{le:SimiPME}
The function $w_\mafo{PME}$ defined via \linebreak[4] $w_\mafo{PME}(t,x)= c (t{+}t_*)^{-d\delta}  
\bbW_\sigma(b(t{+}t_*)^{-\delta} x)$ is a solution of the PME$_\beta$ given by  
\[
\dot w = \DIV(\kappa_0 w^\beta \nabla w) , \qquad \int_{\R^d} w(t,x) \dd x = E_0
\]
if and only if the we choose the parameters such that
\begin{equation}
  \label{eq:ParaChoice}
  E_0 = c\,b^{-d}\, \Theta\big(1/\beta,d\big), \quad \sigma =1/\beta, \quad 
\delta= 1/({2{+}d\beta}), \quad 2\kappa_0 bc^\beta =\delta\beta .
\end{equation}
\end{lemma}
\begin{proof} The first relation stems from the integral constraint for $w$. 

A direct calculation using $y=bt^{-\delta}x$ gives
\begin{align*}
\dot w&= -c\delta t^{-d\delta-1} \big(d \bbW_\sigma(y)+ y\cdot\nabla \bbW_\sigma(y) \big),
\\
\DIV(\kappa_0w^\beta\nabla w)&= -2\sigma \kappa_0 b^2c^{\beta+1}
t^{-(\beta+1)d\delta - 2\delta} \big(d \bbW_{\beta\sigma+\sigma-1} (y)+ y\cdot\nabla
\bbW_{\beta\sigma+\sigma-1}(y) \big). 
\end{align*}
Thus, we first see that matching the two lines needs $\beta\sigma=1$ and
$\delta(2{+}d\beta)=1$. Finally, we compare the prefactor which gives the last
relation.   
\end{proof}

The convergence of all non-negative solutions $w$ to the self-similar profile 
is one of the major achievement in the theory of PME, see
\cite[Cha.\,16]{Vazq07PMEM}, \cite{CarTos00ALDS}, or \cite[Sec.\,3]{Otto01GDEE}.  
We believe that a corresponding result on asymptotic self-similarity for our
system, at least if $\alpha=\beta$. However, in the following we only present
conjectures together with some supporting observations, including a numerical
simulation.

\subsection{Conjectured longtime behavior for $\beta=\alpha$}
\label{su:CONJ.al=be}

To substantiate our conjecture we use the scaling properties of the coupled
system for $\beta=\alpha$. Slightly generalizing \eqref{eq:ScalTrafo} we
consider the two possible transformations, namely  
\begin{equation}
  \label{eq:ScalTra.a=b}
\begin{aligned}
& \tau=\log(t{+}1),  \ \  y=(t{+}1)^{-\delta}x \quad 
 \text{and} \quad \\
&\big(v(t,x),w(t,x)\big) = \big( (t{+}1)^{-d\theta\delta} \,\wt
 v(\tau,y)\, , \, (t{+}1)^{-d\delta}\, \wt w (\tau , y) \big),
\end{aligned}
\end{equation}
where $\delta =1/(2{+}d\beta)$ as above, and the new parameter $\theta$ either
equals \\[0.3em]
\textbullet\ $\theta=1/2$: the ``\emph{energy-conserving scaling}'' or \\[0.3em]
\textbullet\ $\theta=1$: the ``\emph{momentum-conserving scaling}''. 

The reason for these two choices is that the variable $v$ occurs with different
powers in the two conserved quantities quantities
$\calV(v,w)=\int_{\R^d} v(x) \dd x$ and $\calE(v,w)=\int_{\R^d}\big( \frac12
v^2{+}w\big) \dd x $. Thus, a given scaling can conserve at most one of the
two functionals. 

The transformed coupled PDE system for $(\wt
v,\wt w)$ reads 
\begin{subequations}
  \label{eq:CONJ.traSys}
\begin{align}
 \label{eq:CONJ.tS.wtv}
\pl_\tau \wt v &  =\DIV\!\big(\delta \,\wt v \: y + \eta_0\wt w^\beta \nabla
\wt v \big) - (1{-}\theta) \delta d \,v,
\\ 
 \label{eq:CONJ.tS.wtw}
\pl_\tau \wt w &=  \DIV\!\big(\delta \wt w \, y + \kappa_0\wt w^\beta \nabla \wt
  w\big)  + \ee^{(1{-}2\theta)\delta d\, \tau} \eta_0\wt w^\beta|\nabla \wt v|^2.
\end{align}
\end{subequations}
The conserved total momentum $\calV$ and total energy $\calE$ transform as
follows: 
\begin{subequations}
\begin{align}
  \label{eq:CONJ.calV}
 &\wt\calV(\tau,\wt v,\wt w):= \int_{y\in \R^d} \ee^{(1-\theta) \delta d\tau} 
 \wt v(y) \dd y = \int_{x\in
   \R^d} v(x) \dd x =  \calV(v,w)=\calV(\INIv,\INIw)=V_0,
\\
& \label{eq:CONJ.calE}
  \wt\calE\big(\tau,\wt v(\tau),\wt w(\tau)\big):= \int_{\R^d}\!\! \big(
  \tfrac12 \,\ee^{(1-2\theta)d\delta\,\tau}\,\wt v(\tau,y)^2 {+} \wt w(\tau,y) \big) \dd y
\\ \nonumber
& \hspace*{11.5em} 
 = \calE\big(v(\ee^\tau{-}1), w(\ee^\tau{-}1)\big) = \calE(\INIv,\INIw)= E_0.
\end{align}
\end{subequations}
Thus, the linear momentum $\calV$ is remains a conserved quantity (i.e.\
$\calV(\tau,\cdot) = \calV$) for $\theta=1$
only. In contrast, the energy $\calE$ remains a conserved quantity (i.e.\
$\wt\calE(\tau, \cdot) = \calE$) for $\theta=1/2$ only. Moreover, system
\eqref{eq:CONJ.traSys} is a autonomous system if and only if $\theta=1/2$.

\subsubsection{Energy-conserving scaling with $\theta=1/2$}  
\label{suu:theta1/2}

For $\theta=1/2$ the transformed system \eqref{eq:CONJ.traSys} takes the
explicit form 
\begin{subequations}
  \label{eq:CONJ.traSys1/2}
\begin{align}
 \label{eq:CONJ.tS1/2.wtv}
\pl_\tau \wt v &  =\DIV\!\big(\delta \,\wt v \: y + \eta_0\wt w^\beta \nabla
\wt v \big) - \frac12 \delta d \,v,
\\ 
 \label{eq:CONJ.tS1/2.wtw}
\pl_\tau \wt w &=  \DIV\!\big(\delta \wt w \, y + \kappa_0\wt w^\beta \nabla \wt
  w\big)  +  \eta_0\wt w^\beta|\nabla \wt v|^2,
\end{align}
\end{subequations}
which is an autonomous evolutionary system with the conserved quantity
$\calE(\wt v, \wt w)$. To understand the longtime behavior, we can test
\eqref{eq:CONJ.tS1/2.wtv} by $|\wt v|^{p-2} \wt v$ for $p\in {]1,2]}$ and find
\begin{equation}
  \label{eq:wtv.Lp}
  \frac\rmd{\rmd \tau} \int_{\R^d} \frac1p|\wt v|^p \dd y 
= -(p{-}1)\int_{\R^d} \eta_0\wt w^\beta \,|\wt v|^{p-2} \,|\nabla \wt v|^2 \dd y 
+ \delta d\Big(\frac{p{-}1}p -\frac12\Big) \int_{\R^d} |\wt v|^p \dd y .
\end{equation}
This implies exponential decay of $\|\wt v(\tau,\cdot)\|_{\rmL^p}$ for $p \in
{]1,2[}$.  Moreover, for $p=2$ all steady states (where the left-hand side in
the above relation is $0$) satisfy $\nabla \wt v= 0$. As $\calE(\wt v,\wt
w)=E_0 < \infty$ we conclude $\wt v =0$ and obtain the following result.

\begin{corollary}[Unique non-negative steady state]
\label{co:SteadyState}
Given arbitrary $E_0\geq 0$ there is exactly one non-negative steady state for 
\eqref{eq:CONJ.traSys}$_{\theta=1/2}$, namely $(\wt v^{E_0}_\st,\wt w^{E_0}_\st)= \big(
0,\, c\bbW_{1/\beta}(b\cdot)\big)$ where $b$ and $c$ depend on $E_0$ as given
in \eqref{eq:ParaChoice}.
\end{corollary} 

We conjecture that all solutions of \eqref{eq:CONJ.traSys}$_{\theta=1/2}$ with
initial condition satisfying $\calE(\INIv,\INIw)=E_0$ converge to 
$(\wt v^{E_0}_\st,\wt w^{E_0}_\st)$. The reasoning is simple because we know the
convergence of $\wt v$ to $0$, so we expect that for large times the dynamics
of the PME for $\wt w$ (now in transformed variables) and $\wt v\equiv 0$ will
determine the longtime behavior. However, to show this one needs to show that
\emph{all} kinetic energy is converted to heat, i.e.\ $\| \wt
v(\tau,\cdot)\|_{\rmL^2} \to 0$ for $\tau \to \infty$, which would need a more
careful analysis than the simple estimate \eqref{eq:wtv.Lp}.

\subsubsection{Momentum-conserving scaling with $\theta=1$}  
\label{suu:theta1}

For $\theta=1$ 
the transformed system \eqref{eq:CONJ.traSys} takes the
explicit form 
\begin{subequations}
  \label{eq:CONJ.traSys1}
\begin{align}
 \label{eq:CONJ.tS1.whv}
\pl_\tau \wh v &  =\DIV\!\big(\delta \,\wh v \: y + \eta_0\wt w^\beta \nabla
\wh v \big) ,
\\ 
 \label{eq:CONJ.tS1/2.wtw}
\pl_\tau \wt w &=  \DIV\!\big(\delta \wt w \, y + \kappa_0\wt w^\beta \nabla \wt
  w\big)  +  \ee^{-\delta d\,\tau} \eta_0\wt w^\beta|\nabla \wh v|^2,
\end{align}
\end{subequations}
We now use $\wh v$ which is related to $\wt v$ in Section \ref{suu:theta1/2} by
$\wh v(\tau)=\ee^{\delta d \tau /2} \wt v(\tau)$. It is important to keep in
mind that the component $\wt w$ remains exactly the same in both systems,
\eqref{eq:CONJ.traSys1/2}  and \eqref{eq:CONJ.traSys1}.

The new system is no longer autonomous but the time dependence occurs via an
exponentially decaying term, which hopefully does not influence the longtime
behavior even if $\wh v$ is not decaying but stays suitably bounded. 
From the previous Section  \ref{suu:theta1/2} we still expect that $\wt
w(\tau)$ converges to $w_\st^{E_0}$. 

However, the new feature of \eqref{eq:CONJ.traSys1} is that $\wh v$ again
satisfies a PDE in divergence form, namely \eqref{eq:CONJ.tS1.whv}, which
implies that $\calV(\wh w(\tau)= V_0$ for all $\tau\geq 0$. Moreover, as we
already expect $\wt w(\tau) \to  w_\st^{E_0}$ we may expect that $\wh v(\tau)$
converges to a steady state of the linear diffusion equation obtained from
\eqref{eq:CONJ.tS1.whv} by replacing $\wt w$ by its limit $ w_\st^{E_0}$.  
The surprising fact is that this equation has a unique steady state which is
given by an explicit formula.

\begin{proposition}[Similarity solutions $(\wh v_\st,\wt w_\st)$]
\label{pr:CONJ.StSt} 
For each pair $(V_0,E_0)\in \R\ti{]0,\infty[}$ there is a unique solution
$(\wh v_\st,\wt w_\st)$ of the system (where $\delta=1/(2{+}d\beta )\,$)
\[
0 =\DIV\!\big(\delta \wh v \, y + \eta_0\wt w^\beta \nabla \wh v \big) , \quad 
0=  \DIV\!\big(\delta \wt w \, y + \kappa_0\wt w^\beta \nabla \wt w\big) , 
\quad \int_{\R^d} \!\! \wh v \dd y =V_0, \quad \int_{\R^d}\!\!  \wt w \dd y =E_0 .
\]
This solution is given explicitly in the form 
\begin{align*}
&\big(\wh v_\st(y),\wt w_\st(y)\big)  = \Big(a \bbW_\sigma(by)\,,\:  
 c \bbW_{1/\beta}(by)\Big) \quad \text{with } \sigma =
 \frac{\kappa_0}{\eta_0 \beta}, \quad  b^2=\frac{\beta
  \delta}{2c^\beta \kappa_0}  ,
\\
& \hspace{12em} 
V_0 =\frac a{b^d} \Theta(\kappa_0/(\beta\eta_0),d), \quad E_0= \frac c{b^d}
\Theta(1/\beta,d). 
\end{align*}
\end{proposition}
\begin{proof} The result follows easily with \eqref{eq:DIV.bbW.b}. Indeed, we
  even have $ \delta \wh v \, y + \eta_0\wt w^\beta \nabla \wh v \equiv 0$ and 
$ \delta \wt w \, y + \kappa_0\wt w^\beta \nabla \wt w  \equiv 0 $ on $\R^d$. 
\end{proof}

Of course, justifying the convergence to the steady states
$(\wh v_\st,\wt w_\st)$ for a suitable solution class $(\wh v,\wt w)$ of the
full system \eqref{eq:CONJ.traSys} is a nontrivial task. This can be seen by
inserting the steady state into the exponentially decaying term which gives the
perturbation 
\[
\ee^{-d\delta \tau} \eta_0 \wt w_\st^\alpha |\nabla \wh v_\st|^2= \wt c
\,\ee^{-d\delta \tau} \bbW_{2\kappa_0/(\beta\eta_0)-1}(by) |y|^2. 
\]
This term lies in $\rmL^\infty(0,T;\rmL^1(\R^d))$ but is singular at $b|y|=1$ for $\beta
>2\kappa_0/\eta_0$. 

\subsubsection{A numerical simulation showing convergence}
\label{suu:Simulation}

\begin{figure}
\includegraphics[width=0.3\textwidth,trim=10 25 20 80,clip=true]{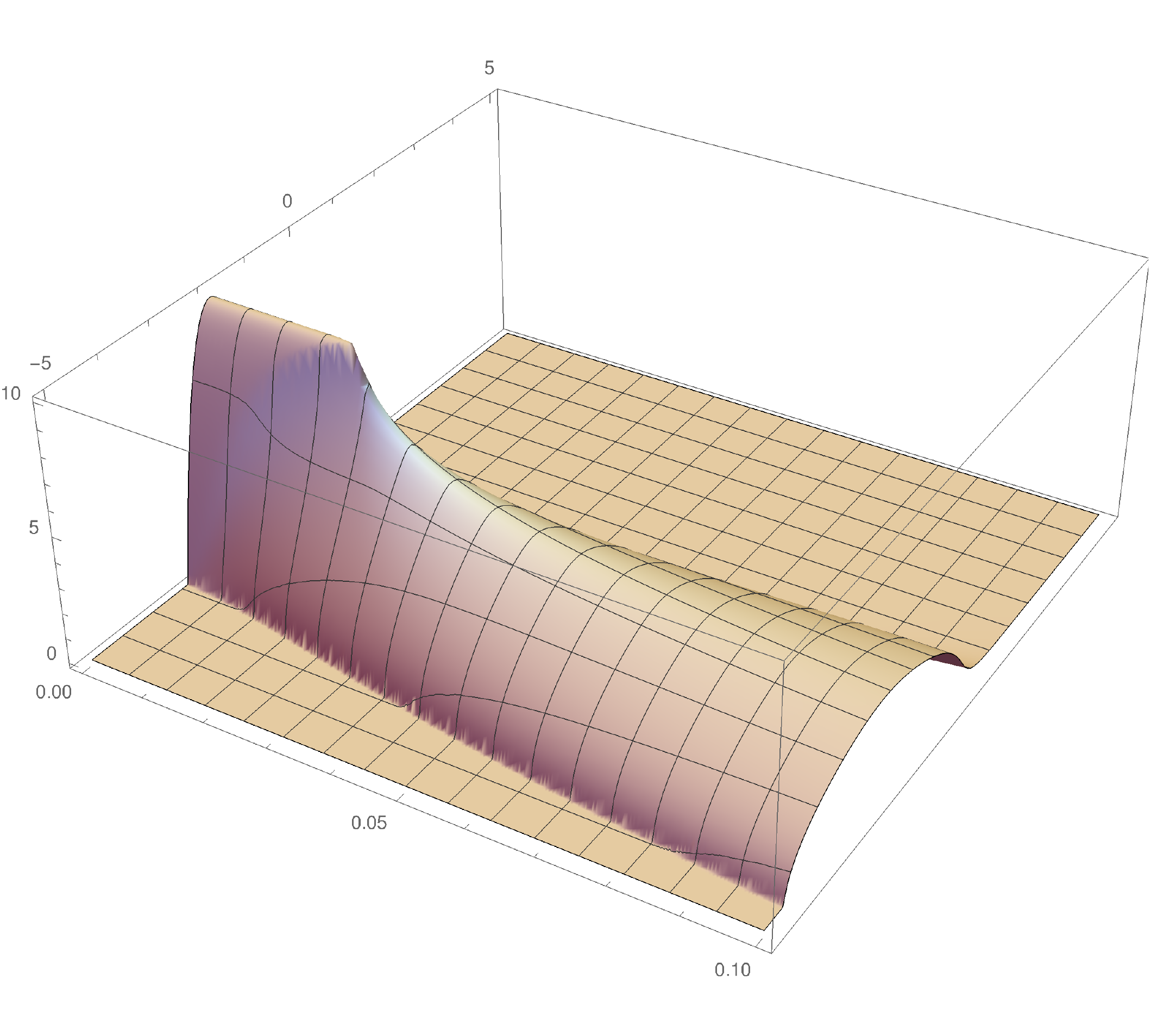} 
\quad
\includegraphics[width=0.3\textwidth,trim=15 20 30 40,clip=true]{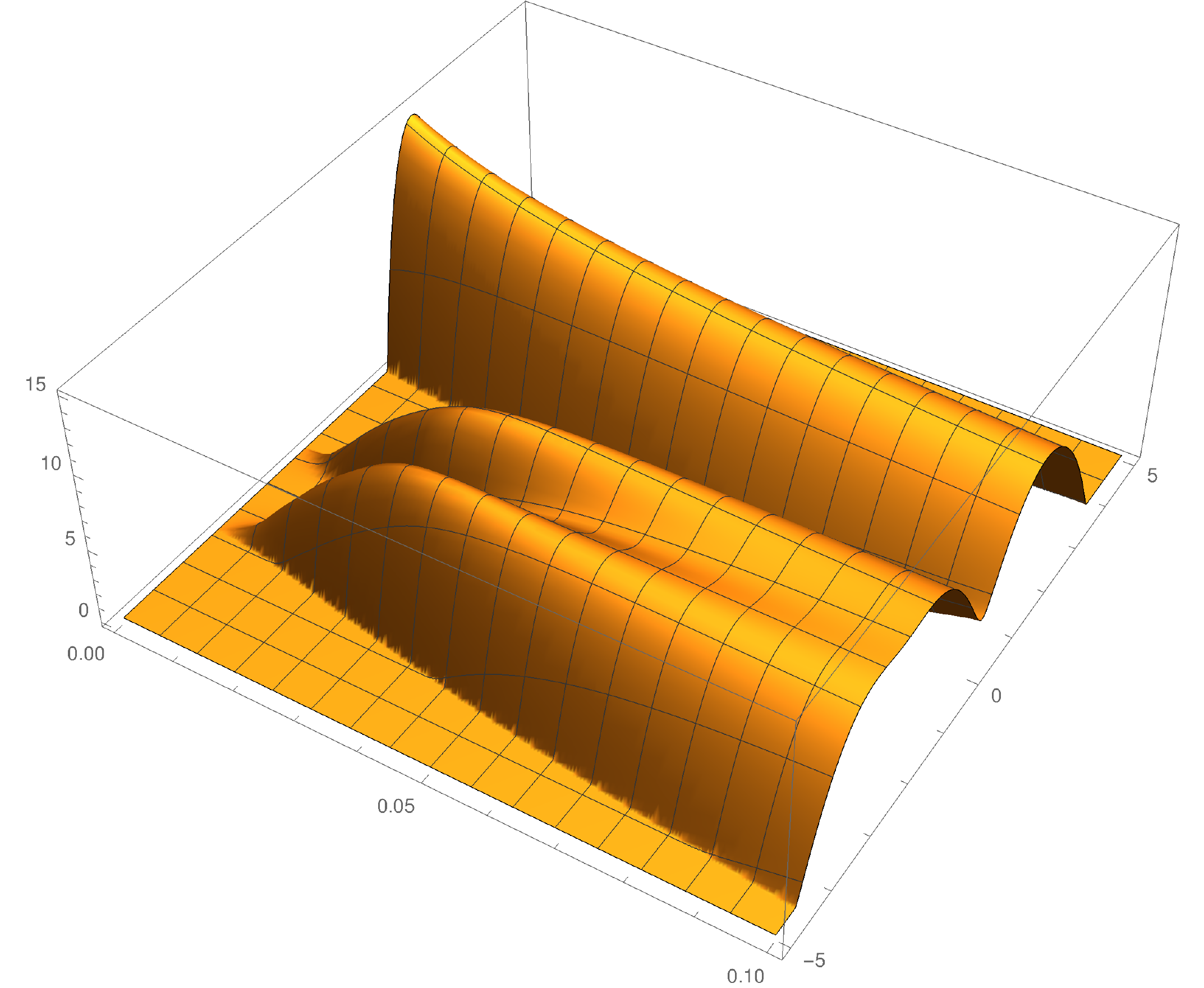} 
\hfill
\begin{minipage}[b]{0.3\textwidth}
  \caption{\small The solution $v(t,x)$ (left) and $w(t,x)$ (right) for
    $t\in [0,0.1]$ and $x\in [-5,5]$.}
\label{fig:initial}
\end{minipage}
\end{figure}

A simple numerical experiment covers the one-dimensional case $\Omega = \R^1$ with
$\alpha=\beta=1$, $\eta_0=2$, and $\kappa_0=1/2$. We start with initial
conditions $\INIv(x)= \max\{0, 10-10(x{+}2)^2\}$ and $\INIw(x)= \max\{0,
15-15(x{-}2)^2\}$. Fin Figure \ref{fig:initial} we show that solution for the
short initial time interval $t\in [0, 0.1]$ which shows that kinetic energy is
dissipated fast and turned into heat.  
% obtained via Mathematica notebood UntransfoCDPS_008.nb   

\begin{figure}
\includegraphics[width=0.3\textwidth,trim=0 20 20 30,clip=true]{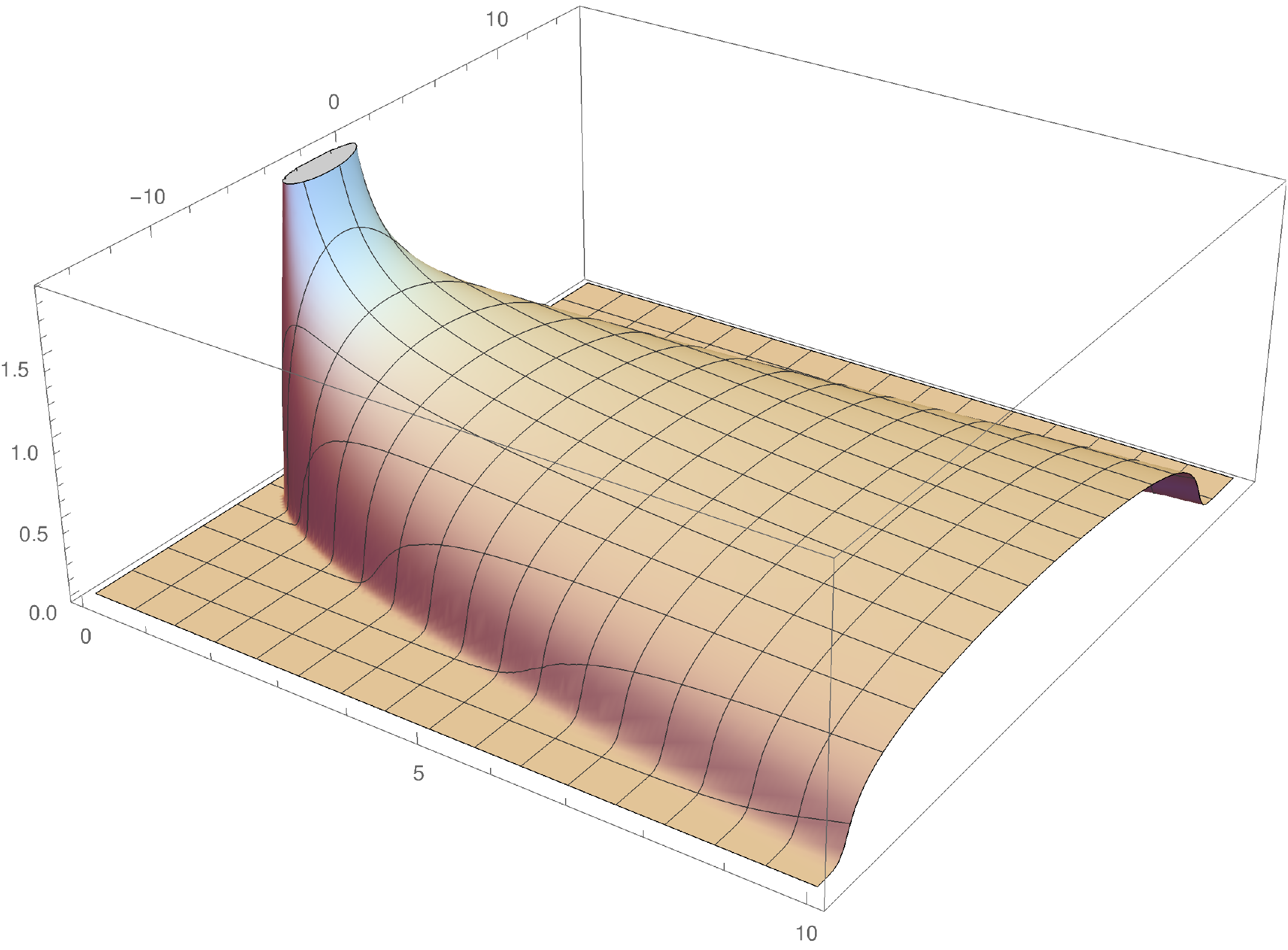} 
\quad
\includegraphics[width=0.3\textwidth,trim=10 20 20 40,clip=true]{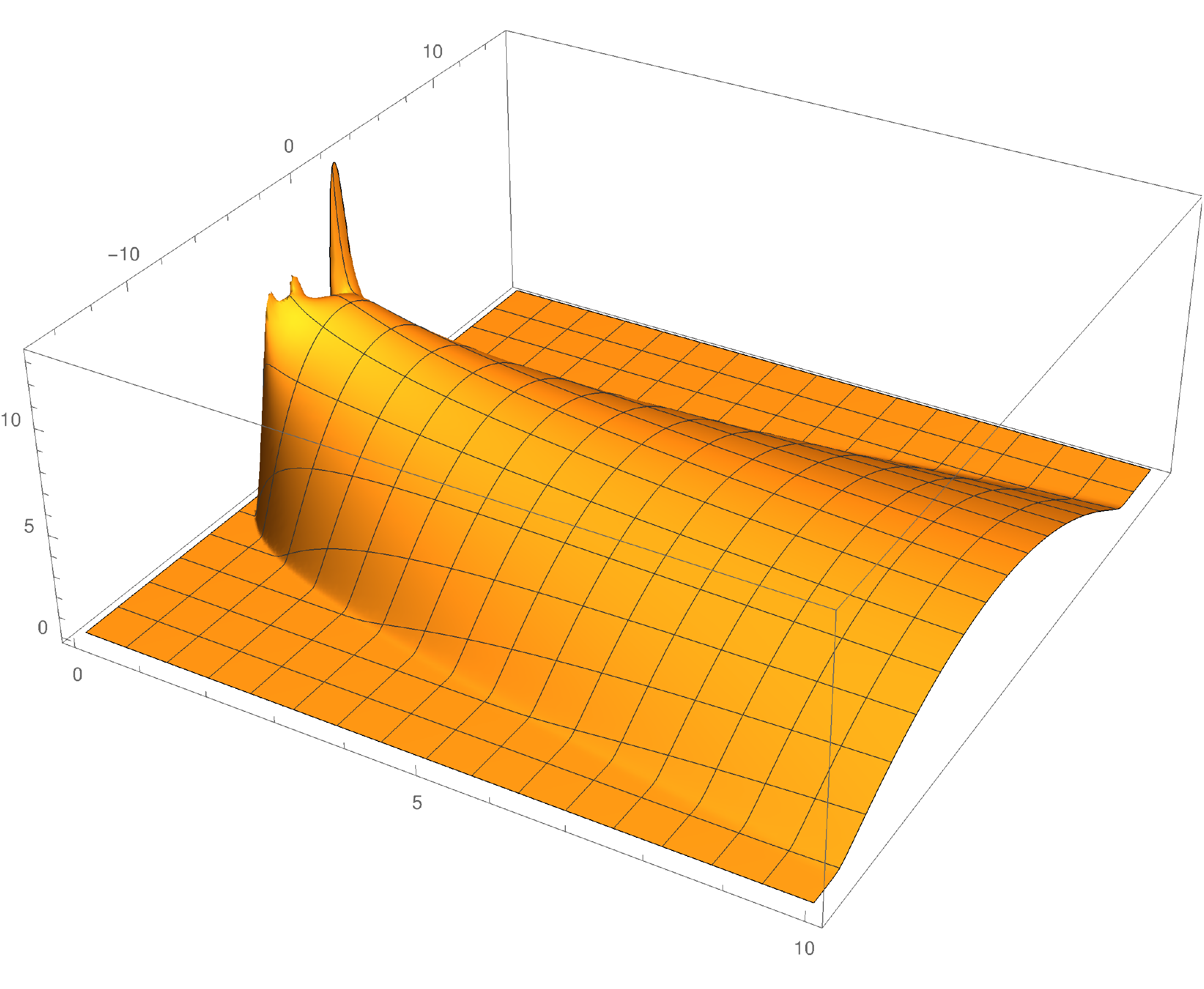} 
\hfill
\begin{minipage}[b]{0.3\textwidth}
  \caption{\small The solution $v(t,x)$ (left) and $w(t,x)$ (right) for
    $t\in [0,0.1]$ and $x\in [-5,5]$.}
\label{fig:global}
\end{minipage}
\end{figure}

In Figure \ref{fig:global} we show that (unscaled) solutions on the longer time interval
$t\in [0,10]$, where the self-similar behavior becomes visible. 

\begin{figure}
\includegraphics[width=0.33\textwidth,trim=12 10 30 50,clip=true]{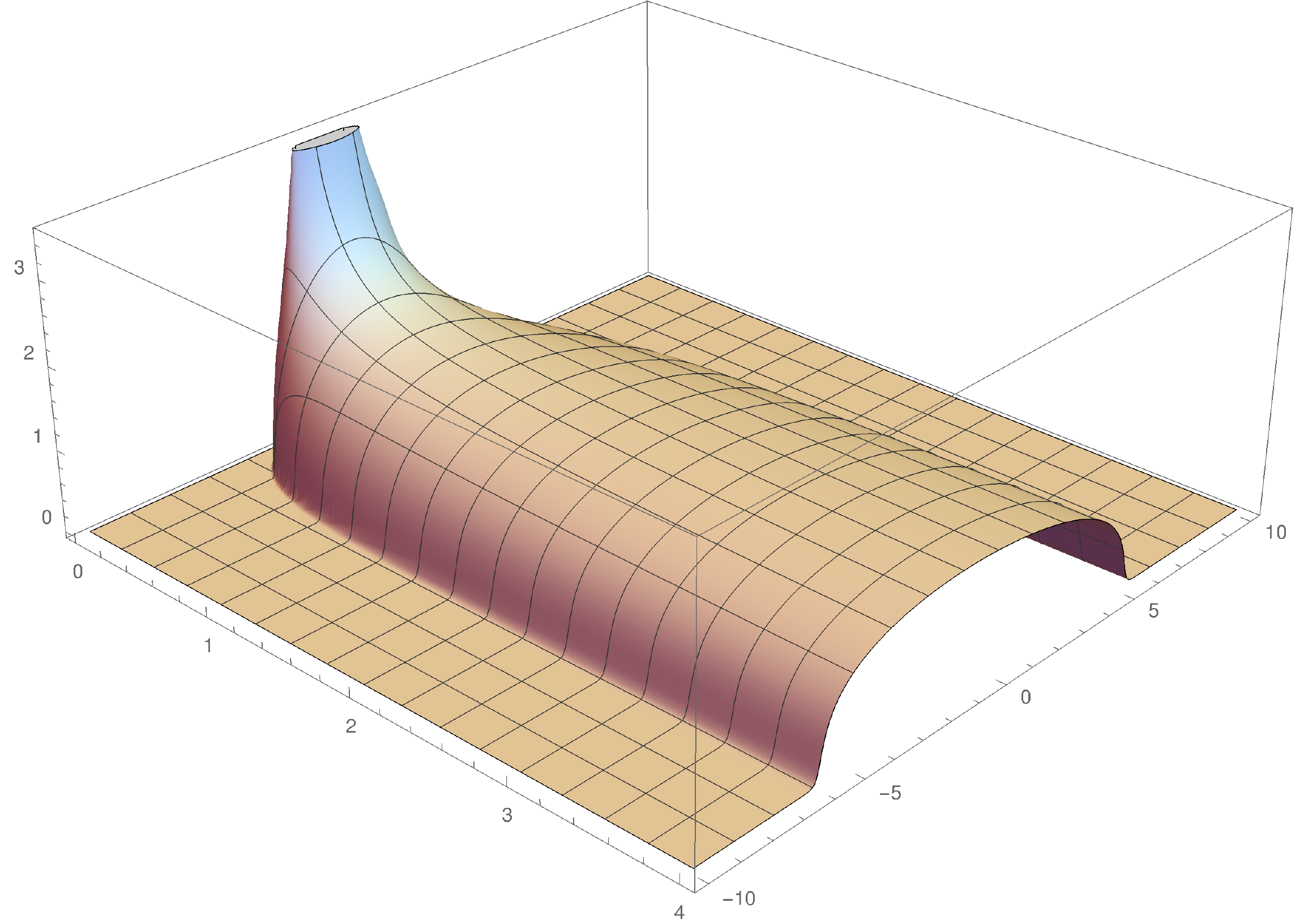}
\includegraphics[width=0.33\textwidth,trim=15 20 63 70,clip=true]{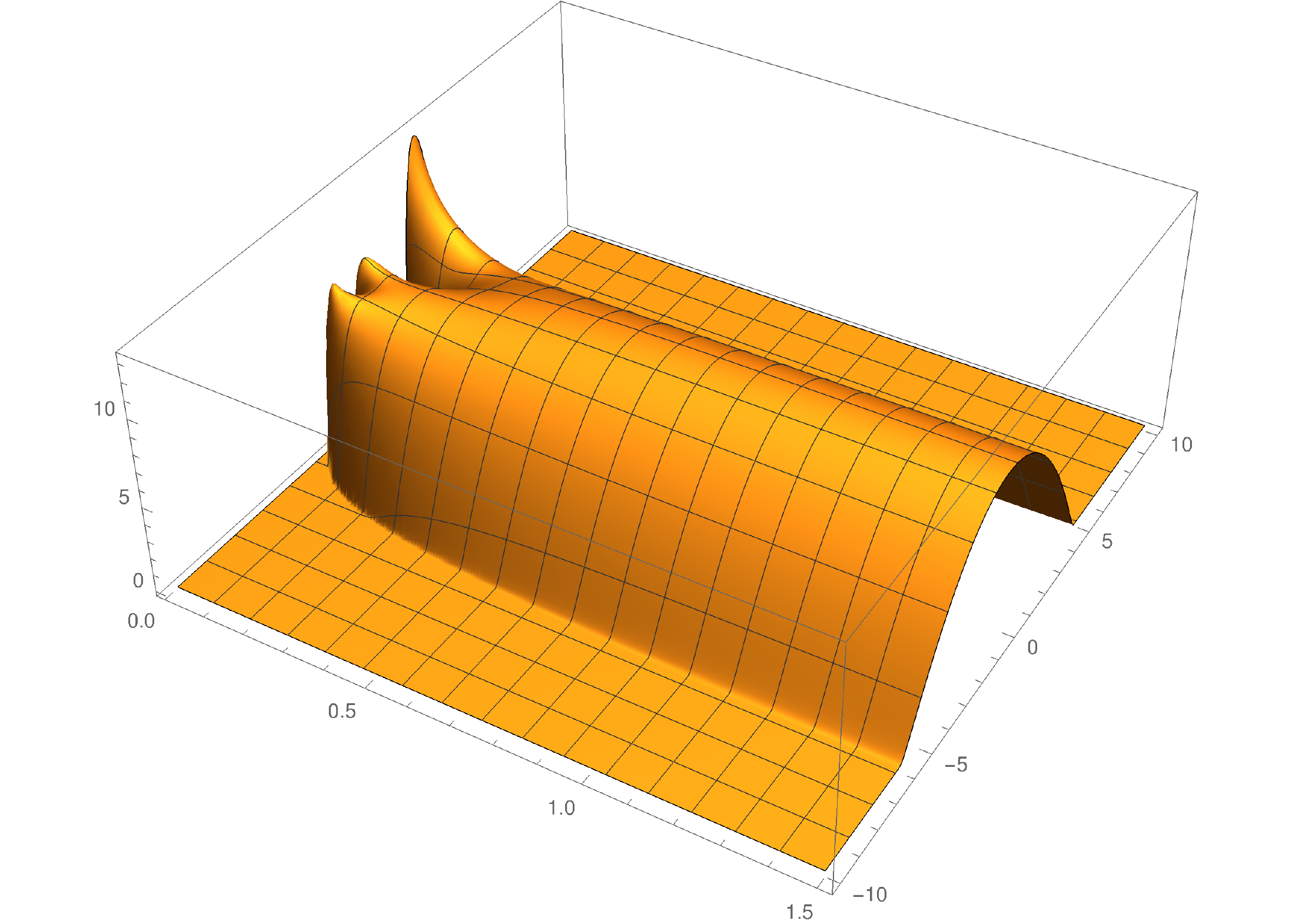}
\hfill
\begin{minipage}[b]{0.32\textwidth}
\caption{\small Scaled solutions $(1{+}t)^\delta v(t,(1{+}t)^\delta y) $ (left) and 
 $(1{+}t)^\delta w(t,(1{+}t)^\delta y) $ (right) for $t\in [0,1.5]$ and $y\in [-10,10]$.} 
\label{fig:scaled}
\end{minipage}
\end{figure}

In Figure \ref{fig:scaled} we show that the rescaled solutions
$(t,y)\mapsto (1{+}t)^\delta v(t,(1{+}t)^\delta y)$ and $(t,y)\mapsto 
(1{+}t)^\delta w(t,(1{+}t)^\delta y)$ for $t\in [0,1.5]$ where convergence into a self-similar profile is
    already evident. It is clearly seen that $w$ develops the simple Barenblatt
    profile $c_w\max\{0, b{-}y^2\}$ ,  whereas $v$ develops the more singular profile
    $c_v\max\{0, b{-}y^2\}^{1/4}$ because of $\kappa_0/\eta_0=1/4$.  

\subsection{Conjectured longtime behavior for $\beta \neq\alpha$}
\label{su:CONJ.a.neq.b}

With the transformation \eqref{eq:ScalTra.a=b} with
$\delta=1/(2{+}d\beta)$ we now find the transformed system 
\begin{subequations}
  \label{eq:CONJ.tS.a.neq.b}
\begin{align}
 \label{eq:CONJ.a.n.b.V}
\pl_\tau \wt v\: &  =\DIV\!\big(\delta \,\wt v \: y + \ee^{\gamma \tau} 
 \eta_0\wt w^\beta \nabla \wt v \big) \qquad \text{with }\gamma =
 (\beta{-}\alpha) d \delta,
\\ 
 \label{eq:CONJ.a.n.b.W}
\pl_\tau \wt w &=  \DIV\!\big(\delta \wt w \, y + \kappa_0\wt w^\beta \nabla \wt
  w\big)  + \ee^{(\gamma-d\delta)\tau} \eta_0\wt w^\beta|\nabla \wt v|^2.
\end{align}
\end{subequations}
We consider the cases $\alpha >\beta>0$ and 
$\beta>\alpha>0$ separately.\medskip

\STEP{\underline{Case $\gamma=(\beta{-}\alpha) d \delta >0$:}} The chosen scaling
leads to the exponentially growing prefactor for the diffusion in 
\eqref{eq:CONJ.a.n.b.V}. Hence it is dominating the concentration term such
that we expect for large $\tau$ that $\nabla v(\tau)$ is small. This can also
be seen by testing \eqref{eq:CONJ.a.n.b.V} by $\ee^{-d\delta \tau}v(\tau)$
which leads to 
\[
\iint_{Q_T}\!\!\ee^{(\gamma-d\delta)\tau} \eta_0\wt w^\beta|\nabla \wt v|^2 \dd
y \dd \tau  \leq \frac12\|\wt v(0)\|_{\rmL^2(R^d)}^2.  
\]  
Thus, it is reasonable to expect that $\wt v(t)$ converges to $c \,
\bm1_{B_{1/b}(0)}(\cdot)$ while $\wt w$ converges to $\wt w_\st$ as for the
case $\alpha=\beta$. 

A proof of this conjecture will be even more difficult, because $\nabla
v(\tau)$ must develop a singularity near the boundary of $ B_{1/b}(0)\subset
\R^d$ to allow for the convergence $\wt v(\tau) \to c\,
\bm1_{B_{1/b}(0)}(\cdot)$ in $\rmL^2(\R^d)$. 
\medskip

\STEP{\underline{Case $\gamma=(\beta{-}\alpha) d \delta <0$:}} In this case we
expect a different behavior because of the decaying diffusion in the equation
for $\wt v$. For large $\tau$ the dominating term in \eqref{eq:CONJ.a.n.b.V} is
the transport term $\DIV\big( \delta \wt v y\big)$ which is simply due to the
similarity rescaling. Hence, for $\tau >\tau_0\gg 1$ one expects
$\wt v(\tau,y) \approx \ee^{d\delta(t-\tau_0)} \wt
v\big(\tau_0,\ee^{\delta(\tau-\tau_0)} y\big) $, which in the original physical
coordinates means the same as $v(t,x)\approx v(t_0,x)$ for $t > t_0 \gg
1$. Thus, it is likely that
$\| v(t)\|_{\rmL^2}^2 = \ee^{-d\delta\tau}\|\wt v(\tau)\|_{\rmL^2}^2 $ does not
converge to $0$ for $ \tau=\log(t{+}1)\to \infty$. Consequently, the
nondecreasing function $\wt E(t):=\int_{\R^d} \wt w(\tau,y)\dd y$ is converging
to some $\wt E_\infty \in {[E(0),E_0]}$ with $E_0=\calE(\INIv,\INIw)$. There is
still a chance that $\wt w$ converges to some $\wt w_\st$ but now with mass
$\wt E_\infty$ instead of $E_0$.

\section{Parabolic systems in turbulence modeling}
\label{su:DerivFluids}

The author's motivation for studying the given class of coupled degenerate
parabolic equation comes from turbulence modeling. However, for the sake of a
concise presentation the model was simplified considerably, but still keeping
the main feature of a velocity-type variable with a degenerate viscosity
depending on the energy-like variable, which is now the mean turbulent kinetic energy
$k\geq 0$ playing the role of $w$.  We refer to
\cite{Lewa97MACT,GLLMT03TSUE,LedLew07RANS,DruNau09EWSS,Naum13DPPT,ChaLew14MNFT} for a 
mathematical exposition of the ideas behind of such a turbulent modeling and to
\cite[Ch.\,4]{Wilc93TMCF} for a fluid mechanical approach.
  
We first introduce Prandtl's model for turbulence and then 
Kolmogorov's two-equation model. In both cases, we introduce the full model
with the macroscopically mean velocity $\bfu$ satisfying the incompressibility
$\Div \bfu=0$, the pressure $p$, and the turbulent kinetic energy $k\geq
0$. Then, we show that restricting to simple shear flows of the form  $\bfu(t,x)=
\big(0,...,0,v(t,x_1,...,x_{d})\big) \in \R^{d+1}$  we obtain a system that has the form of
our coupled system \eqref{eq:SM01} with some simple extra terms.  

\subsection{Kolmogorov's one-equation model = Prandtl's model}
\label{suu:Kolmo1eqn}

Prandtl's one-equation model for turbulence was developed from 1925 to 1945,
see \cite{Pran25BUAT,Pran46UNFA} and the historical remarks
about the development of the model in \cite[p.\,20]{Naum13DPPT}. 
\begin{align*}
&\dot\bfu +\bfu\cdot\ol\nabla \bfu -\ol\nabla p= \ol\Div\big(\ell \sqrt k \,
\ol\rmD(\bfu)\big) +f, \quad \ol\Div\bfu=0,
\\
&\dot k+ \bfu\cdot \ol\nabla k = \wh\kappa  \,
\ol\Div\big(\ell \sqrt k  \,\ol\nabla k\big) +\ell \sqrt k \,|\ol\rmD(\bfu)|^2 -
a\,\frac{k\sqrt k}{\ell},
\end{align*}
where $a\geq 0$ is a dimensional parameter and $\ell>0$ is a constant
characteristic length that may depend on position
(e.g.\ via the distance to the wall). Moreover,  $\ol\nabla$, $\ol\rmD$, and
$\ol\Div$ stand for the differential operators in $\R^{d+1}$ with $d\in\{1,2\}$.

Looking for shear flows  as defined above and assuming that the external force
vanishes, i.e.\  $f\equiv 0$, and that $k$ is independent of $x_{d+1}$ 
we obtain the coupled system 
\begin{equation}
  \label{eq:ShearPrandtl}
  \dot v= \Div\!\big(\eta_0k^{1/2} \nabla v\big), \qquad 
\dot k= \Div\!\big( \kappa_0 k^{1/2}\nabla k\big) + \eta_0 k^{1/2} |\nabla v|^2  -
\frac{a}{\ell}\,k^{3/2} 
\end{equation}
with $\eta_0 =\ell/2$ and $\kappa_0=\wh\kappa \ell$.
Clearly, when neglecting the term involving $a$, we arrive at our coupled system
in the case $\alpha=\beta=1/2$. 

The case $a>0 $ does not pose any additional problem when constructing
solutions. In particular, the term reduces the total energy in the form
$\frac\rmd{\rmd t} \calE(v,k)=- \frac a\ell \int_\Omega k^{3/2}\dd x \leq 0$. 
Hence, the a priori estimates of this equation remain the same, such that existence
of weak solutions with $k(t,x)\geq \underline k >0$ was shown already in
\cite{Naum13DPPT}. Our existence result for weak solutions in Theorem
\ref{th:ExistWeakSol} works for $p_0=2$ because of $\alpha=\beta=1/2$. 
  
The general scaling \eqref{eq:ScalTrafo} leading to \eqref{eq:TraS2ParaSyst} 
can still be applied, but to have the same scaling behavior for the terms $\dot
k$ and $k^{3/2}$ we have only one choice, namely $\gamma=2$ and with
$\alpha=1/2$ and $2\delta+\alpha\gamma=1$ we find $\delta=0$. Hence, for
$(\wh v,\wh w)=\big( \frac1t v,\frac1{t^2} w\big)$ we
find the rescaled equation
\begin{equation}
  \label{eq:PrandtlTraS2}
 \begin{aligned}
 \pl_\tau \wh v\: -\ \wh v\ &= \DIV\!\big( \eta_0 \wh k{}^{1/2}\, \nabla
\wh v\big),  
\\
\pl_\tau \wh k- 2 \wh k \,&= 
 \DIV\!\big( \kappa_0\wh k{}^{1/2}\nabla \wh k\big) + \eta_0\wh k{}^{1/2}
 \big|\nabla \wh v\big|^2 - \frac a\ell\, \wh k{}^{3/2}. 
\end{aligned} 
\end{equation}
We may still consider the conserved quantities and find the relations 
\[
\frac\rmd {\rmd \tau} \int_\Omega\! \wh v \dd y = \int_\Omega \wh v \dd y 
 \ \text{ and } \
\frac\rmd {\rmd \tau} \int_\Omega\! \big(\frac12\wh v^2 {+} \wh k \big) \dd y =
 2\int_\Omega\! \big(\frac12\wh v^2 {+} \wh k \big) \dd y - \int_\Omega \frac
a\ell\,\wh k^{3/2} \dd y .
\]
Thus, we may expect that \eqref{eq:PrandtlTraS2} has steady states 
$(v_\st,k_\st)$ with nontrivial $k_\st\geq 0$ (but necessarily with $\int_\Omega v_\st
\dd y =0$) that correspond to solutions of  \eqref{eq:ShearPrandtl} that
decay, namely $(v(t,x),k(t,x))=(\frac1t v_\st(x), \frac1{t^2} k_\st(x))$.

\subsection{Kolmogorov's $k$-$\omega$ model} 
\label{suu:Model}

Kolmogorov's model \cite{Kolm42ETMI} (see \cite{Spal91KTEM} for a translation,
sometimes also called Wilcox $k$-$\omega$ model because of \cite{Wilc93TMCF})
consists of the Navier-Stokes equation for the velocity $\bfv \in \R^d$ and the
pressure $p$ coupled to two scalar equations, namely
for the specific dissipation rate $\omega>0$ (or better dissipation per unit
turbulent kinetic energy) and the turbulent energy $k\geq 0$:
\begin{subequations}
  \label{eq:KoS}
\begin{align}
 \label{eq:KoS.NS}
  \dot\bfu + \bfu\cdot \ol\nabla \bfu + \ol\nabla p &= \mu_1\, \ol\Div\big(
                    \tfrac k\omega \ol\rmD(\bfu)\big), 
\quad \ol\Div \,\bfu = 0, 
\\
\label{eq:KoS.om}
 \dot \omega + \bfu\cdot \ol \nabla \omega &= \mu_2 \, \ol\Div \big(\tfrac k\omega
  \ol\nabla \omega\big) - \alpha_1 \omega^2,\\ 
\label{eq:KoS.k}
  \dot k + \bfu\cdot \ol\nabla k  &= \mu_3 \, \ol\Div  \big( \tfrac k \omega 
  \ol\nabla k\big) + \mu_1\tfrac k\omega\, \big|\ol\rmD(\bfu)\big|^2 - \alpha_2 k\omega . 
\end{align}
\end{subequations}
The main point is that all quantities are convected with the fluid
velocity $\bfv$ and that all quantities diffuse with ``viscosities''
that are proportional to $k/\omega$ by the dimensionless factors $\mu_j>0$. 
There are sink terms in the equations for $\omega$ and $k$, namely
$-\alpha_1 \omega^2$ and $-\alpha_2 k\omega$ with dimensionless nonnegative
constants $\alpha_j$. 
The nonlinearities in system \eqref{eq:KoS} are devised in a specific way to
allow for a rich scaling group, see 
\cite[Sec.\,2]{MieNau18?EGTW}. 

The existence of weak solutions to the above coupled system is studied in 
\cite{MieNau15GTEW,BulMal19LDAK,MieNau18?EGTW}.  However, in both of these
papers the construction of the solutions for \eqref{eq:KoS} strongly relies on
lower bounds on $k(t,x)$. In \cite{MieNau15GTEW,MieNau18?EGTW} the assumption
$k_0(x)=k(0,x)\geq k_* >0$ a.e.\ in $\Omega$ was used, whereas
\cite{BulMal19LDAK} shows that the weaker assumption $\int_\Omega \big|\log
k_0(x)  \big| \dd x < \infty$ is sufficient. 

Here we are interested in solutions that may have compact
support, i.e.\ $k(t,x)=0$ is allowed on a set of positive measure. In those
regions the regularizing diffusion terms $\mu_j \DIV(\frac k\omega \nabla \,\cdot\,)$
disappear. Thus, we also restrict to shear flows as above: 
\begin{equation}
  \label{eq:KolmoShear}
  \begin{aligned}
\dot v&= \tfrac{\mu_1}2 \DIV\! \big(\tfrac k\omega \nabla v\big), 
\qquad
\dot \omega= \mu_2 \DIV\! \big(\tfrac k\omega \nabla \omega\big) - \omega^2,
\\
\dot k&= \mu_3 \DIV \!\big(\tfrac k\omega \nabla k\big) + \tfrac{\mu_1}2 \, 
 \tfrac k\omega|\nabla v|^2 - a k\omega. 
\end{aligned} 
\end{equation}

The special symmetry of the nonlinearities allows us for more scalings than the
Prandtl equation \eqref{eq:ShearPrandtl}, because \eqref{eq:KolmoShear} only
contains dimensionless parameters $\mu_j$ and $a$ whereas
\eqref{eq:ShearPrandtl} contains the Prandtl length $\ell$. 
We can look for solutions in
the following form (cf.\ \eqref{eq:ScalTrafo}), where now $\omega$ is given by
an explicit solution that is independent of $v$ and $k$: 
\begin{align*}
& \tau= \log(t{+}1), \quad y=(t{+}1)^{-\delta} x, \quad \text{and} 
\\
&\big(v(t,x),\omega(t,x),k(t,x)) = \big( {(t{+}1)^{-\gamma/2}} \:\wh
  v(\tau,y),\ \frac1{t{+}1}, \ {(t{+}1)^{-\gamma}}\: \wh k 
  (\tau,y)\big). 
\end{align*}
If the exponents $\gamma$ and $\delta$ are chosen with $\gamma+2\delta=2$ we
are lead to the rescaled system  
\begin{subequations}
  \label{eq:KolmoScaled}
\begin{align}
  \label{eq:KolmoScal.a}
\pl_\tau\wh v - \frac\gamma2 \,\wh v - \delta\,y \cdot \nabla \wh v&= \frac{\mu_1}2\DIV 
   (\wh k \nabla \wh v),
\\
  \label{eq:KolmoScal.b}
\pl_\tau \wh k - \gamma \,\wh k - \delta\,y \cdot \nabla \wh k &=  
 \mu_3\DIV(\wh k \nabla \wh k) + \frac{\mu_1}2 \wh k|\nabla \wh v|^2 - \alpha\, \wh k. 
\end{align}
\end{subequations}
The case $(\gamma,\delta)=(2,0)$ can be applied in bounded domains $\Omega$ and
it corresponds to the the case in the Prandtl equation. Clearly, we are in the
case $\alpha=\beta=1$ with two additional linear terms, namely
$ \wt v$ and $(2{-}a)\wt k$.  Thus, the existence theory of very weak solutions
in Theorem \eqref{th:Exi.eta=kappa} only applies in the case $\mu_1=2\mu_2$.

For the case $\delta>0$ we consider $\Omega=\R^d$, and the rescaled momentum
$V(\wh v,\wh k)=\int_{\R^d} v(y) \dd y $ and the 
 rescaled energy $E(\wh v, \wh k)=\int_{\R^d} \big(\frac12 \wh v^2 {+} \wh
k\big) \dd y$ now satisfies, along solutions $(\wh v(\tau),\wh k(\tau))$, the
linear relations 
\begin{align}
  \frac{\rmd}{\rmd \tau} V(\wh v,\wh k)= \big(\frac\gamma2- d\delta) \big) V(\wh
v,\wh k) \quad \text{and} \quad 
\frac{\rmd}{\rmd \tau} E(\wh v, \wh k) = (\gamma{-}d\delta) E(\wh v, \wh k ) - a
\int_{\R^d}  \wh k \dd y.  
\end{align}
Because of $2=\gamma+2\delta$ the case $(\gamma,\delta)
=\big(\frac{2d}{1{+}d},\ \frac1{1{+}d}\big)$ is special as it leads to 
\[
\tbinom\gamma\delta=\tfrac1{1{+}d} \tbinom{2d}1
\quad
\Longrightarrow \quad
\frac{\rmd}{\rmd \tau} V(\wh v,\wh k)= 0 \ \text{ and } \ 
\frac{\rmd}{\rmd \tau} E(\wh v, \wh k) = \frac{d}{1{+}d} E(\wh v, \wh k ) - a
\int_{\R^d}  \wh k \dd y.  
\]
For $a< \frac d{1{+}d}$ we obtain $\frac{\rmd}{\rmd \tau} E(\wh v, \wh k) >0$
as long as $E(\wh v, \wh k)>0$, so no steady state can exist. However, for $a>
\frac d{1{+}d}$ one can expect the existence of a family of state states 
$(v_\st,k_\st) $ such that $V(v_\st,k_\st)=V_0 \in \R$.

\footnotesize
%\bibliographystyle{alpha_AMs}
%\bibliography{alex_pub,bib_alex}

\newcommand{\etalchar}[1]{$^{#1}$}
\def\cprime{$'$}
\providecommand{\bysame}{\leavevmode\hbox to3em{\hrulefill}\thinspace}

\end{document}